\font\tenmib=cmmib10
\font\sevenmib=cmmib10 scaled 800
\font\titolo=cmbx12
\font\titolone=cmbx10 scaled\magstep 2

\font\cs=cmcsc10
\font\sc=cmcsc10
\font\css=cmcsc8

\font\ninerm=cmr9
\font\ottorm=cmr8
\textfont5=\tenmib\scriptfont5=\sevenmib\scriptscriptfont5=\fivei

\font\euftw=eufm9 scaled\magstep1
\font\euftww=eufm7 scaled\magstep1
\font\msytw=msbm9 scaled\magstep1

\font\msytwww=msbm7 scaled\magstep1

\font\indbf=cmbx10 scaled\magstep2

\font\ottorm=cmr8\font\ottoi=cmmi8\font\ottosy=cmsy8%
\font\ottobf=cmbx8\font\ottott=cmtt8%
\font\ottocss=cmcsc8%
\font\ottosl=cmsl8\font\ottoit=cmti8%
\font\sixrm=cmr6\font\sixbf=cmbx6\font\sixi=cmmi6\font\sixsy=cmsy6%
\font\fiverm=cmr5\font\fivesy=cmsy5\font\fivei=cmmi5\font\fivebf=cmbx5%

\def\ottopunti{\def\rm{\fam0\ottorm}%
\textfont0=\ottorm\scriptfont0=\sixrm\scriptscriptfont0=\fiverm%
\textfont1=\ottoi\scriptfont1=\sixi\scriptscriptfont1=\fivei%
\textfont2=\ottosy\scriptfont2=\sixsy\scriptscriptfont2=\fivesy%
\textfont3=\tenex\scriptfont3=\tenex\scriptscriptfont3=\tenex%
\textfont4=\ottocss\scriptfont4=\sc\scriptscriptfont4=\sc%
\textfont\itfam=\ottoit\def\it{\fam\itfam\ottoit}%
\textfont\slfam=\ottosl\def\sl{\fam\slfam\ottosl}%
\textfont\ttfam=\ottott\def\tt{\fam\ttfam\ottott}%
\textfont\bffam=\ottobf\scriptfont\bffam=\sixbf%
\scriptscriptfont\bffam=\fivebf\def\bf{\fam\bffam\ottobf}%
\setbox\strutbox=\hbox{\vrule height7pt depth2pt width0pt}%
\normalbaselineskip=9pt\let\sc=\sixrm\normalbaselines\rm}
\let\nota=\ottopunti%

%
%
%
%
%
%
%

\global\newcount\numsec\global\newcount\numapp
\global\newcount\numfor\global\newcount\numfig
\global\newcount\numsub
\numsec=0\numapp=0\numfig=0
\def\veroparagrafo{\number\numsec}\def\veraformula{\number\numfor}
\def\veraappendice{\number\numapp}\def\verasub{\number\numsub}
\def\verafigura{\number\numfig}

\def\section(#1,#2){\advance\numsec by 1\numfor=1\numsub=1\numfig=1%
\SIA p,#1,{\veroparagrafo} %
\write15{\string\Fp (#1){\secc(#1)}}%
\write16{ sec. #1 ==> \secc(#1)  }%
\hbox to \hsize{\titolo\hfill \number\numsec. #2\hfill%
\expandafter{\alato(sec. #1)}}\*}

\def\appendix(#1,#2){\advance\numapp by 1\numfor=1\numsub=1\numfig=1%
\SIA p,#1,{A\veraappendice} %
\write15{\string\Fp (#1){\secc(#1)}}%
\write16{ app. #1 ==> \secc(#1)  }%
\hbox to \hsize{\titolo\hfill Appendix A\number\numapp. #2\hfill%
\expandafter{\alato(app. #1)}}\*}

\def\senondefinito#1{\expandafter\ifx\csname#1\endcsname\relax}

\def\SIA #1,#2,#3 {\senondefinito{#1#2}%
\expandafter\xdef\csname #1#2\endcsname{#3}\else
\write16{???? ma #1#2 e' gia' stato definito !!!!} \fi}

\def \Fe(#1)#2{\SIA fe,#1,#2 }
\def \Fp(#1)#2{\SIA fp,#1,#2 }
\def \Fg(#1)#2{\SIA fg,#1,#2 }

\def\etichetta(#1){(\veroparagrafo.\veraformula)%
\SIA e,#1,(\veroparagrafo.\veraformula) %
\global\advance\numfor by 1%
\write15{\string\Fe (#1){\equ(#1)}}%
\write16{ EQ #1 ==> \equ(#1)  }}

\def\etichettaa(#1){(A\veraappendice.\veraformula)%
\SIA e,#1,(A\veraappendice.\veraformula) %
\global\advance\numfor by 1%
\write15{\string\Fe (#1){\equ(#1)}}%
\write16{ EQ #1 ==> \equ(#1) }}

\def\getichetta(#1){\veroparagrafo.\verafigura%
\SIA g,#1,{\veroparagrafo.\verafigura} %
\global\advance\numfig by 1%
\write15{\string\Fg (#1){\graf(#1)}}%
\write16{ Fig. #1 ==> \graf(#1) }}

\def\getichettaa(#1){A\veraappendice.\verafigura%
\SIA g,#1,{A\veraappendice.\verafigura} %
\global\advance\numfig by 1%
\write15{\string\Fg (#1){\graf(#1)}}%
\write16{ Fig. #1 ==> \graf(#1) }}

\def\etichettap(#1){\veroparagrafo.\verasub%
\SIA p,#1,{\veroparagrafo.\verasub} %
\global\advance\numsub by 1%
\write15{\string\Fp (#1){\secc(#1)}}%
\write16{ par #1 ==> \secc(#1)  }}

\def\etichettapa(#1){A\veraappendice.\verasub%
\SIA p,#1,{A\veraappendice.\verasub} %
\global\advance\numsub by 1%
\write15{\string\Fp (#1){\secc(#1)}}%
\write16{ par #1 ==> \secc(#1)  }}

\def\Eq(#1){\eqno{\etichetta(#1)\alato(#1)}}
\def\eq(#1){\etichetta(#1)\alato(#1)}
\def\Eqa(#1){\eqno{\etichettaa(#1)\alato(#1)}}
\def\eqa(#1){\etichettaa(#1)\alato(#1)}
\def\eqg(#1){\getichetta(#1)\alato(fig. #1)}
\def\sub(#1){\0\palato(p. #1){\bf Definition \etichettap(#1).}}
\def\asub(#1){\0\palato(p. #1){\bf Definition \etichettapa(#1).}}
\def\lm(#1){\0\palato(p. #1){\bf Lemma \etichettap(#1).}}
\def\lma(#1){\0\palato(p. #1){\bf Lemma \etichettapa(#1).}}
\def\equv(#1){\senondefinito{fe#1}$\clubsuit$#1%
\write16{eq. #1 non e' (ancora) definita}%
\else\csname fe#1\endcsname\fi}
\def\grafv(#1){\senondefinito{fg#1}$\clubsuit$#1%
\write16{fig. #1 non e' (ancora) definito}%
\else\csname fg#1\endcsname\fi}
\def\secv(#1){\senondefinito{fp#1}$\clubsuit$#1%
\write16{par. #1 non e' (ancora) definito}%
\else\csname fp#1\endcsname\fi}

\def\equ(#1){\senondefinito{e#1}\equv(#1)\else\csname e#1\endcsname\fi}
\def\graf(#1){\senondefinito{g#1}\grafv(#1)\else\csname g#1\endcsname\fi}
\def\figura(#1){{\css Figure} \getichetta(#1)}
\def\figuraa(#1){{\css Figure} \getichettaa(#1)}
\def\secc(#1){\senondefinito{p#1}\secv(#1)\else\csname p#1\endcsname\fi}
\def\sec(#1){{\S\secc(#1)}}
\def\refe(#1){{[\secc(#1)]}}

\def\BOZZA{\bz=1
\def\alato(##1){\rlap{\kern-\hsize\kern-1.2truecm{$\scriptstyle##1$}}}
\def\palato(##1){\rlap{\kern-1.2truecm{$\scriptstyle##1$}}}
}

\def\alato(#1){}
\def\galato(#1){}
\def\palato(#1){}


{\count255=\time\divide\count255 by 60 \xdef\hourmin{\number\count255}
        \multiply\count255 by-60\advance\count255 by\time
   \xdef\hourmin{\hourmin:\ifnum\count255<10 0\fi\the\count255}}

\def\oramin{\hourmin }

\def\data{\number\day/\ifcase\month\or gennaio \or febbraio \or marzo \or
aprile \or maggio \or giugno \or luglio \or agosto \or settembre
\or ottobre \or novembre \or dicembre \fi/\number\year;\ \oramin}
\footline={\rlap{\hbox{\copy200}}\tenrm\hss \number\pageno\hss}


\newcount\driver 
\newdimen\xshift \newdimen\xwidth
\def\ins#1#2#3{\vbox to0pt{\kern-#2 \hbox{\kern#1
#3}\vss}\nointerlineskip}

\def\insertplot#1#2#3#4#5{\par%
\xwidth=#1 \xshift=\hsize \advance\xshift by-\xwidth \divide\xshift by 2%
\yshift=#2 \divide\yshift by 2%
\line{\hskip\xshift \vbox to #2{\vfil%
\ifnum\driver=0 #3
\special{ps: plotfile #4.eps} 
\fi
\ifnum\driver=1 #3 \includegraphics{#4.eps}\fi
\ifnum\driver=2 #3
\ifnum\mgnf=0\special{#4.eps 1. 1. scale} \fi
\ifnum\mgnf=1\special{#4.eps 1.2 1.2 scale}\fi
\fi }\hfill \raise\yshift\hbox{#5}}}

\def\insertplotttt#1#2#3{\par%
\xwidth=#1 \xshift=\hsize \advance\xshift by-\xwidth \divide\xshift by 2%
\yshift=#2 \divide\yshift by 2%
\line{\hskip\xshift \vbox to #2{\vfil%
\includegraphics{#3.eps}}\hfill}}

\def\insertplott#1#2#3{\par%
\xwidth=#1 \xshift=\hsize \advance\xshift by-\xwidth \divide\xshift by 2%
\yshift=#2 \divide\yshift by 2%
\line{\hskip\xshift \vbox to #2{\vfil%
\ifnum\driver=0
\special{ps: plotfile #3.eps} 
 \fi
\ifnum\driver=1 \includegraphics{#3.eps}\fi
\ifnum\driver=2
\ifnum\mgnf=0\special{#3.eps 1. 1. scale} \fi 
\ifnum\mgnf=1\special{#3.eps 1.2 1.2 scale}\fi
\fi }\hfill}}

\newdimen\xshift \newdimen\xwidth \newdimen\yshift
\def\eqfig#1#2#3#4#5{
\par\xwidth=#1 \xshift=\hsize \advance\xshift
by-\xwidth \divide\xshift by 2
\yshift=#2 \divide\yshift by 2
\line{\hglue\xshift \vbox to #2{\vfil
\ifnum\driver=0 #3
\special{ps: plotfile #4.ps} 
\fi
\ifnum\driver=1 #3 \includegraphics{#4.ps}\fi
\ifnum\driver=2 #3 \special{
\ifnum\mgnf=0 #4.ps 1. 1. scale \fi
\ifnum\mgnf=1 #4.ps 1.2 1.2 scale\fi}
\fi}\hfill\raise\yshift\hbox{#5}}}

\let\a=\alpha \let\b=\beta  \let\g=\gamma  \let\d=\delta \let\e=\varepsilon
  \let\h=\eta   \let\th=\theta  \let\l=\lambda
\let\m=\mu    \let\n=\nu         \let\p=\pi    \let\r=\rho
\let\s=\sigma \let\t=\tau    
   \let\o=\omega
\let\G=\Gamma \let\D=\Delta  \let\Th=\Theta\let\L=\Lambda 
\let\P=\Pi          
\let\O=\Omega 
\def\R{{\cal{R}}}
\def\vb{\vec {\bf b}}
\def\ka{{(k)}}
\def\\{\hfill\break} \let\==\equiv

\let\io=\infty 

\let\0=\noindent

\let\dpr=\partial
\def\der{{\rm d}}

\def\tende#1{\,\vtop{\ialign{##\crcr\rightarrowfill\crcr
 \noalign{\kern-1pt\nointerlineskip}
 \hskip3.pt${\scriptstyle #1}$\hskip3.pt\crcr}}\,}
\def\otto{\,{\kern-1.truept\leftarrow\kern-5.truept\to\kern-1.truept}\,}

\def\PPP{{\cal P}}\def\MMM{{\cal M}} \def\VV{{\cal V}}
\def\CC{{\cal C}}\def\WW{{\cal W}}
\def\TT{{\cal T}}\def\BBB{{\cal B}}
\def\RR{{\cal R}} \def\OO{{\cal O}}
\def\DDD{{\cal D}}\def\AAA{{\cal A}}\def\GG{{\cal G}}

\def\MM{\hbox{\euftw{M}}} \def\EE{\hbox{\euftw{E}}}
\def\DD{\hbox{\euftw{D}}} \def\CCCCCC{\hbox{\euftw{A}}}
\def\SS{\hbox{\euftw{S}}} \def\III{\hbox{\euftw{I}}}
 \def\KKK{\hbox{\euftw{K}}}
\def\EEsub{\hbox{\euftww{E}}}

\def\T#1{{#1_{\kern-3pt\lower7pt\hbox{$\widetilde{}$}}\kern3pt}}
\def\VVV#1{{\underline #1}_{\kern-3pt
\lower7pt\hbox{$\widetilde{}$}}\kern3pt\,}
\def\W#1{#1_{\kern-3pt\lower7.5pt\hbox{$\widetilde{}$}}\kern2pt\,}

\def\indica{\leaders \hbox to 0.5cm{\hss.\hss}\hfill}
\def\guida{\leaders\hbox to 1em{\hss.\hss}\hfill}

   \def\V0{{\bf 0}}

  \def\bb{{\bf b}}
\def\ul{\underline}

\def\oo{{\omega}}

\mathchardef\aa   = "050B
\mathchardef\bbb  = "050C
\mathchardef\xxx  = "0518
\mathchardef\hhh  = "0511
\mathchardef\zzzzz= "0510
\mathchardef\oo   = "0521
\mathchardef\tt   = "051C
\mathchardef\lll  = "0515
\mathchardef\mmm  = "0516
\mathchardef\kkkk = "056B
\mathchardef\Dp   = "0540
\mathchardef\H    = "0548
\mathchardef\FFF  = "0546
\mathchardef\ppp  = "0570
\mathchardef\nnnnn= "056E
\mathchardef\jjjj = "056A
\mathchardef\pps  = "0520
\mathchardef\XXX  = "0504
\mathchardef\FFF  = "0508
\mathchardef\LLLL = "054C
\mathchardef\PPPP = "0550
\mathchardef\JJJJ = "054A
\mathchardef\QQQQ = "0551
\mathchardef\DDDDD= "0544
\mathchardef\TTTT = "0502
\mathchardef\RRRR = "0552
\mathchardef\UUUU = "0555
\mathchardef\minore  = "053C
\mathchardef\maggiore  = "053E

\def\to{\rightarrow}
\def\la{\left\langle}
\def\ra{\right\rangle}

\def\qed{\hfill\raise1pt\hbox{\vrule height5pt width5pt depth0pt}}

 \def\Val{{\rm Val}}
\def\indic{\hbox{\raise-2pt \hbox{\indbf 1}}}

\def\DDDD{\hbox{\msytw D}}

\def\RRR{\hbox{\msytw R}} 

\def\CCC{\hbox{\msytw C}} 

\def\NNN{\hbox{\msytw N}} 
\def\nnn{\hbox{\msytwww N}}
\def\ZZZ{\hbox{\msytw Z}} 
\def\zzz{\hbox{\msytwww Z}}

\def\GGG{\hbox{\msytw G}}


\newcount\mgnf  
\mgnf=0 

\ifnum\mgnf=0
\def\openone{\leavevmode\hbox{\ninerm 1\kern-3.3pt\tenrm1}}%
\def\*{\vglue0.3truecm}\fi
\ifnum\mgnf=1
\def\openone{\leavevmode\hbox{\ninerm 1\kern-3.63pt\tenrm1}}%
\def\*{\vglue0.5truecm}\fi


\newcount\tipobib\newcount\bz\bz=0\newcount\aux\aux=1
\newdimen\bibskip\newdimen\maxit\maxit=0pt


\tipobib=0
\def\9#1{\ifnum\aux=1#1\else\relax\fi}

\newwrite\bib
\immediate\openout\bib=\jobname.bib
\global\newcount\bibex
\bibex=0
\def\verabib{\number\bibex}

\ifnum\tipobib=0
\def\cita#1{\expandafter\ifx\csname c#1\endcsname\relax
\hbox{$\clubsuit$}#1\write16{Manca #1 !}%
\else\csname c#1\endcsname\fi}
\def\rife#1#2#3{\immediate\write\bib{\string\raf{#2}{#3}{#1}}
\immediate\write15{\string\C(#1){[#2]}}
\setbox199=\hbox{#2}\ifnum\maxit < \wd199 \maxit=\wd199\fi}
\else
\def\cita#1{%
\expandafter\ifx\csname d#1\endcsname\relax%
\expandafter\ifx\csname c#1\endcsname\relax%
\hbox{$\clubsuit$}#1\write16{Manca #1 !}%
\else\probib(ref. numero )(#1)%
\csname c#1\endcsname%
\fi\else\csname d#1\endcsname\fi}%
\def\rife#1#2#3{\immediate\write15{\string\Cp(#1){%
\string\immediate\string\write\string\bib{\string\string\string\raf%
{\string\verabib}{#3}{#1}}%
\string\Cn(#1){[\string\verabib]}%
\string\CCc(#1)%
}}}%
\fi

\def\Cn(#1)#2{\expandafter\xdef\csname d#1\endcsname{#2}}
\def\CCc(#1){\csname d#1\endcsname}
\def\probib(#1)(#2){\global\advance\bibex+1%
\9{\immediate\write16{#1\verabib => #2}}%
}

\def\C(#1)#2{\SIA c,#1,{#2}}
\def\Cp(#1)#2{\SIAnx c,#1,{#2}}

\def\SIAnx #1,#2,#3 {\senondefinito{#1#2}%
\expandafter\def\csname#1#2\endcsname{#3}\else%
\write16{???? ma #1,#2 e' gia' stato definito !!!!}\fi}

\bibskip=10truept
\def\hboxto{\hbox to}

\catcode`\{=12\catcode`\}=12
\catcode`\<=1\catcode`\>=2
\immediate\write\bib<
        \string\halign{\string\hboxto \string\maxit%
        {\string #\string\hfill}&%
        \string\vtop{\string\parindent=0pt\string\advance\string\hsize%
        by -.55truecm%
        \string#\string\vskip \bibskip
        }\string\cr%
>
\catcode`\{=1\catcode`\}=2
\catcode`\<=12\catcode`\>=12

\def\raf#1#2#3{\ifnum \bz=0 [#1]&#2 \cr\else
\llap{${}_{\rm #3}$}[#1]&#2\cr\fi}

\newread\bibin

\catcode`\{=12\catcode`\}=12
\catcode`\<=1\catcode`\>=2
\def\chiudibib<
\catcode`\{=12\catcode`\}=12
\catcode`\<=1\catcode`\>=2
\immediate\write\bib<}>
\catcode`\{=1\catcode`\}=2
\catcode`\<=12\catcode`\>=12
>
\catcode`\{=1\catcode`\}=2
\catcode`\<=12\catcode`\>=12

\def\makebiblio{
\ifnum\tipobib=0
\advance \maxit by 10pt
\else
\maxit=1.truecm
\fi
\chiudibib
\immediate \closeout\bib
\openin\bibin=\jobname.bib
\ifeof\bibin\relax\else
\raggedbottom
\input \jobname.bib
\fi
}

\openin13=#1.aux \ifeof13 \relax \else
\input #1.aux \closein13\fi
\openin14=\jobname.aux \ifeof14 \relax \else
\input \jobname.aux \closein14 \fi
\immediate\openout15=\jobname.aux

\def\biblio{\*\*\centerline{\titolo References}\*\nobreak\makebiblio}


\ifnum\mgnf=0
   \magnification=\magstep0
   \hsize=16.truecm\vsize=21.6truecm\voffset0.0truecm\hoffset-.2truecm
   \parindent=0.3cm\baselineskip=0.43cm\fi
\ifnum\mgnf=1
   \magnification=\magstep1\hoffset=0.truecm
   \hsize=16.truecm\vsize=21.8truecm
   \baselineskip=18truept plus0.1pt minus0.1pt \parindent=0.9truecm
   \lineskip=0.5truecm\lineskiplimit=0.1pt      \parskip=0.1pt plus1pt\fi


\mgnf=0   
\driver=1 
\openin14=\jobname.aux \ifeof14 \relax \else
\input \jobname.aux \closein14 \fi
\openout15=\jobname.aux




\centerline{\titolone Periodic solutions for the
Schr\"odinger equation with}
\vskip.1cm
\centerline{\titolone nonlocal smoothing nonlinearities
in higher dimension}

\vskip1.truecm \centerline{{\titolo
Guido Gentile and Michela Procesi}}
\vskip.2truecm
\centerline{{} Dipartimento di Matematica,
Universit\`a di Roma Tre, Roma, I-00146}

\vskip1.truecm

\line{\vtop{
\line{\hskip1.1truecm\vbox{\advance \hsize by -2.3 truecm
\0{\cs Abstract.}
{\it We consider the nonlinear Schr\"odinger equation in higher
dimension with Dirichlet boundary conditions and with a non-local
smoothing nonlinearity. We prove the existence of small amplitude
periodic solutions. In the fully resonant case we find solutions
which at leading order are wave packets, in the sense that
they continue linear solutions with an arbitrarily large number
of resonant modes. The main difficulty in the proof
consists in solving a ``small divisor problem'' which we do
by using a renormalisation group approach.  } \hfill} }}}

\*\*
\section(1,Introduction and results)

\0In this paper we prove the existence of small amplitude
periodic solutions for a class of nonlinear Schr\"odinger equations
in $D$ dimensions
$$ iv_t-\Delta v +\m v = f(x,\Phi(v),\Phi(\bar v)) :=
|\Phi(v)|^{2}\Phi(v)+F(x,\Phi(v),\Phi(\bar v)) ,
\Eq(1.1) $$
with Dirichlet boundary conditions on the square $[0,\pi]^{D}$.
Here $D\ge 2$ is an integer, $\m$ is a real parameter,
$\Phi$ is a smoothing operator, which in Fourier space acts as
$$ (\Phi(u))_k = |k|^{-2s} u_k ,
\Eq(1.2) $$
for some positive $s$, and $F$ is an analytic odd function,
real for real $u$, such that $F(x,u,\bar u)$ is of order higher
than three in $(u,\bar u)$, i.e.
$$ F(x,u,\bar u) =
\sum_{p=4}^{\io} \sum_{p_{1}+p_{2}=p}
a_{p_{1},p_{2}}(x) u^{p_{1}} \bar u^{p_{2}} ,
\qquad F(-x,-u,-\bar u) = - F(x,u,\bar u) .
\Eq(1.3) $$
In particular this implies that the functions $a_{p_{1},p_{2}}$
must be even for odd $p$ and odd for even $p$, and real for all $p$.
The reality condition is assumed to simplify the analysis.

For $D=2$ we do not impose any further condition on $f$,
whereas for $D\ge 3$ we shall consider a more restrictive class
of nonlinearities, by requiring
$$ f(x,u,\bar u) = {\dpr\over \dpr \bar u} H(x,u, \bar u)
+ g(x,\bar u) , \qquad \overline {H(x,u,\bar u)} = H(x,u,\bar u) ,
\Eq(1.4) $$
i.e. with $H$ a real function and $g$ depending explicitly only
on $\bar u$ (besides $x$) and not on $u$.

In general when looking for small periodic solutions for PDE's one
expects to find a ``small divisor problem'' due to the fact that the
eigenvalues of the linear term   accumulate to zero in the space of
$T-$periodic solutions, for any $T$ in a positive measure set.

The case of one space dimension was widely studied in the '90 for
non-resonant equations by using KAM theory by Kuksin-P\"oshel \cita{K1},
\cita{KP} and Wayne \cita{W}, and by using Lyapunov-Schmidt decomposition
by Craig-Wayne \cita{CW} and Bourgain \cita{Bo1}, \cita{Bo4}.
The two techniques are somehow complementary.
The Lyapunov-Schmidt decomposition is more flexible:
it can be successfully adapted to non-Hamiltonian equations
and to ``resonant''equations, i.e. where the linear frequencies are
not rationally independent \cita{LS}, \cita{Bo3}, \cita{GMP}.
On the other hand KAM theory provides more information, for instance
on the stability of the solutions.

Generally speaking  the main feature which is used to solve
the small divisor problem (in all the above mentioned techniques) 
is the ``separation of the resonant sites''.
Such a feature can be described as follows. For instance
for $D=1$ consider an equation $\DDDD[u]=f(u)$, where $\DDDD$
is a linear differential operator and $f(u)$
a smooth super-linear function; let $\l_{k}$ with $k\in \ZZZ^{2}$
be the linear eigenvalues in the space of $T$-periodic solutions,
so that after rescaling the amplitude and in Fourier
space the equation has the form
$$ \l_{k} u_{k}= \e f_{k}(u) ,
\Eq(1.5) $$
with $\inf_{k} |\l_{k}|=0$. The separation property
for Dirichlet boundary conditions requires:

\*

\01. if $|\l_{k}|<\a$ then $|k|> C \a^{-\d_{0}}$
(this is generally obtained by restricting $T$ to a Cantor set).

\02. if both  $|\l_{k}|<\a$ and $|\l_{h}|<\a$ then either $h=k$ or
$|h-k|\ge C (\min \{|h|,|k|\})^{\d}$.

\*

Here $\d_{0}$ and $\d$ are model-dependent parameters, and $C$
is some positive constant. In the case
of periodic boundary conditions, 2. should be suitably modified.

It is immediately clear that 2. cannot be satisfied by our equation
\equ(1.1) as the linear eigenvalues are
$$ \l_{n,m} = - \o n + |m|^{2} + \mu , \qquad
\o={2\pi\over T} ,
\Eq(1.6) $$
so that all the eigenvalues $\l_{n_{1},m_{1}}$
with $n_{1}=n$ and $|m_{1}|=|m|$ are equal to $\l_{n,m}$.

The existence of periodic solutions for $D>1$ space dimensions
was first proved by Bourgain in \cita{Bo2} and \cita{Bo4},
by using a Lyapunov-Schmidt decomposition and a technique by
Spencer and Fr\"olich to solve the small divisor problem. 
Again the separation properties are crucial: 1. is assumed and 2. is
weakened in the following way:

\*

\02$'$. the sets of $k\in \ZZZ^{D+1}$ such that $|\l_{k}|<1$ and
$R<|k|<2R$ are separated in clusters, say $C_{j}$ with $j\in \NNN$,
such that each cluster contains at most $R^{\d_{1}}$ elements and
dist$(C_{i},C_{j})\ge R^{\d_2}$, with $0<\d_2\leq \d_{1}\ll 1$.

\*

Now, in order to apply Spencer and Fr\"olich's method, one has to
control the eigenvalues of appropriate matrices
of dimension comparable to $|C_{j}|$. Such dimension goes to infinity
with $R$ and at the same time the linear eigenvalues go to zero,
so that achieving such estimates is a rather delicate question.

\*

Recently Bourgain also proved the existence of quasi-periodic
solutions for the nonlinear Schr\"odinger equation, with
local nonlinearities, in any dimensions \cita{Bo5}.
Still more recently in \cita{EK}, Eliasson and Kuksin
proved the same result by using KAM techniques.
We can also mention a very recent preprint by Yuan \cita{Y},
where a variant of the KAM approach was provided to show
the existence of quasi-periodic solutions: in this version,
stability of the solutions is not obtained, but, conversely,
the proof rather simplifies with respect to that given in \cita{EK}.

In this paper we use a Lyapunov-Schmidt decomposition
and then the so-called ``Lindstedt series method''
\cita{GM2} to solve the small divisor problem.  The main purpose
of this paper is to reobtain Bourgain's result \cita{Bo2} with
the Lindstedt series method, on the simplest possible model
which still carries the main difficulties of the $D$ space dimensions.
Recently Geng and You \cita{GY} have proved, via KAM theory,
the existence of quasi-periodic solutions for the NLS with
a non-local smoothing non-linearity and with periodic boundary
conditions; in such case they show the existence of a symmetry,
which greatly simplifies the analysis. In the case of Dirichlet
boundary condition this symmetry is broken, so that the results
of \cita{GY} do not apply to the equation \equ(1.1) with
Dirichlet boundary conditions. None the less the
regularisation provides some nice simplifications. This motivates
our choice of equation \equ(1.1), since the main purpose of the paper
is to establish appropriate techniques and notation in the
simplest (non-trivial) possible case. 

Moreover, we are able to find periodic solutions also in some
non-Hamiltonian and in resonant cases, where 
the result was not known in the literature.
In particular in the completely resonant case ($\m=0$ in \equ(1.1))
we find solutions which reduce to wave packets
(i.e. linear combinations of harmonics centred around
suitable frequencies) in the absence of the perturbation.

\*

Let us now describe the general lines of the Lindstedt series
approach, which were originally developed by Eliasson \cita{E}
and Gallavotti \cita{Ga} in the context
of KAM theory for finite dimensional systems.

The main idea is to consider a ``renormalisation'' of equation
\equ(1.5) which can be proved to have solutions. More precisely we
consider a new, vector-valued, equation  with unknowns
$U_{j}:= \{u_{k}: k\in C_{j}\}$
$$ \left( \DDDD_{j} (\omega) +M_{j} \right) U_{j}=
\e F_{j}(U)+ L_{j}U_{j}  ,
\Eq(1.7) $$
where $\DDDD_{j}(\omega)$ is the diagonal matrix of the eigenvalues 
$\l_{k}$ with $k\in C_{j}$, $ F_{j}(U)$ is the vector $\{f_{k}(u):
k\in C_{j}\}$ defined in \equ(1.5)  and  $M_{j},L_{j}$
are matrices of free parameters. Equation \equ(1.7) coincides
with \equ(1.5) provided $M_{j}=L_{j}$ for all $j\in \NNN$.

The aim then is to proceed as in the one dimensional
renormalisation scheme proposed in \cita{GM2} and \cita{GMP};
namely we restrict $(\o,\{M_{j}\})$ to a Cantor set and construct
both the solution $U_{j}(\e,\o,\{M_{h}\})$ and $L_{j}(\e,\o,\{M_{h}\})$
as convergent power series in $\e$.
Then one solves the compatibility equation
$M_{j}=L_{j}(\e,\omega,\{M_{h}\})$; essentially  this is done  by the
implicit function theorem but with the additional complication  that
$L_{j}$  is defined for $(\omega,\{M_{h}\})$ in a Cantor set.

\*

We look for periodic solutions of frequency $\o= D +\m-\e$,
with $\e>0$, which continue the unperturbed one ($\e=0$)
with frequency $\o_{0}= D+\m$.
Note that the choice of this particular unperturbed  frequency is
made only for the sake of definiteness: any other linear frequency
would yield the same type of results.

For $\e\neq0$ we perform the change of variables
$$ \sqrt{\e} u(x,t)= \Phi(v(x,\o t)) ,
\Eq(1.8) $$
so that \equ(1.1) becomes
$$ \Phi^{-1}(i\o u_t-\D u+ \m u)= \e |u|^{2}u +
{1\over\sqrt{\e}} F(x,\sqrt{\e}u,\sqrt{\e}\bar u) \=
\e f(x,u,\bar u,\e) ,
\Eq(1.9) $$
with a slight abuse of notation in the definition of $f$.

We start by considering explicitly the case $F=0$, for
simplicity, so that $f(x,u,\bar u,\e)=f(u,\bar u)=|u|^{2}u$.
In that case the problem of the existence of periodic solutions
becomes trivial, but the advantage of proceeding this way is that
the construction that we are going to envisage extends
easily to more general $f$, with some minor technical adaptations.

We pass to the equation for the Fourier coefficients, by writing
$$ u(x,t)= \sum_{n\in \zzz, m\in \zzz^{D}} u_{n,m}
{\rm e}^{i (n  t+ m\cdot x)} ,
\Eq(1.10) $$
so that \equ(1.9) gives
$$ |m|^{2s} \left( -\o n+ |m|^{2}+\m \right) u_{n,m} =
\e \sum_{n_{1}+n_2-n_3=n\atop m_{1}+m_2-m_3= m}
u_{n_{1},m_{1}}u_{n_2,m_2} \bar u_{n_3,m_3}
\equiv \e f_{n,m}(u,\bar u) ,
\Eq(1.11) $$
and the Dirichlet boundary conditions spell
$$ u_{n,m}= u_{n,S_{i}(m)}\,,\qquad S_{i}(e_{j})=\left( 1 - 2 \d(i,j)
\right) e_{j} \qquad \forall i=1,\ldots, D ,
\Eq(1.12) $$
where $\d(i,j)$ is Kronecker's delta and $S_{i}(m)$ is the linear
operator that changes the sign of the $i$-th component of $m$.

\*

We proceed as follows. We perform a Lyapunov-Schmidt decomposition
separating the $P$-$Q$ supplementary subspaces. By definition
$Q$ is the space of Fourier labels $(n,m)$ such that
$u_{n,m}$ solves \equ(1.11) at $\e=0$. If $\m\neq0$ we impose
an irrationality condition on $\m$, i.e. $\o_{0}n-p \neq 0$,
so that $Q$ is defined as
$$ Q:= \left\{ (n,m) \in \ZZZ \times \ZZZ^{D} :
n=1, m_{i}=\pm 1\,\forall i \right\} .
\Eq(1.13) $$
By the Dirichlet boundary conditions, calling $V=\{1,1,\dots,1\}$, for
all $(1,m)\in Q$ we have that $u_{1,m}= \pm u_{1,V}$; see \equ(1.12).
Then \equ(1.11) naturally splits into two sets of equations:
the $Q$ equations, for $(n,m)$ such that $n=1$ and $|m|=\sqrt{D}$,
and the $P$ equations, for all the other values of $(n,m)$.
We first solve the $P$ equation keeping  $q:=u_{1,V}$ as a
parameter. Then we consider the $Q$ equations and solve them
via the implicit function theorem.

\*   

We look for solutions of \equ(1.11) such that $u_{n,m}\in \RRR$ for all
$(n,m)$; this is possible as one can find real solutions for the
bifurcation equations in $Q$, and then the recursive $P$-$Q$ equations
are closed on the subspace of real $u_{n,m}$. The same condition
can be imposed also in the more general case \equ(1.3), provided
the functions $a_{p_{1},p_{2}}$ are real, as we are assuming.

For $\m \neq 0$ we shall construct periodic solutions which are
analytic both in time and space, and not only sub-analytic,
as usually found \cita{Bo2}. 
This is due to the presence of the smoothing non-linearity.

\*

\0{\bf Theorem 1.}
{\it Consider equation \equ(1.9), with $\Phi$ defined by \equ(1.2)
for arbitrary $s>0$ and $F$ given by \equ(1.3) if
$D=2$ and by \equ(1.3) and \equ(1.4) if $D\ge 3$. There exist a
Cantor set $\MM \subset (0,\mu_{0})$ and a constant $\e_{0}$
such that the following holds. For all $\m\in \MM$
there exists a Cantor set $\EE(\m) \subset (0,\e_{0})$, such that
for all $\e\in \EE(\m)$ the equation admits a solution $u(x,t)$,
which is $2\pi$-periodic in time,
analytic in time and in space, such that
$$ \left| u(x,t) - q_{0} \, {\rm e}^{it}
\prod_{i=1}^{D} \sin x_{i} \right|\le C \e, \qquad
q_{0} = \sqrt{ D^{s}3^{-D} } ,
\Eq(1.14) $$
uniformly in $(x,t)$. The set $\MM$ has full measure and for all
$\m\in\MM$ the set $\EE=\EE(\mu)$ has positive Lebesgue measure and
$$ \lim_{\e\to 0^{+}} {{\rm meas}(\EE \cap [0,\e]) \over \e}=1 ,
\Eq(1.15) $$
where ${\rm meas}$ denotes the Lebesgue measure.}

\*

For $\m=0$ the following result extends Theorem 1 of \cita{GP}
to the higher dimensional case.

\*

\0{\bf Theorem 2.}
{\it Consider equation \equ(1.9) with $\m=0$, $D\ge 2$,
$\Phi$ defined by \equ(1.2) and $F$ given \equ(1.3) and \equ(1.4).
There exist a constant $\e_{0}$ and a Cantor set $\EE
\subset (0,\e_{0})$, such that for all $\e\in \EE$ the equation
admits a solution $u(x,t)$, which is $2\pi$-periodic in time,
sub-analytic in time and in space, satisfying
\equ(1.14) and \equ(1.15).}

\*

\0{\bf Remark.} For $\mu\neq 0$ we could consider other unperturbed periodic 
solutions and we would obtain the same kind of results as in Theorem 1, 
with only some trivial changes of notation in the proofs. 

For $\m=0$ and if the functions
$a_{p_1,p_2}(x)$ in \equ(1.3) are constant, we could easily
extend Theorem 2 to other unperturbed solutions. Considering 
non-constant $a_{p_1,p_2}$'s would require some extra work.
\*

For $D=2$ the following result extends Theorem 2 of \cita{GP}.

\*

\0{\bf Theorem 3.}
{\it Consider equation \equ(1.9) with $\m=0$, $D=2$,
$\Phi$ defined by \equ(1.2) and $F$ given by \equ(1.3) and \equ(1.4).
Let $\KKK$ any interval in $\RRR_{+}$. For $N>4$ there exist
sets $\MMM_+$ of $N$ vectors in $\ZZZ_{+}^{2}$
and sets of amplitudes $a_m$ with $m\in \MMM_+$ such that
the following holds. Define
$$ q_{0}(x,t)= \sum_{m\in \MMM_+} a_{m} \, {\rm e}^{i|m|^{2} t}
\sin (m_{1} x_{1})\sin (m_{2} x_{2}) .
\Eq(1.16) $$
There are a finite set $\KKK_{0}$ of points in $\KKK$, a positive
constant $\e_{0}$ and a set $\EE\in (0,\e_{0})$ (all depending
on $\MMM_+$), such that for all $s\in\KKK\setminus\KKK_{0}$ and
$\e\in \EE$, equation \equ(1.9) admits a solution $u(x,t)$, which
is $2\pi$-periodic in time, sub-analytic in time and space, such that
$$ |u(x,t)-q_{0}(x,t)|\le C \e,
\Eq(1.17) $$
uniformly in $(x,t)$. Finally
$$\lim_{\e\to 0^{+}} {{\rm meas}(\EE\cap [0,\e]) \over \e}=1 ,
\Eq(1.18) $$
where ${\rm meas}$ denotes the Lebesgue measure.}

\*

In the case $D>2$ we can still find a solution of the leading order
of the $Q$ equations of the form \equ(1.16); however in order to
prove the existence  of a solution $u(x,t)$ of the full equation
we need a ``non-degeneracy condition'', namely that some finite
dimensional matrix (denoted by $J_{1,1}$ and defined in
Section \secc(8)) is invertible. 

\*
\0{\bf Theorem 4.}
{\it Consider equation \equ(1.9) with $\m=0$, $D\ge 2$,
$\Phi$ defined by \equ(1.2) and $F$ given \equ(1.3) and \equ(1.4).
There exist sets  $\MMM_+$ of $N$ vectors in $\ZZZ_{+}^{D}$
and sets of amplitudes $a_m$ with $m\in \MMM_+$ such that
the $Q$ equations at $\e=0$ have the solution
$$ q_{0}(x,t)= \sum_{m\in \MMM_+} a_{m} \, {\rm e}^{i|m|^{2} t}
\prod_{i=1}^{D} \sin (m_{i} x_{i}).
\Eq(1.19) $$ 
The set $\MMM_+$ identifies a finite order matrix $J_{1,1}$ (depending
analytically on the parameter $s$).
For $N>1$ if ${\rm det} J_{1,1}=0$
is not an identity in $s$ then the following holds. 
There are a finite set $\KKK_{0}$ of points in $\KKK$, a positive
constant $\e_{0}$ and a set $\EE\in (0,\e_{0})$ (all depending
on $\MMM_+$), such that for all $s\in\KKK\setminus\KKK_{0}$ and
$\e\in \EE$, equation \equ(1.9) admits a solution $u(x,t)$, which is
$2\pi$-periodic in time, sub-analytic in time and space, such that
$$ \left| u(x,t)-q_{0}(x,t) \right| \le C \e,
\Eq(1.20) $$
uniformly in $(x,t)$, and $\EE$ satisfies the property \equ(1.18).}

\*\*
\section(2, Technical set-up and propositions)

\0{\bf 2.1. Separation of the small divisors}\*

\0Let us require that $\m$ is strongly non-resonant (and in a full
measure set), i.e. that there exist $1 \gg \g_{0}>0$ and $\t_{0}>1$
such that
$$ \left| \left( D+\m \right) n - p - a \m \right| \geq
{\g_{0}\over |n|^{\t_{0}}} \quad \forall a=0,1, \quad
(n,p) \in \ZZZ^{2} ,\quad (n, p)\neq (1,D) , \quad n \neq 0 .
\Eq(2.1) $$
We shall denote by $\MM$ the set of values $\m\in(0,\m_{0})$
which satisfy \equ(2.1). For $\m\in\MM$ and $\e_{0}$ small enough
we shall restrict $\e$ to a large relative
measure set $\EE_{0}(\g) \subset (0,\e_{0})$ by imposing
the Diophantine conditions (recall that $\o=D+\m-\e$)
$$ \EE_{0}(\g):= \left\{ \e \in (0,\e_{0}) : \left| \o n - p \right|
\geq {\g\over n^{\t_{1}}} \quad \forall (n,p)\in
\NNN^{2} \right\}
\Eq(2.2)$$
for some $\t_{1}>\t_{0}+1$ and $\g\leq \g_{0}/2$; see
Appendix \secc(A1). These conditions guarantee the ``separation
of the resonant sites'', due to the regularising non-linearity,
for all pairs $(n,m)$ and $(n',m')$ such that $n\neq n'$;
indeed we have the following result.

\*
\0\lm(L1)
{\it Fix $s_{0}\in\RRR$. For all $\e\in \EE_{0}(\g)$ if for some 
$p\ge p_{1},n,n_{1}\in \NNN$ one has
$$ p^{s_{0}}|\o n-p-\m |\le \g/2 , \qquad
p_{1}^{s_{0}}|\o n_{1}-p_{1}-\m |\le \g/2 ,
\Eq(2.3)$$
then either $n=n_{1}$ and $p=p_{1}$ or $|n-n_{1}|\ge p_{1}^{s_{0}/
\t_{1}}$ and $n+n_{1}\ge B_{0} p_{1}$ for some constant $B_{0}$.}
\*

\0{\bf Proof.} If $n-n_{1}\neq 0$ one has $\g/|n-n_{1}|^{\t_{1}}
\le|\o(n-n_{1})-(p-p_{1})| \le \g/p_{1}^{s_{0}}$,
so that  one obtains $p_{1}^{s_{0}} \le |n-n_{1}|^{\t_{1}}$. If
$n=n_{1}$ then $|p-p_{1}| \le \g/p_{1}^{s_{0}}$, hence $p=p_{1}$.
Finally the inequality $n+n_{1}\ge B_{0} p_{1}$
follows immediately from \equ(2.3),
with the constant $B_{0}$ depending on $\o$ and $\m$. \qed

\*

\0{\bf Remark.} Note that if $s_{0}$ is small enough one can always
bound $B_{0}p_{1} \ge p_{1}^{s_{0}/\t_{1}}$.

\*

We shall now use the following lemma \cita{C}
to reorder our space index set $\ZZZ^{D}$.
The proof is deferred to Appendix \secc(A2) (see also \cita{Bo4}).

\*
\0\lm(L2)
{\it For all $\a>0$ small enough one can write
$\ZZZ^{D}= \cup_{j\in \nnn} \Lambda_{j}$ such that\\
(i) all $m\in \L_{j}$ are on the same sphere,
i.e. for all $j\in \NNN$ there exists $p_{j}\in\NNN$
such that $|m|^{2}\equiv  p_{j}$ $\forall m\in\Lambda_{j}$;\\
(ii) $\L_{j}$ has $d_{j}$ elements such
that $|\L_{j}|\equiv d_{j}\le C_{1} p_{j}^\a$,
for some $j$-independent constant $C_{1}$;\\
\\
(iii) for all $i\neq j$ such that $\L_{j}$ and $\L_{i}$ are
on the same sphere (i.e. such that $p_{j}= p_{i}$) one has
$$ {\rm dist}(\L_{i},\L_{j})\geq C_{2} p_{j}^\b, \qquad
\b = {2 \a\over 2D + (D+2)!D^{2}} ,
\Eq(2.4) $$
for some $j$-independent constant $C_{2}$;\\
(iv) if $d_{j}>1$ then for any $m\in\L_{j}$ there exists $m'\in\L_{j}$
such that $|m-m'| < C_{2} p_{j}^{\b}$, so that
one has ${\rm diam}(\L_{j}) \le C_{1} C_{2}p_{j}^{\a+\b}$;\\
If $D=2$ one can take $d_{j}=2$ for all $j$ and $\b=1/3$.}
\*

\0{\bf Remarks.}
(1) Essentially Lemma \secc(L2) assures that the points located on the
intersection of the lattice $\ZZZ^{D}$ with a sphere of any given
radius $r$ can be divided into a finite number of clusters,
containing each just a few elements (that is of order $r^{\a}$,
$\a\ll 1$) and not too closer to each other (that is at a distance
not less than of order $r^{\b}$, $\b>0$; in fact one has $\b<\a$).\\
(2) In fact the proof given in Appendix \secc(A2) shows that 
${\rm diam}(\Lambda_{j})< {\rm const.} p_{j}^{\a/D}$.

\*

By definition we call $\L_{1}$ the list of vectors $m$
such that $m_{i}=\pm 1$ (that is $p_{j}=D$).
In the following we shall take $\a \ll \min\{s,1\}$,
with $s$ given in \equ(1.2).

\*\*
\0{\bf 2.2. Renormalised $\PPPP$-$\QQQQ$ equations}\*

\0For $(n,j)\neq (1,1)$, let us define
$$ U_{n,j}=\{ {u_{n,m}}\}_{m\in \L_{j}} ,
\Eq(2.5) $$
which is a vector in $\RRR^{d_{j}}$.
Recall that $p_{j}=|m|^{2}$ if $m\in\L_{j}$;
the equations for $U_{n,j}$ are then by definition
$$ p_{j}^{s} \d_{n,j}U_{n,j} = \e F_{n,j} ,
\Eq(2.6) $$
where
$$ \d_{n,j}= - \o n + p_{j} + \m, \qquad F_{n,j}=
\{ f_{n,m}\}_{m\in \L_{j}} .
\Eq(2.7) $$

We introduce the $\e$-dependent
$$ y_{n,j} := p_{j}^{s_{2}} \d_{n,j} ,
\Eq(2.8) $$
where the exponent $s_{2}<s$ will be fixed in the forthcoming
Definition \secc(D2) (iv), and we define the {\it renormalised
$P$ equations} (for $(n,j)\neq (1,1)$) as
$$ p_{j}^{s} \left( \d_{n,j}I + p_{j}^{-s}\bar\chi_{1}(y_{n,j}) \,
M_{n,j} \right) U_{n,j}= \h F_{n,j} + L_{n,j}U_{n,j} ,
\Eq(2.9)$$
where $I$ (the identity), $M_{n,j}$ and $L_{n,j}$ are $d_{j}
\times d_{j}$ matrices and $\bar\chi_{1}$ is a $C^{\infty}$
non-increasing function such that (see Figure 2 below)
$$\left\{ \eqalign{
\bar\chi_{1}(x)=1 , & \qquad \hbox{if }
|x|<\g/8 , \cr
\bar\chi_{1}(x)=0 , & \qquad \hbox{if }
|x|>\g/4 , \cr} \right.
\Eq(2.10) $$
and $\bar\chi_{1}'(x)<C\g^{-1}$ for some positive constant $C$
(the prime denotes derivative with respect to the argument).

Clearly \equ(2.9)  coincides
with \equ(2.6), hence with \equ(1.11), provided
$$ \h=\e , \qquad \bar\chi_{1}(y_{n,j}) M_{n,j}=L_{n,j} ,
\Eq(2.11) $$
for all $(n,j)\neq(1,1)$. The matrices $L_{n,j}$ will be called
the {\it counterterms}.

We complete the renormalised $P$ equations
with the {\it renormalised $Q$ equations}
$$ D^{s} q = \sum_{n_{1}+n_{2}-n_3=1 \atop n_{i}=1}
\sum_{m_{1}+m_{2}-m_3=V \atop m_{i}\in \L_{1}}
u_{n_{1},m_{1}}u_{n_{2},m_{2}}u_{n_3,m_3} +
{\mathop{\sum}_{n_{1}+n_{2}-n_3=1 \atop
m_{1}+m_{2}-m_3=V}}^{\hskip-0.7truecm*}
u_{n_{1},m_{1}}u_{n_{2},m_{2}}u_{n_3,m_3} ,
\Eq(2.12) $$
where the symbol $\sum^*$ implies the restriction to the triples
of $(n_{i},m_{i})$ such that at least one has not $n_{i}=|m_{i}|^{2}=1$.
It should be noticed that the second sum vanishes at $\h=0$.

\*\*
\0{\bf 2.3. Matrix spaces}\*

\0Here we introduce some notations and properties that we shall need
in the following.

\*
\0\sub(D1)
{\it Let $A$ be a $d\times d$ real-symmetric matrix,
and denote with $A(i,j)$ and $\l^{(i)}(A)$ its entries and
its eigenvalues, respectively.
Given a list $\ul{m}:=\{m_{1},\dots,m_d\}$ with $m_{i}\in \ZZZ^{D}$
and a positive number $\s$, we define the norms
$$ \eqalign{
& \left| A \right|_{\io} := \max_{i,j\leq d}
|A(i,j)| , \qquad \left| A \right|_{\s,\ul{m}} :=
\max_{i,j\leq d} |A(i,j)| \, {\rm e}^{\s|m_{i}-m_{j}|^{\r}} , \cr
& \Vert A \Vert := {1\over \sqrt{d}} \sqrt{ {\rm tr}(A^{T}A) } = 
\sqrt{ {1\over d} \sum_{i,j=1}^{d} A(i,j)^{2}} , \cr}
\Eq(2.13) $$
with $\r$ depending on $D$. For fixed $\ul{m}=\{m_{1},\ldots,m_{d}\}
\in \ZZZ^{dD}$ we call $\AAA(\ul{m})$ the space of $d\times d$
real-symmetric matrices $A$ with norm $|A|_{\s,\ul{m}}$.}

\*
\0\lm(L3)
{\it Given a matrix $A\in \AAA(\ul{m})$, the following properties hold.\\
\0 (i) The norm $\Vert A\Vert$ is a smooth function in
the coefficients $A(i,j)$.\\
(ii) One has ${1\over \sqrt{d}}\Vert A\Vert\le |A|_\io\le
\sqrt{d}\Vert A\Vert $.\\
(iii) One has
${1\over \sqrt{d}} \max_{i} \sqrt{\l^{(i)}(A^{T}A)} \le
\Vert A\Vert \le \max_{i} \sqrt{\l^{(i)}(A^{T}A)}$.\\
(iv) For invertible $A$ one has
$$ \dpr_{A(i,j)}A^{-1}(h,l)= - A^{-1}(h,i)A^{-1}(j,l) , \qquad
\dpr_{A(i,j)}\Vert A\Vert = {A(i,j)\over d \Vert A\Vert} .
\Eq(2.14) $$
}

\0{\bf Proof.}
Item (i) follows by the invariance of the characteristic polynomial
under change of coordinates.

Items (iii) and (iv) are trivial.

The first relation in item (iv) follows by the definition of
differential as
$$ D_A f(A)[B]\equiv \dpr_\e f(A+\e B)\vert_{\e=0} .
\Eq(2.15) $$
Now by Taylor expansion we get $ D_A (A)^{-1}[B]= - A^{-1}BA^{-1}$.
The second relation is trivial. \qed

\*

\0{\bf Remark.} Note that for $A$ symmetric
one has $\sqrt{\l^{(i)}(A^{T}A)}= |\l^{(i)}(A)|$.

\*
\0\sub(D2)
{\it Let $\{\L_{j}\}_{j=1}^{\io}$ be the partition of $\ZZZ^{D}$
introduced in Lemma \secc(L2). Fix $\a$ small enough with
respect to $\min\{s,1\}$, with $s$ given in \equ(1.2).
Call $\O \subset \ZZZ \times \NNN$
the set of indexes $(n,j)\neq(1,1)$ such that
$$ - {1\over2} + \left( D+\m-\e_{0} \right) n < p_{j} <
\left( D+\m \right)n + {1\over 2} .
\Eq(2.16) $$
For $\e_{0}$ small enough \equ(2.16) in particular implies $n>0$,
hence $\O \subset \NNN^{2}$.
With each $(n,j)\neq (1,1)$
we associate the list $\L_{j}=\{m_{j}^{(1)},\dots,m_{j}^{(d_{j})}\}$,
with $d_{j} \le C_{1} p_{j}^{\a}$, and a $d_{j}\times d_{j}$
real-symmetric matrix $M_{n,j} \in \AAA(\Lambda_{j})$ (see
Definition \secc(D1)), such that $M_{n,j}=0$ if $(n,j)\notin \O$.\\
(i) We call $\MMM$ the space of all matrices which belong to
a space $\AAA(\Lambda_{j})$ for some $j\in\NNN$, and for $A\in
\AAA(\Lambda_{j})$ we set $|A|_{\s}=|A|_{\s,\Lambda_{j}}$.\\
(ii) We denote the eigenvalues of $\bar\chi_{1}(y_{n,j})M_{n,j}$
with $p_{j}^{\a} \n_{n,j}^{(i)}$, so that $\n_{n,j}^{(i)}\le
C |M_{n,j}|_{\io} \le C |M_{n,j}|_{\s}$,
for some constant $C$.\\
(iii) For invertible $\d_{n,j}I +p_{j}^{-s}\bar\chi_{1}(y_{n,j})M_{n,j}$
we define $x_{n,j} $ and $\n_{n,j}$ by setting
$$ x_{n,j}=\left|\d_{n,j}+p_{j}^{-s+ 2\a}\n_{n,j}\right|=
\left\Vert (\d_{n,j}I + p_{j}^{-s} \bar \chi_{1}(y_{n,j})M_{n,j})^{-1}
\right\Vert^{-1} ,
\Eq(2.17) $$
where the norm $\Vert A\Vert$ is introduced in Definition \secc(D1) --
notice that $\n_{n,j}$, hence $x_{n,j}$, depends both on $\e$ and $M$;\\
(iv) We call $s_{1}=s-2\a$ and set $s_{2}=s_{1}/4$ in \equ(2.8).}
\*

\0{\bf Remark.} Note that the eigenvalues $\n_{n,j}^{(i)}$ 
are proportional to  $\bar\chi_{1}(y_{n,j})$, hence vanish
for $|y_{n,j}|>\g/4$.

\*
\0\lm(L4) {\it There exists a positive constant $C$ such that
one has $|\n_{n,j}|\le C|M_{n,j}|_{\io} \le C |M_{n,j}|_{\s}$.}
\*

\0{\bf Proof.} For notational simplicity set $M_{n,j}=M$, $\d_{n,j}=\d$,
$p_{j}=p$, $d_{j}=d$, $x_{n,j}=x$, $\n_{n,j}=\n$, $\n^{(i)}_{n,j}=
\n_{i}$, and define $\l_{i}=\d+p^{-s+\a} \n_{i}$, with $|\n_{i}|
\le C |M|_{\io}$ (see Definition \secc(D2) (ii)). Then one has
$$ x = \left| \d + p^{-s+2\a} \n \right| = \left( {1 \over d}
\sum_{i=1}^{d} {1 \over \l_{i}^{2}} \right)^{-1/2} \le
C_{1}^{1/2} p^{\a/2} \min_{i} |\l_{i}| \le C_{1}^{1/2}
p^{\a/2} \left( |\d| + p^{-s+\a} \min_{i} |\n_{i}| \right) . $$
We distinguish between two cases.
\vskip.2truecm
\01. If there exists $i=i_{0}$ such that $|\d|< 2p^{-s+\a}|\n_{i_{0}}|$
then one obtains
$$ x \le  2 C_{1}^{1/2} p^{-s+3\a/2} |\n_{i_{0}}| + p^{-s+3\a/2}
\min_{i}|\n_{i}| \le 4 C_{1}^{1/2} p^{-s+2\a}|\n_{i_{0}}| . $$
Therefore, if $|\d|<p^{-s+2\a}|\n|/2$ one has
$$ p^{-s+2\a}|\n|/2 < x < 4 C_{1}^{1/2} p^{-s+2\a}|\n_{i_{0}}| \le
4 C C_{1}^{1/2} p^{-s+2\a}|M|_{\io}, $$
hence $|\n| \le {\rm const.}|M|_{\io}$.
If $|\d| \ge p^{-s+2\a}|\n|/2$ one has, by the assumption on $\d$,
$p^{-s+2\a}|\n|/2 \le |\d| < 2p^{-s+\a}|\n_{i_{0}}| \le
4p^{-s+2\a}|\n_{i_{0}}|$, and the same bound follows.
\vskip.2truecm
\02. If $|\d| \ge 2p^{-s+\a}|\n_{i}|$ for all $i=1,\ldots,d$, then one has
$$ x = |\d | \left( {1 \over d}
\sum_{i=1}^{d} {1 \over (1 + \d^{-1}p^{-s+\a}\n_{i})^{2}} \right)^{-1/2} 
= |\d| + O(p^{-s+\a}\max_{i}\n_{i}) , $$
so that $|\n| \le {\rm const.}p^{-\a}C |M|_{\io}$. \qed

\*

\0{\bf Remark.} The space of lists $M=\{M_{n,j}\}_{(n,j)\in
\nnn^{2}}$ such that $M_{n,j}\in\MMM$ (cf. Definition \secc(D2) (i))
and $|M|_{\s}= \sup_{n,j}|M_{n,j}|_{\s}<
\infty$ is a Banach space, that we denote with $\BBB$.

\*

\0\sub(D3) 
{\it We define $\DD_{0} =\{(\e,M) : 0<\e \le \e_{0} , \;
|M|_\s\le C_{0} \e_{0} \}$, for a suitable positive constant $C_{0}$,
and $ \DD(\g)\subset \DD_{0}$ as the set of all
$(\e,M)\in \DD_{0}$ such that $\e\in \EE_{0}(\g)$ and
$$ \left| \o n - \left( p_{j}+\m+ {\n_{n,j}\over p_{j}^{s_{1}}}
\right) \right| \ge {\g\over |n|^\t} \qquad \forall
(n,j) \in \O , \quad (n,j) \neq (1,1) , \quad n \neq 0 ,
\Eq(2.18) $$
for some $\t>\t_{0}+1+D$.}
\*

\0{\bf Remark.}
We shall call {\it Melnikov conditions} the Diophantine conditions
in \equ(2.2) and \equ(2.18). We shall call \equ(2.2) the
{\it second Melnikov conditions}, as they will be used to bound
the difference of the momenta of comparable lines of the
forthcoming tree formalism.

\*\*
\0{\bf 2.4. Main propositions}\*

\0We state the propositions which represent our main technical
results. Theorem 1  is an immediate consequence of
Propositions 1 and 2 below.

\*
\0{\bf Proposition 1.}
{\it Assume that $(\e,M)\in \DD(\g)$. There exist positive constants
$c_{0},K_{0},K_{1},\s,\h_{0},Q_{0}$ such that
the following holds true. It is possible to find 
a sequence of matrices $L\in \BBB$,
$$ L:= \left\{ L_{n,j}(\h,\e,M;q) \right\}_{(n,j)\in
\nnn^{2}\setminus \{( 1,1)\}} ,
\Eq(2.19) $$  
such that the following holds.\\
(i) There exists a unique solution $U_{n,j}(\h,M,\e;q)$, with
$(n,j)\in \ZZZ\times \NNN\setminus \{(1,1)\}$,
of equation \equ(2.9) which is analytic in $\h,q$ for
$|\h|\leq \h_{0}$, $|q|\le Q_{0}$, $\h_{0} Q_{0}^{2} \le c_{0}$
and such that
$$ \left| U_{n,j}(\h,M,\e;q)(a) \right|
\le |\h|q^3 K_{0} {\rm e}^{- \s (|n|+|p_{j}|^{1/2})} .
\Eq(2.20) $$
(ii) The sequence $L_{n,j}(\h,\e,M;q)$ is analytic in $\h$ and
uniformly bounded for $(\e,M)\in \DD(\g)$ as
$$ |L(\h,\e,M;q)|_\s \leq  K_{0}|\h|q^{2} .
\Eq(2.21) $$
(iii) The functions $U_{n,j}(\h,\e,M;q)$ and $L_{n,j}(\h,\e,M;q)$  can be
extended on the set $\DD_{0}$ to $C^1$ functions,
denoted by $U^{E}_{n,j}(\h,\e,M;q)$ and $L^{E}_{n,j}(\h,\e,M;q)$,
such that
$$ L_{n,j}^{E}(\h,\e,M;q)=  L_{n,j}(\h,\e,M;q) , \qquad
U_{n,j}^{E}(\h,\e,M;q)=  U_{n,j}(\h,\e,M;q) ,
\Eq(2.22)$$
for all $(\e,M)\in \DD(2\g)$.\\
(iv) The extended Q-equation, obtained from \equ(2.12) by substituting
$U_{n,j}(\h,\e,M;q)$ with $U^{E}_{n,j}(\h,\e,M;q)$, has a solution
$q^{E}(\h,\e,M)$, which is a true solution of \equ(2.12)
for $(\e,M)\in \DD(2\g)$; with an abuse of notation we shall call 
$$U^{E}_{n,j}(\h,\e,M)= U^{E}_{n,j}(\h,\e,M; q^{E}(\h,\e,M)) , \qquad
L^{E}_{n,j}(\h,\e,M)= L^{E}_{n,j}(\h,\e,M; q^{E}(\h,\e,M)) . $$\\
(v) The functions $L_{n,j}^{E}(\h,\e,M) $ satisfy the bounds
$$ \eqalign{
& |L^{E}(\h,\e,M)|_\s \leq |\h| K_{1}\,,
\quad |\dpr_\e L^{E}_{n,j}(\h,\e,M)|_{\s}
\leq |\h| K_{1} |n|^{1+s_{2}} , \cr
& \sum_{(n,j)\in \O}\sum_{a,b=1}^{d_{j}}
\left| \dpr_{M_{n,j}(a,b)} L^{E}(\h,\e,M) \right|_\s
{\rm e}^{-\s|m_{a}-m_{b}|^{\r}} \le |\h| K_{1} , \cr}
\Eq(2.23) $$
with $\r$ depending on $D$, and one has
$$ \left| U^{E}_{n,j}(\h,\e,M) \right|
\le |\h|  K_{1} {\rm e}^{-\s (|n|+|p_{j}|^{1/2})} ,
\Eq(2.24) $$
uniformly for $(\e,M)\in \DD_{0}$.}
\*

Once we have proved Proposition 1, we solve the compatibility equation
for the extended counterterm function $L_{n,m}^{E}(\h=\e,\e,M)$, which is
well defined provided we choose $\e_{0}$ so that $\e_{0}<\h_{0}$.

\*
\0{\bf Proposition 2.}
{\it For all $(n,j)\in \O$, there exist  $C^1$ functions 
$M_{n,j}(\e):(0,\e_{0})\to \DD_{0}$ (with an appropriate choice of
$C_0$) such that\\
(i) $M_{n,j}(\e)$ verifies
$$\bar\chi_{1}(y_{n,j})M_{n,j}(\e)=  L^{E}_{n,j}(\e,\e,M(\e)) , 
\Eq(2.25) $$ 
and is such that
$$ \left| M_{n,j}(\e) \right|_\s\leq K_{2} \e , \qquad
\left| \dpr_\e M_{n,j}(\e) \right|_{\s}
\le K_{2} \left( 1 + |\e n | \right) |n|^{s_{2}} ,
\Eq(2.26) $$ 
for a suitable constant $K_{2}$;\\
(ii) the set $\CCCCCC \equiv \CCCCCC(2\g)$, defined as
$$ \CCCCCC= \left\{ \e \in \EE_{0}(\g) : (\e,M(\e))\in \DD(2\g) \right\} ,
\Eq(2.27) $$
has large relative Lebesgue measure, namely
$\lim_{\e\to 0^+}\e^{-1}{\rm meas}(\CCCCCC \cap (0,\e))=1$.
}
\*

\0{\bf Proof of Theorem 1.} 
By proposition 1 (i) for all $(\e,M)\in \DD(\g)$ we can find a
sequence $L_{n,j}(\h,\e,M)$ so that there  exists a unique
solution $U_{n,j}(\h,\e,M)$ of \equ(2.6) for all $|\h|\le \h_{0}$,
where $\h_{0}$ depends only on $\g$ for $\e_{0}$ small enough.
By Proposition 1 (iii) the sequence $L_{n,j}(\h,\e,M)$ and the
solution $U_{n,j}(\h,\e,M)$ can be extended 
to $C^1$ functions (denoted by $L^{E}(\h,\e,M)$ and
$U^{E}(\h,\e,M)$) for all $(\e,M)\in D$.
Moreover  $L^{E}_{n,j}(\h,\e,M)= L_{n,j}(\h,\e,M)$ and 
$U^{E}_{n,j}(\h,\e,M)= U_{n,j}(\h,\e,M) $ for all $(\e,M)\in \DD(2\g)$.

Equation  \equ(2.8) coincides with our original \equ(2.6)
provided the compatibility equations \equ(2.10) 
are satisfied. Now we fix $\e_{0}<\h_{0}$ so that
$L^{E}_{n,m}(\h=\e,\e,M)$ and  $U^{E}_{n,j}(\h=\e,\e,M)$ are well defined.
By Proposition 2 (i) there exists a sequence of matrices $M_{n,j}(\e)$
which satisfies the extended compatibility equations \equ(2.24).
Finally by Proposition 2 (ii) the Cantor set 
$\CCCCCC(2\g)$ is well defined and of large relative measure.
 
For all $\e\in \CCCCCC(2\g)$ the pair $(\e,M(\e))$ is by definition
in $\DD(2\g)$ so  that by Proposition 1 (iii) one has
$$ L_{n,j}(\e,\e,M(\e))=L^{E}_{n,j}(\e,\e,M(\e)) ,
\qquad u(\e,\e,M(\e);x,t)=u^{E}(\e,\e,M(\e);x,t) ,
\Eq(2.28) $$
so that  $U_{n,j}(\e,\e,M(\e))$ solves \equ(2.8) for $\h=\e$.
So by Proposition 2 (i) $M(\e)$ solves the true compatibility
equations \equ(2.10), $\bar\chi_{1}(y_{n,j})M_{n,j}(\e)=L_{n,j}
(\e,\e,M(\e))$, for all $\e\in \CCCCCC(2\g)$. Then $u(\e,\e,M(\e);x,t)$ 
is a true nontrivial solution of our \equ(1.9) in $\CCCCCC(2\g)$.
Then by setting $\EE(\mu)=\CCCCCC(2\g)$ the result follows.\qed

\*\*
\section(3,Recursive equations and tree expansion)

\0In this section we find a formal solution
$U_{n,j}$ of \equ(2.9) as a power series on $\h$;
the solution $U_{n,j}$ is parameterised by the matrices 
$L_{n,j}$ and it will be written in the form of a tree expansion.

We assume for  $L_{n,j}(\h,\e,M)$ and $U_{n,j}(\h,\e,M)$,
with $(n,j)\neq (1,1)$, a formal series expansion in $\h$, i.e.
$$ L_{n,j}(\h,\e,M)= \sum_{k=1}^{\io}\h^{k} L_{n,j}^\ka , \qquad
U_{n,j}(\h,\e,M)= \sum_{k=1}^{\io}\h^{k} U_{n,j}^{(k)} ,
\Eq(3.1) $$
for all $(n,j)\neq(1,1)$. Note that \equ(3.1) naturally 
defines the vector components $u^{(k)}_{n,m}$, $m\in\L_{j}$.

By definition we set
$$ U_{1,1}^{(0)}=\{ u_{1,m} \,: m\in \L_{1}\},\qquad
u_{1,V}=q , \qquad U_{1,1}^\ka=0 , \; k\neq 0 ,
\Eq(3.2) $$
where $V=(1,1,\ldots,1)$.
Inserting the series expansion in \equ(2.9)
we obtain for all $(n,j)\neq (1,1)$ the recursive equations
$$ p_{j}^{s} \left( \d_{n,j}I +p_{j}^{-s}
\bar\chi_{1}(y_{n,j}) M_{n,j} \right)
U^\ka_{n,j}=  F^\ka_{n,j} +
\sum_{r=1}^{k-1}L^{(r)}_{n,j}U^{(k-r)}_{n,j} ,
\Eq(3.3) $$
while for $(n,j)=(1,1)$ we have
$$ q= f_{1,V} .
\Eq(3.4) $$
In \equ(3.3), for $m_{a} \in \L_{j}$, where $a=1,\ldots,d_{j}$,
$F^\ka_{n,j}(a)$ is defined as
$$ F^\ka_{n,j}(a)  = \sum_{k_{1}+k_{2}+k_3=k-1}
\sum_{n_{1}+n_{2}-n_3=n\atop m_{1}+m_{2}-m_3= m_a}
u^{(k_{1})}_{n_{1},m_{1}} u^{(k_{2})}_{n_{2},m_{2}} u^{(k_3)}_{n_3,m_3},
\Eq(3.5) $$
where each $u^{(k_{i})}_{n_{i},m_{i}}$ is a component of
some $U^{(k_{i})}_{n_{i},j_{i}}$. Recall that we are
assuming for the time being $f(u,\bar u)=|u|^{2}u$ and we are
looking for solutions with real Fourier coefficients $u_{n,m}$.

\*\*
\0{\bf 3.1. Multiscale analysis}\*

\0It is convenient to rewrite  \equ(3.3) introducing the following
scale functions.

\*
\0\sub(D4)
{\it Let $\chi(x)$ be a $C^{\io}$ non-increasing function such that
$\chi(x)=0$ if $|x|\ge 2\g$ and $\chi(x)=1$ if $|x|\le \g$;
moreover, if the prime denotes derivative with respect to the argument,
one has $|\chi'(x)|\le C\g^{-1}$ for some positive constant $C$.
Let $\chi_{h}(x)=\chi(2^{h}x)-\chi(2^{h+1}x)$ for
$h \ge 0$, and $\chi_{-1}(x)=1-\chi(x)$; see Figure 1. Then
$$ 1 = \chi_{-1}(x)+\sum_{h=0}^{\io} \chi_{h}(x)=
\sum_{h=-1}^{\io} \chi_{h}(x) .
\Eq(3.6) $$
We can also write
$$ 1 = \bar\chi_{1}(x)+\bar\chi_{0}(x)+\bar\chi_{-1}(x) ,
\Eq(3.6a) $$
with $\bar\chi_{1}(x)=\chi(8x)$ (cf. \equ(2.8) and Figure 2),
$\bar\chi_{-1}(x)=1-\chi(4x)$,
and $\bar\chi_{0}(x)=\chi_{2}(x)=\chi(4x)-\chi(8x)$.}

\midinsert
\*
\insertplotttt{300pt}{144pt}{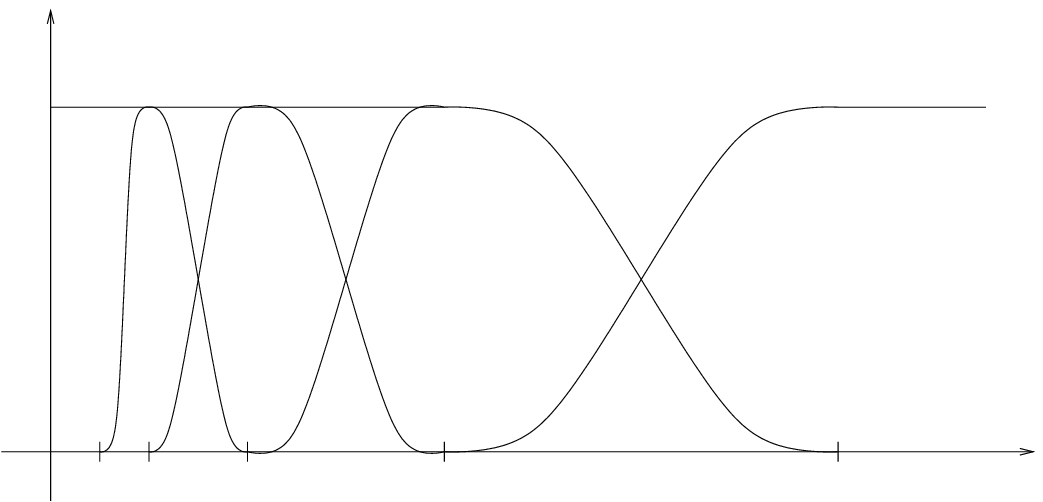}
\ins{370pt}{20pt}{$x$}
\ins{314pt}{20pt}{$2\g$}
\ins{203pt}{20pt}{$\g$}
\ins{144pt}{20pt}{$\g/2$}
\ins{116pt}{20pt}{$\g/4$}
\ins{094pt}{20pt}{$\g/8$}
\ins{150pt}{144pt}{$\chi_{1}(x)$}
\ins{118pt}{144pt}{$\chi_{2}(x)$}
\ins{210pt}{144pt}{$\chi_{0}(x)$}
\ins{320pt}{144pt}{$\chi_{-1}(x)$}
\ins{095pt}{144pt}{$\chi(x)$}
\line{\vtop{\line{\hskip1.2truecm\vbox{\advance\hsize by -2.5 truecm
\noindent{\nota {\bf Figure 1.} Graphs of some of the $C^{\io}$ compact
support functions $\chi_{h}(x)$ partitioning the unity.
The function $\chi(x)$ is given by the envelope of all functions
but $\chi_{-1}(x)$.\vfil
}} \hfill} }}
\*
\endinsert

\midinsert
\*
\insertplotttt{300pt}{144pt}{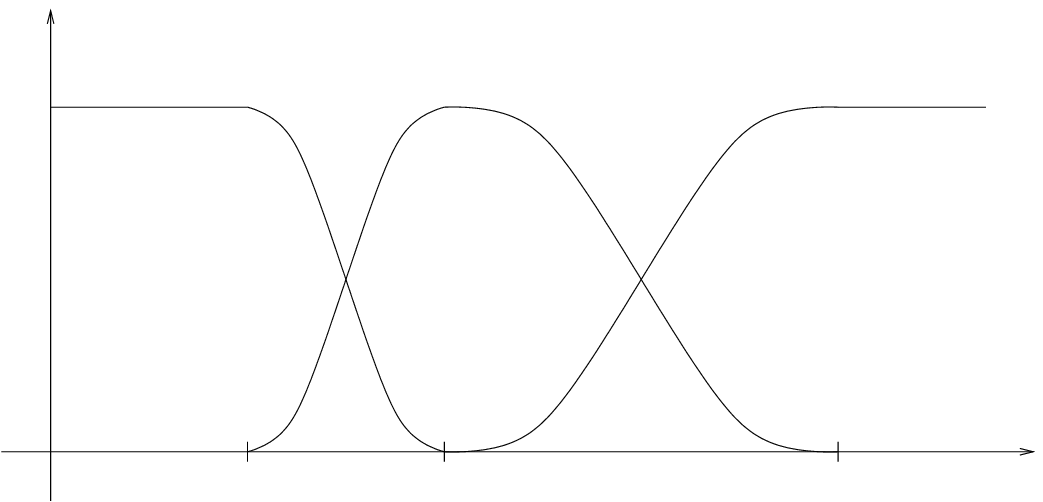}
\ins{370pt}{20pt}{$x$}
\ins{310pt}{20pt}{$\g/2$}
\ins{200pt}{20pt}{$\g/4$}
\ins{140pt}{20pt}{$\g/8$}
\ins{120pt}{144pt}{$\bar\chi_{1}(x)$}
\ins{210pt}{144pt}{$\bar\chi_{0}(x)$}
\ins{320pt}{144pt}{$\bar\chi_{-1}(x)$}
\centerline{
\noindent{{\nota {\bf Figure 2.} Graphs of the $C^{\io}$ functions
partitioning the unity $\bar \chi_{-1}(x)$, $\bar\chi_{0}(x)$
and $\bar\chi_{1}(x)$.}
}}
\*
\endinsert

\*

\0{\bf Remark.} Note that $\chi_{h}(x) \neq 0$ implies
$2^{-h-1}\g < |x| < 2^{-h+1}\g$ if $h\ge 0$ and
$\g < |x|$ if $h=-1$. In particular if $\chi_{h}(x) \neq 0$
and $\chi_{h'}(x) \neq 0$ for $h \neq h'$ then $|h-h'|=1$.

\*
\0\sub(D5)
{\it We denote (recall \equ(2.17) and that $s_{1}=s-2\a$)
$$ x_{n,j} \equiv x_{n,j}(\e,M) = \left|
\d_{n,j}+{\n_{n,j}\over p_{j}^{s_{1}}} \right| .
\Eq(3.7)$$
For $h=-1,0,1,2,\ldots,\io$ and $i=-1,0,1$ we define
$G_{n,j,h,i}(\e,M)$ as follows:\\
\0(i) for $i=-1,0$, we set $G_{n,j,h,i}=0$ for $h \neq -1$
and $G_{n,j,-1,i}=0$ for all $(\e,M)$
such that $\bar\chi_{i}(y_{n,j})=0$;\\
\0(ii) similarly we set $G_{n,j,h,1}=0$ for all $(\e,M)$
such that  $\chi_{h}(x_{n,j})=0$;\\
\0(iii) otherwise we set
$$ \left\{
\eqalign{
G_{n,j,-1,i} & = \bar\chi_{i}(y_{n,j})
p_{j}^{-s}\left( \d_{n,j}I+{\bar\chi_{1}(y_{n,j})\,M_{n,j}\over
p_{j}^{s}} \right)^{-1} , \qquad i=-1,0 , \cr
G_{n,j,h,1} & = \bar\chi_{1}(y_{n,j}) \chi_{h}(x_{n,j})
p_{j}^{-s}\left( \d_{n,j}I+{\bar\chi_{1}(y_{n,j})\,M_{n,j}\over
p_{j}^{s}} \right)^{-1} , \qquad h \ge -1 . \cr} \right.
\Eq(3.8)$$
Then $G_{n,j,h,i}$  will be called the {\rm propagator} on scale $h$.}
\*

\0{\bf Remarks.}
(1) If $p_{j}^{\a}\n^{(i)}_{n,j}$ are the eigenvalues of
$\bar\chi_{1}(y_{n,j}) M_{n,j}$ (cf. Definition \secc(D2))
one has by Lemma \secc(L3)
$$ \min_{i} \left| \d_{n,j}I+ p_{j}^{-s+\a} \n^{(i)}_{n,j} \right|
\le x_{n,j} \le \min_{i} \sqrt{d_{j}} \left| \d_{n,j}I
+ p_{j}^{-s+\a} \n^{(i)}_{n,j} \right| ,
\Eq(3.9) $$
so that $\d_{n,j}I+p_{j}^{-s} \bar\chi_{1}(y_{n,j})\,M_{n,j}$ is
invertible where $G_{n,j,h,i}(\e,M)$ is not identically zero;
this implies that $G_{n,j,h,i}(\e,M)$ is well defined (and
$C^{\io}$) on all $\DD_{0}$ (as given in Definition \secc(D3)).\\
(2) If $i=-1,0$, then for $(\e,M)\in\DD_{0}$ the denominators are
large. Indeed $i\neq1$ implies $|y_{n,j}|\ge \g/8$,
hence $|\d_{n,j}|\ge p_{j}^{-s_{2}}\g/8$, whereas $|p_{j}^{-s_{1}}
\n_{n,j}| \le p_{j}^{-s_{1}} C C_{0} |\e_{0}|\le {\rm const.}
p_{j}^{-s_{2}} \e_{0}$ in $\DD_{0}$ (with $C$ as in Lemma \secc(L4)
and $C_{0}$ as in Definition \secc(D3)), so that
$x_{n,j}=|\d_{n,j}+ p_{j}^{-s_{1}} \n_{n,j}| \ge |\d_{n,j}|/2$. Then
$$ \eqalign{
\left| G_{n,j,-1,i} \right|_{\io} & =
p_{j}^{-s} \Big| \Big( \d_{n,j}I +
{\bar\chi_{1}(y_{n,j})\,M_{n,j}\over p_{j}^{s}} \Big)^{-1} \Big|_\io
\cr
& \le C_{1}^{1/2} p_{j}^{-s+\a/2}
\Big| \d_{n,j}+ {\n_{n,j}\over p_{j}^{s_{1}}}
\Big|^{-1} \le 2 C_{1}^{1/2} p_{j}^{-s+\a/2+s_{2}} |y_{n,j}|^{-1}
\le {16 \over \g} C_{1}^{1/2} p_{j}^{-3s/4} , \cr}
\Eq(3.10) $$
where we have also used Lemma \secc(L3) (ii).\\
(3) Notice that  $G_{n,j,-1,-1}$ is a diagonal matrix
(cf. \equ(3.8) and notice that $\bar\chi_{-1}
(y_{n,j})\bar\chi_{1}(y_{n,j})=0$ identically).

\*


Inserting the multiscale decomposition \equ(3.6) and \equ(3.6a)
into \equ(3.3) we obtain
$$ U^\ka_{n,j} = \sum_{i=-1,0,1}
\sum_{h=-1}^\infty U^\ka_{n,j,h,i} ,
\Eq(3.11) $$
with
$$ U^\ka_{n,j,h,i}= G_{n,j,h,i} F^\ka_{n,j}+ \d(i,1)\,G_{n,j,h,1}
\left(\sum_{h_{1}=-1}^{\io}\sum_{i_{1}=0,1}
\sum_{r=1}^{k-1}L^{(r)}_{n,j,h}U^{(k-r)}_{n,j,h_{1},i_{1}} \right) ,
\Eq(3.12) $$
where $\d(i,j)$ is Kronecker's delta, and
we have used that $h=-1$ for $i \neq 1$ and written
$$ L^{(r)}_{n,j} = \sum_{h=-1}^{\io}
\bar\chi_{1}(y_{n,j}) \chi_{h}(x_{n,j}) \, L^{(r)}_{n,j,h} ,
\Eq(3.13) $$
with the functions $L^{(r)}_{n,j,h}$ to be determined.

\*\*
\0{\bf 3.2. Tree expansion}\*

\0The equations \equ(3.12) can be applied recursively until we
obtain the Fourier components $u_{n,m}^\ka$ 
as (formal) polynomials in the variables $G_{n,j,h,i}$,
$q$ and $L^{(r)}_{n,j,h} $ with $r<k$. It turns out that 
$u_{n,m}^\ka$ can be written as sums over {\it trees}
(see Lemma \secc(L6) below), defined in the following way.

A (connected) graph $\GG$ is a collection of points (vertices)
and lines connecting all of them. The points of a graph
are most commonly known as graph vertices, but may also be
called {\it nodes} or points. Similarly, the lines connecting
the nodes of a graph are most commonly
known as graph edges, but may also be called branches or
simply {\it lines}, as we shall do. We denote with
$V(\GG)$ and $L(\GG)$ the set of nodes and the set of lines,
respectively. A path between two nodes is the minimal subset of
$L(\GG)$ connecting the two nodes. A graph is planar if it
can be drawn in a plane without graph lines crossing.

\*
\0\sub(D6)
{\it A {\rm tree} is a planar graph $\GG$ containing no closed loops.
One can consider a tree $\GG$ with a single special node $v_{0}$:
this introduces a natural partial ordering on the set
of lines and nodes, and one can imagine that each line
carries an arrow pointing toward the node $v_{0}$.
We can add an extra (oriented) line $\ell_{0}$ exiting the special
node $v_{0}$; the added line will be called the {\rm root line}
and the point it enters (which is not a node) will
be called the {\rm root} of the tree. In this way we obtain
a {\rm rooted tree} $\th$ defined by $V(\th)=V(\GG)$
and $L(\th)=L(\GG)\cup\ell_{0}$. A {\rm labelled tree} is a
rooted tree $\th$ together with a label function defined on
the sets $L(\th)$ and $V(\th)$.}
\*

We shall call {\it equivalent} two rooted trees which can be transformed
into each other by continuously deforming the lines in the plane
in such a way that the latter do not cross each other
(i.e. without destroying the graph structure).
We can extend the notion of equivalence also to labelled trees,
simply by considering equivalent two labelled trees if they
can be transformed into each other in such a way that also
the labels match. An example of tree is illustrated in Figure 3.

\midinsert
\insertplotttt{290pt}{120pt}{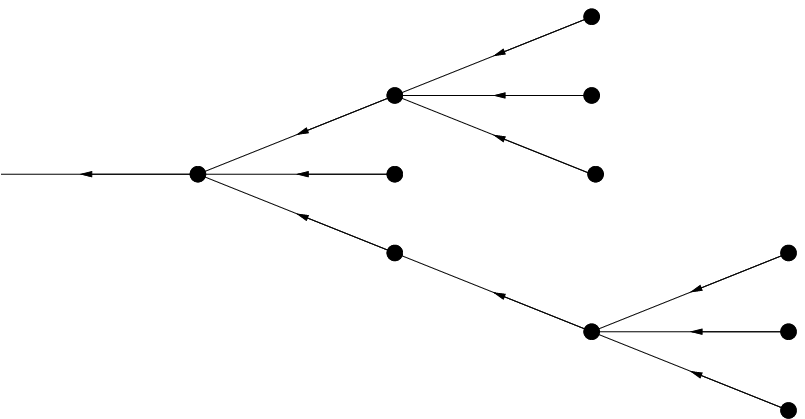}
\line{\vtop{\line{\hskip1.2truecm\vbox{\advance\hsize by -2.5 truecm
\*
\noindent{\nota {\bf Figure 3.}
Example of an unlabelled tree (only internal nodes with
1 and 3 entering lines are taken into account, according to
the diagrammatic rules in Section 3.3).\vfil}} \hfill} }}
\*
\endinsert

Given two nodes $v,w\in V(\th)$, we say that $w \prec v$
if $v$ is on the path connecting $w$ to the root line.
We can identify a line with the nodes it connects;
given a line $\ell=(v,w)$ we say that $\ell$
enters $v$ and exits (or comes out of) $w$. 
Given two comparable lines $\ell$ and $\ell_{1}$,
with $\ell_{1} \prec \ell$, we denote with $\PPP(\ell_{1},\ell)$
the path of lines connecting $\ell_{1}$ to $\ell$; by definition
the two lines $\ell$ and $\ell_{1}$ do not belong to $\PPP(\ell_{1},
\ell)$. We say that a node $v$ is along the path $\PPP(\ell_{1},
\ell)$ if at least one line entering or exiting $v$ belongs to the path.
If $\PPP(\ell_{1},\ell)=\emptyset$ there is only one node $v$
along the path (such that $\ell_{1}$ enters $v$ and $\ell$ exits $v$).

In the following we shall deal mostly with labelled trees:
for simplicity, where no confusion can arise, we shall call them just
trees. 

We call {\it internal nodes} the nodes such that there is at least
one line entering them; we call {\it internal lines} the lines
exiting the internal nodes. We call {\it end-points} the nodes
which have no entering line. We denote with $L(\th)$, $V_{0}(\th)$
and $E(\th)$ the set of lines, internal nodes and end-points,
respectively. Of course $V(\th)=V_{0}(\th)\cup E(\th)$.

\*\*
\0{\bf 3.3. Diagrammatic rules}\*

\0We associate with the nodes (internal nodes and end-points)
and lines of any tree $\th$ some labels, according to the
following rules; see Figure 4 for reference.

\midinsert
\*
\insertplotttt{360pt}{90pt}{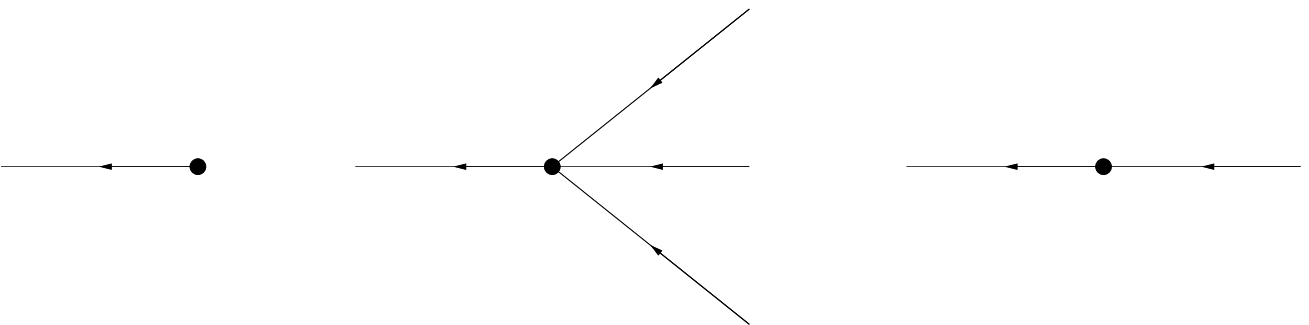}
\ins{030pt}{64pt}{(a)}
\ins{131pt}{64pt}{(b)}
\ins{291pt}{64pt}{(c)}
\ins{091pt}{53pt}{$v$}
\ins{202pt}{53pt}{$v$}
\ins{363pt}{53pt}{$v$}
\ins{067pt}{56pt}{$\ell$}
\ins{170pt}{56pt}{$\ell$}
\ins{334pt}{56pt}{$\ell$}
\ins{390pt}{56pt}{$\ell_{1}$}
\line{\vtop{\line{\hskip1.2truecm\vbox{\advance\hsize by -2.5 truecm
\noindent{\nota {\bf Figure 4.}
Labels associated to the nodes and lines of the trees.
(a) The line $\ell$ exits the end-point $v$:
one associate with $\ell$ the labels
$i_{\ell}$, $h_{\ell}$, $n_{\ell}$ and $m_{\ell}$,
and with $v$ the labels $n_{v}$, $m_{v}$ and $k_{v}$, with
the constraints $i_{\ell}=-1$, $h_{\ell}=-1$, $n_{\ell}=n_{v}=1$,
$m_{\ell}=m_{v}\in\L_{1}$, $k_{v}=0$.
(b) The line $\ell$ exits the node $v$ with $s_{v}=3$:
one associate with $\ell$ the labels
$i_{\ell}$, $h_{\ell}$, $n_{\ell}$, $j_{\ell}$, $m_{\ell}$,
$m_{\ell}'$, $a_{\ell}$, $b_{\ell}$, and with $v$ the label $k_{v}$,
with the constraints $(n_{\ell},j_{\ell})\neq(1,1)$,
$m_{\ell}=\L_{j_{\ell}}(a_{\ell})$,
$m_{\ell}'=\L_{j_{\ell}}(b_{\ell})$, $k_{v}=1$.
(c) The line $\ell$ exits the node $v$ with $s_{v}=1$:
one associate with $\ell$ the labels
$i_{\ell}$, $h_{\ell}$, $n_{\ell}$, $j_{\ell}$, $m_{\ell}$,
$m_{\ell}'$, $a_{\ell}$, $b_{\ell}$, and with $v$ the
labels $k_{v}$, $a_{v}$, $b_{v}$, $j_{v}$ and $n_{v}$,
with the constraints $(n_{\ell},j_{\ell})\neq(1,1)$,
$m_{\ell}=\L_{j_{\ell}}(a_{\ell})$,
$m_{\ell}'=\L_{j_{\ell}}(b_{\ell})$, $k_{v}\ge1$,
$a_{v}=b_{\ell}$, $b_{v}=a_{\ell_{1}}$,
$n_{\ell}=n_{\ell_{1}}$, $j_{\ell}=j_{\ell_{1}}$.
Other constraints are listed in the text.\vfil
}} \hfill} }}
\*
\endinsert

\vskip.2truecm
\0(1) For each node $v$ there are $s_v$ entering lines,
with $s_{v}\in\{0,1,3\}$; if $s_v=0$ then $v\in E(\th)$.
\vskip.1truecm
\0(2) With each end-point $v\in E(\th)$ one associates the {\it mode}
labels $(n_{v},m_{v})$, with $m_{v}\in \L_{1}$
and $n_{v}=1$. One also associates with each end-point an {\it order}
label $k_v=0$, and a {\it node factor} $\h_v=\pm q$, with the
sign depending on the sign of the permutation from $m_v$ to $V$:
one can write $\h_{v}=(-1)^{|m_{v}-V|_{1}/2} q$, where
$|x|_{1}$ is the $l_{1}$-norm of $x$.
\vskip.1truecm
\0(3) With each line $\ell\in L(\th)$ not exiting an end-point, one
associates the {\it index} label $j_{\ell}\in\NNN$ and the {\it momenta}
$(n_{\ell},m_{\ell},m'_{\ell}) \in \ZZZ \times \ZZZ^{D} \times \ZZZ^{D}$
such that $(n_{\ell},j_{\ell})\neq (1,1)$ and $m_{\ell},m'_{\ell}\in
\L_{j_{\ell}}$. One has $p_{j_{\ell}}=|m_{\ell}|^{2}=|m_{\ell}'|^{2}$
(see Lemma \secc(L2) (ii) for notations).
The momenta  define $a_{\ell},b_{\ell}\in \{1,\ldots,d_{j}\}$,
with $d_{j_{\ell}}=|\L_{j_{\ell}}| \le C_{1} p_{j_{\ell}}^{\a}$, such that
$m_{\ell}= \L_{j_{\ell}}(a_{\ell})$, $m'_{\ell}=\L_{j_{\ell}}(b_{\ell})$. 
\vskip.1truecm
\0(4) With each line $\ell \in L(\th)$ not exiting an end-point
one associates a {\it type} label $i_{\ell}=-1,0,1$.
If $i_{\ell}=-1$ then $m_{\ell}=m'_{\ell}$.
\vskip.1truecm
\0(5) With each line $\ell\in L(\th)$ not exiting an end-point
one associates the scale label $h_{\ell}\in \NNN \cup\{-1,0\}$.
If $i_{\ell}=0,-1$ then $h_{\ell}=-1$; if two lines $\ell,\ell'$
have $(n_{\ell},j_{\ell})=(n_{\ell'},j_{\ell'})$,
then $|i_{\ell}-i_{\ell'}| \le 1$ and if moreover
$i_{\ell}=i_{\ell'}=1$ then also $|h_{\ell}-h_{\ell'}| \le 1$.
\vskip.1truecm
\0(6) If $\ell\in L(\th)$ exits an end-point $v$
then $h_{\ell}=-1$, $i_{\ell}=-1$,
$j_{\ell}=1$, $n_{\ell}=1$ and $m_{\ell}=m_{v}$.
\vskip.1truecm
\0(7) With each line $\ell\in L(\th)$ except the root line one
associates a sign $\s(\ell)=\pm 1$ such that for
all $v\in V_{0}(\th)$ one has
$$ 1 = \sum_{\ell\in L(v)}\s(\ell),
\Eq(3.14) $$
where $L(v)$ is the set of the $s_{v}$ lines entering $v$.
One does not associate any label $\s$ to the root line $\ell_{0}$.
\vskip.1truecm
\0(8) If $s_v=1$ the labels $n_{\ell_{1}},j_{\ell_{1}}$ of the line
entering $v$ are the same as the labels $n_{\ell},j_{\ell}$
of the line $\ell$ exiting $v$, and one defines
$j_{v}=j_{\ell}$, $a_{v}= b_{\ell}$, $b_{v}=a_{\ell_{1}}$.
With such $v$ one associates an order label $k_{v}\in\NNN$.
\vskip.1truecm
\0(9) If $s_{v}=3$ then $k_{v}=1$. If $\ell$
is the line exiting $v$ and $\ell_{1},\ell_{2},\ell_{3}$
are the lines entering $v$ one has
$$  n_{\ell} = \s(\ell_{1}) n_{\ell_{1}} + \s(\ell_{2}) n_{\ell_{2}} +
\s(\ell_3) n_{\ell_{3}} =
\sum_{\ell'\in L(v)} \s(\ell') n_{\ell'} 
\Eq(3.15) $$
and
$$ m_{\ell}'=
\s(\ell_{1}) m_{\ell_{1}} + \s(\ell_{2}) m_{\ell_{2}} +
\s(\ell_3) m_{\ell_{3}} =
\sum_{\ell' \in L(v)} \s(\ell') m_{\ell'} ,
\Eq(3.16) $$
with $L(v)$ defined after \equ(3.14).
\vskip.1truecm
\0(10) With each line $\ell\in L(\th)$ one associates the {\it propagator}
$$ g_{\ell} := G_{n_{\ell} ,j_{\ell},h_{\ell},i_{\ell}}
(a_{\ell},b_{\ell}) .
\Eq(3.17) $$
if $\ell$ does not exit an end-point and $g_{\ell}=1$ otherwise.
\vskip.1truecm
\0(11) With each internal node $v\in V_{0}(\th)$ one associates a
{\it node factor} $\h_{v}$ such that $\h_v=1/3$ for
$s_v=3$ and $\h_{v}=L^{(k_v)}_{n_{\ell},j_{\ell}}(a_v,b_v)$
for $s_{v}=1$.
\vskip.1truecm
\0(12) Finally one defines the {\it order} of a tree as
$$ k(\th)= \sum_{v\in V(\th)} k_v .
\Eq(3.18) $$

\*
\0\sub(D7)
{\it We call  $\Theta^{(k)}$ the set of all the nonequivalent
trees of order $k$ defined according to the diagrammatic rules.
We call  $\Theta^{(k)}_{n,m}$ the set of all the nonequivalent
trees of order $k$ and with labels $(n,m)$ associated to the root line.}
\*

\0\lm(L5)
{\it  For all $\th\in \Th^\ka$ and for all lines $\ell\in L(\th)$
one has $ |n_{\ell}|,|m_{\ell}|,|m'_{\ell}| \le B k$,
for some constant $B$.}
\* 

\0{\bf Proof.} By definition of order one has $|V_{0}(\th)|\le k$
and by induction one proves $|E(\th)| \le 2|V_{0}(\th)|+1$
(by using that $s_{v}\le 3$ for all $v\in V_{0}(\th)$).
Hence $|E(\th)|\le 2k+1$. Each end-point $v$ contributes
$n_{v}=\pm 1$ to the momentum $n_{\ell}$ of any line $\ell$ following
$v$, so that $|n_{\ell}|\le 2k+1$ for all lines $\ell\in L(\th)$.

Let $\th_{\ell}$ be the tree with root line $\ell$ and let
$k(\th_{\ell})$ be its order. Then the bounds $|m_{\ell}|,|m_{\ell}'|
\le 2 k(\th_{\ell})+1$ can be proved by induction on $k(\th_{\ell})$
as follows. If $v$ is the internal node which $\ell$ exits and
$s_{v}=3$, call $\ell_{1},\ell_{2},\ell_{3}$ the lines entering $v$
(the case $s_{v}=1$ can be discussed in the same way, and
it is even simpler) and for $i=1,\ldots,3$ denote by $\th_{i}$
the tree with root line $\ell_{i}$ and by $k_{i}$ the corresponding
order. Then $k_{1}+k_{2}+ k_{3}=k(\th_{\ell})-1$,
so that by the inductive hypothesis one has
$$ m_{\ell}' = m_{\ell_{1}} + m_{\ell_{2}} + m_{\ell_{3}}
\quad \Longrightarrow \quad
|m_{\ell}'| \le \sum_{i=1}^{3} \left( 2k_{i}+1 \right)
\le 2 k(\th_{\ell}) + 1 , $$
and hence also $|m_{\ell}|=|m_{\ell}'| \le 2k(\th_{\ell})+1$. \qed

\*

The coefficients $u^{(k)}_{ n,m}$ can be represented as sums over the
trees defined above; this is in fact the content of the following lemma.

\*
\0\lm(L6)
{\it The coefficients $u^{(k)}_{n,m}$ can be written as
$$ u^{(k)}_{n,m} = 
\sum_{\th\in\Theta^{(k)}_{n,m}} \Val(\th) ,
\Eq(3.19)$$
where
$$ \Val(\th) = 
\Big( \prod_{\ell \in L(\th)} g_{\ell} \Big)
\Big( \prod_{v \in V(\th)} \h_{{v}} \Big)\;.
\Eq(3.20) $$}
\*

\0{\bf Proof.} 
The proof is done by induction on $k\ge 1$.
For $k=1$ it reduces just to a trivial check.

Now, let us assume that \equ(3.19) holds for $k'<k$, and
use that $u^{(0)}_{n,m}= q \, \d(n, 1)
\prod_{i=1}^{D} (\pm \d(m_{i}, \pm 1))$.
If we set $m=\L_{j}(a)$, we have (see Figure 5)
$$ \eqalign{
u_{n,m}^{(k)} & =\sum_{h=-1}^{\io} \sum_{i=-1,0,1}
\sum_{b=1}^{d_{j}}
G_{n,j,h,i}(a,b) \sum_{k_{1}+k_{2}+k_3=k}
\sum_{n_{1}+n_{2}-n_3=n \atop m_{1}+m_{2}-m_3= \L_{j}(b)}
u_{n_{1},m_{1}}^{(k_{1})} u_{n_{2},m_{2}}^{(k_{2})}u_{n_3,m_3}^{(k_3)} \cr
& + \sum_{h=-1}^{\io} 
\sum_{b,b'=1}^{d_{j}} G_{n,j,h,1}(a,b)
\sum_{r=1}^{k-1} L_{n,j,h}(b,b') \, u_{n,\L_{j}(b')}^{(k-r)} . \cr}
\Eq(3.21) $$
Consider a tree $\th\in \Th_{n,m}^\ka$ such that  $m=\L_{j}(a)$,  
$s_{v_{0}}=3$ and $h_{\ell_{0}}=h$, if $\ell_{0}$ is the root line of
$\th$ and $v_{0}$ is defined in \secc(D6).
Let $\th_{1},\th_{2},\th_3$ be the sub-trees whose root lines
$\ell_{1},\ell_{2},\ell_3$ enter $v_{0}$. By \equ(3.14) one
has $\sum_{j=1}^{3} \s(\ell_{j})m_{\ell_{j}} = m_{\ell_{0}}'$,
with $m_{\ell_{0}}'=\L_{j}(b)$ for $b=b_{\ell_{0}}$. Then we have
$$ \Val(\th) =  G_{n,j,h,i}(a,b) \Val(\th_{1})\Val(\th_{2})\Val(\th_3),
\Eq(3.22) $$
and we reorder the lines so that $\s(\ell_3)=-1$,
which produces a factor $3$.

\midinsert
\*
\insertplotttt{400pt}{150pt}{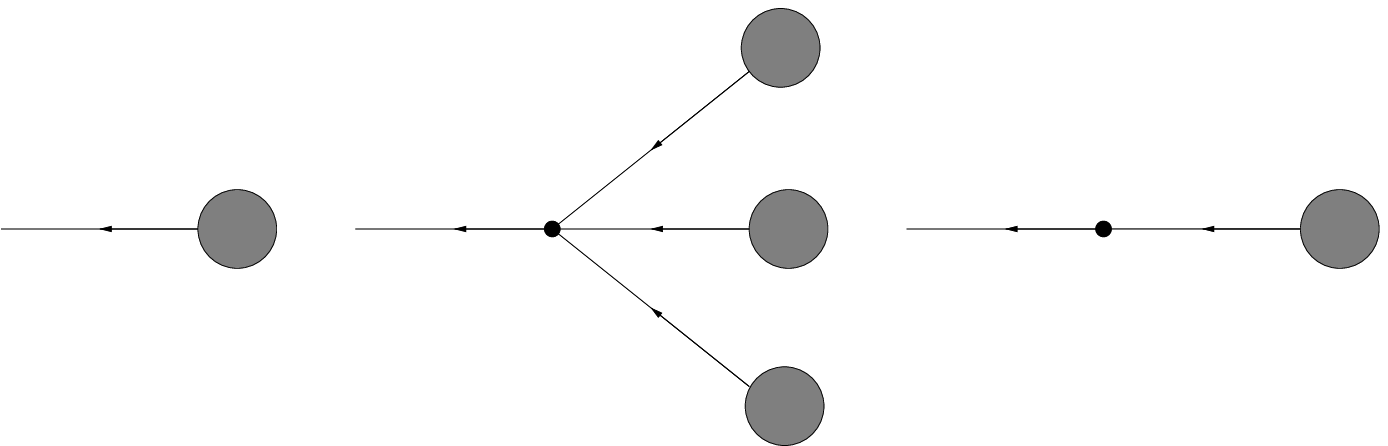}
\ins{078pt}{103pt}{$(k,n,m)$}
\ins{115pt}{077.5pt}{$=$}
\ins{122pt}{070pt}{$(h,i,n,m,m')$}
\ins{234pt}{152pt}{$(k_{1},n_{1},m_{1})$}
\ins{234pt}{100pt}{$(k_{2},n_{2},m_{2})$}
\ins{234pt}{050pt}{$(k_{3},n_{3},m_{3})$}
\ins{282pt}{070pt}{$(h,1,n,m,m')$}
\ins{320pt}{90pt}{$(k_{v},n,m',m_{1})$}
\ins{390pt}{100pt}{$(k_{1},n,m_{1})$}
\ins{276pt}{079pt}{$+$}
\line{\vtop{\line{\hskip1.2truecm\vbox{\advance\hsize by -2.5 truecm
\noindent{\nota {\bf Figure 5.}
Graphical representation of \equ(3.19);
the sums are understood; note that $\sum_{j} \s(\ell_{j}) m_{j}=m'$
in the first summand and $k_{v}+k_{1}=k$ in the second summand.\vfil
}} \hfill} }}
\*
\endinsert

In the same way consider a tree $\th\in \Th_{n,m}^\ka$ such that
$m=\L_{j}(a)$, $s_{v_{0}}=1$ and $h_{\ell_{0}}=h$, with the same
notations as before. Let $\th_{1}$ be the sub-tree whose
root line $\ell_{1}$ enters $v_{0}$. Set $k_{v_{0}}=r$,
$m_{v_{0}}=\L_{j}(b)$, $m'_{v_{0}}=\L_{j}(b')$,
where $b=b_{\ell_{0}}$ and $b'=a_{\ell_{1}}$. Then
$$ \Val(\th)=  G_{n,j,h,1}(a,b) \, L^{(r)}_{n,j,h}(b,b')
\, \Val(\th_{1}) ,
\Eq(3.23) $$
so that the proof is complete. \qed

\*\*
\0{\bf 3.4. Clusters and resonances}\*

\0In the preceding section we have found a power series expansion
for $U_{n,j}$ solving \equ(2.9) and parameterised by
$L_{n,j}$. However for general values of $L_{n,j}$ 
such expansion is not convergent, as one can easily identify
contributions at order $k$ which are $O(k!^\xi)$,
for a suitable constant $\xi$. In this section we show that it is
possible to choose the parameters $L_{n,j}$ in a proper way to
cancel such ``dangerous'' contributions; in order to do this
we have to identify the dangerous contributions and this will be
done through the notion of {\it clusters} and {\it resonances}.

\*
\0\sub(D8)
{\it Given a tree $\th\in\Th^{(k)}_{n,m}$ a {\rm cluster} $T$
on scale $h$ is a connected maximal set of nodes and lines such that
all the lines $\ell$ have a scale label $\le h$ and at least
one of them has scale $h$; we shall call $h_{T}=h$ the scale of
the cluster. We shall denote by $V(T)$, $V_{0}(T)$ and $E(T)$
the set of nodes, internal nodes and the set of end-points,
respectively, which are contained inside the cluster $T$,
and with $L(T)$ the set of lines connecting them.
Finally $k_{T}=\sum_{V(T)}k_{v}$ will be called the order of $T$.}
\*

Therefore an inclusion relation is established between clusters,
in such a way that the innermost clusters are the clusters
with lowest scale, and so on. A cluster $T$ can have an arbitrary
number of lines entering it ({\it entering lines}), but only one
or zero line coming out from it ({\it exiting line} or {\it root line}
of the cluster); we shall denote the latter (when it exists)
with $\ell_{T}^{1}$. Notice that by definition
all the external lines have $i_{\ell}=1$. 

\*
\0\sub(D9)
{\it We call {\rm 1-resonance} on scale $h$ a cluster $T$ of
scale $h_{T}=h$ with only one entering line $\ell_{T}$ and
one exiting line $\ell_{T}^{1}$ of scale $h_{T}^{(e)}>h+1$,
with $|V(T)|>1$ and such that\\
(i) one has
%
$$ n_{\ell_{T}^{1}} =  n_{\ell_{T}}
\ge 2^{(h-2)/\t} ,\qquad m_{\ell_{T}^{1}}' \in \L_{j_{\ell_T}} ,
\Eq(3.24) , $$
(ii) if for some $\ell\in L(T)$ not on the path $\PPP(\ell_{T},
\ell_{T}^{1})$ one has $n_{\ell}=n_{\ell_{T}}$, then
$j_{\ell}\neq j_{\ell_{T}}$.\\
We call {\rm 2-resonance} a set of lines and nodes which
can be obtained from a 1-resonance by setting $i_{\ell_{T}}=0$.\\
Finally we call {\rm resonances} the 1- and 2-resonances.
The line $\ell_{T}^{1}$ of a resonance
will be called the root line of the resonance.
The root lines of the resonances will be also called resonant lines.}
\*

\0{\bf Remarks.}
(1) A 2-resonance is not a cluster, but it is well defined due
to condition (ii) of the 1-resonances. Indeed, such a condition
implies that there is a one to one correspondence between
1-resonances and 2-resonances.\\
(2) The reason why we do not include in the definition of 1-resonances
the clusters which satisfy only condition (i), i.e. such that there
is a line $\ell\in L(T)\setminus \PPP(\ell_{T},\ell_{T}^{1})$
with $n_{\ell}=n_{\ell_{T}}$ and $j_{\ell} = j_{\ell_{T}}$,
is that these clusters do not give any problems and can be easily
controlled, as will become clear in the proof of Lemma \secc(L10);
cf. also the subsequent Remark (1).\\
(3) The 2-resonances are included among the resonances
for the following reason. The 1-resonances are the dangerous
contributions, and we shall cancel them by a suitable choice of the
counterterms. Such a choice automatically cancels out the 2-resonances.

\*

An example of resonance is illustrated in Figure 6.
We associate a numerical value with the resonances as done
for the trees. To do this we need some further notations.

\midinsert
\*
\insertplotttt{220pt}{175pt}{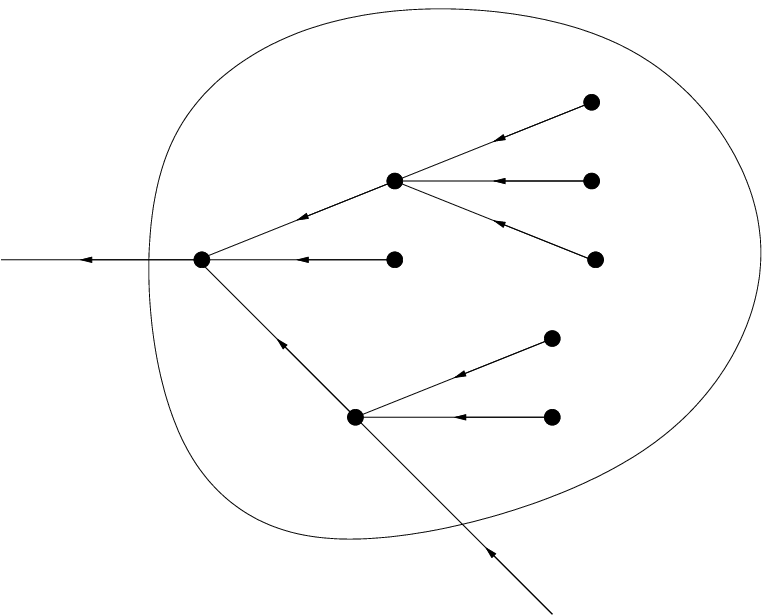}
\ins{106pt}{110pt}{$n,j,m',h_{\ell_{T}^{1}}$}
\ins{138pt}{134pt}{$\ell_{T}^{1}$}
\ins{200pt}{142pt}{$\ell_{1}$}
\ins{190pt}{90pt}{$\ell_{2}$}
\ins{220pt}{024pt}{$n,j,m,h_{\ell_{T}}$}
\ins{268pt}{036pt}{$\ell_{T}$}
\line{\vtop{\line{\hskip1.2truecm\vbox{\advance\hsize by -2.5 truecm
\noindent{\nota {\bf Figure 6.}
Example of resonance $T$. We have set $j_{\ell_{T}^{1}}=j$,
$n_{\ell_{T}^{1}}=n$, $m_{\ell_{T}^{1}}'=m'$,
$m_{\ell_{T}}=m$, so that $n_{\ell_{T}}=n$ and $j_{\ell_{T}}=j$,
by \equ(3.24). Moreover, if $h_{T}=h$ is the scale of $T$,
one has $h_{\ell_{T}} \ge h+1$ by definition of cluster and
$h_{\ell_{T}^{1}} = h^{(e)}_{T} >h+1$ by definition of resonance.
For any line $\ell\in L(T)$ one
has $h_{\ell}\le h$ and there is at least one line on scale $h$.
The path $\PPP(\ell_{T},\ell_{T}^{1})$ consists of
the line $\ell_{1}$. If $n_{\ell_{2}}=n$ then
$j_{\ell_{2}} \neq j$ by the condition (ii).\vfil
}} \hfill} }}
\*
\endinsert

\*
\0\sub(D10)
{\it The trees $\th\in \R^\ka_{h,n,j}$
with $n \ge 2^{(h-2)/ \t}$ and $(n,j)\in \O$  are defined as the
trees $\th\in \Th^\ka_{h,n,m}$ with the following modifications:\\
(a) there is a single end-point, called $e$, carrying the labels
$n_{e},m_{e}$ such that $n_{e} = n$, $m_{e}\in \L_{j}$; if $\ell_{e}$
is the line exiting from $e$ then we associate with it a propagator
$g_{\ell_{e}}=1$, a label $m_{\ell_{e}}=m_{e}$
and a label $\s_{\ell_{e}}\in\{\pm1\}$;\\
(b) the root line $\ell_{0}$ has $i_{\ell_{0}}=1$, $n_{\ell_{0}}=n$
and $m_{\ell_{0}}'\in \L_{j}$ and the corresponding propagator is
$g_{\ell_{0}}=1$;\\
(c) one has $\max_{\ell\in L(\th)\setminus
\{\ell_{0},\ell_{e})\}}h_{\ell}=h$.\\
A cluster $T$ (and consequently a resonance) on scale $h_{T} \le h$ for $\th\in \R^\ka_{h,n,j}$
is defined as a connected maximal set of nodes $v\in V(\th)$ and lines
$\ell\in L(\th) \setminus \{\ell_{0},\ell_{e}\}$ such that
all the lines $\ell$ have a scale label $\le h_{T}$ and at least
one of them has scale $h_{T}$. \\
We define the set $\R^{(k)}$ as the set of trees belonging
to $\R^{\ka}_{h,n,j}$ for some triple $(h,n,j)$.}
\*

\0{\bf Remark.} The entering line $\ell_{e}$ has no label $m_{\ell_{e}}'$,
while the root line has no label $m_{\ell_{0}}$.
Both carry no scale label. Recall that by the diagrammatic
rule (7) the root line $\ell_{0}$ has no $\s$ label.

\*
\0\lm(L7) {\it Let $B$ be the same constant as in Lemma \secc(L5).
For all $\th\in \RR^\ka_{h,n,j}$ and for all $\ell$ not in the
path $\PPP(\ell_{e},\ell_{0})$ one has $|n_{\ell}| \le B k$
and $| m_{\ell}|, |m'_{\ell}| \le B k$.
For $\ell$ on such path one has $\min\{|n_{\ell} - n_{e}|,
|n_{\ell} + n_{e}|\}\le B k$.}
\*

\0{\bf Proof.} For the lines not along the path $\PPP=\PPP(\ell_{e},
\ell_{0})$ the proof is as for Lemma \secc(L5). If a line $\ell$ is
along the path $\PPP$ then one can write $n_{\ell}= n_{\ell}^{0} \pm
n_{e}$, where $n_{\ell}^{0}$ is the sum of the labels $\pm n_{v}$
of all the end-points preceding $\ell$ but $e$. The signs depend
on the labels $\s(\ell')$ of the lines $\ell'$ preceding $\ell$;
in particular the sign in front of $n_{e}$ depends on the
labels $\s(\ell')$ of the lines $\ell'\in\PPP(\ell_{e},\ell)$,
in agreement with to \equ(3.15). Then the last assertion follows
by reasoning once more as in the proof of Lemma \secc(L5). \qed

\*

The definition of value of the trees in $\RR^{(k)}$
is identical to that given in \equ(3.20) for the trees
in $\Th^{(k)}$.

\midinsert
\*
\insertplotttt{339pt}{138pt}{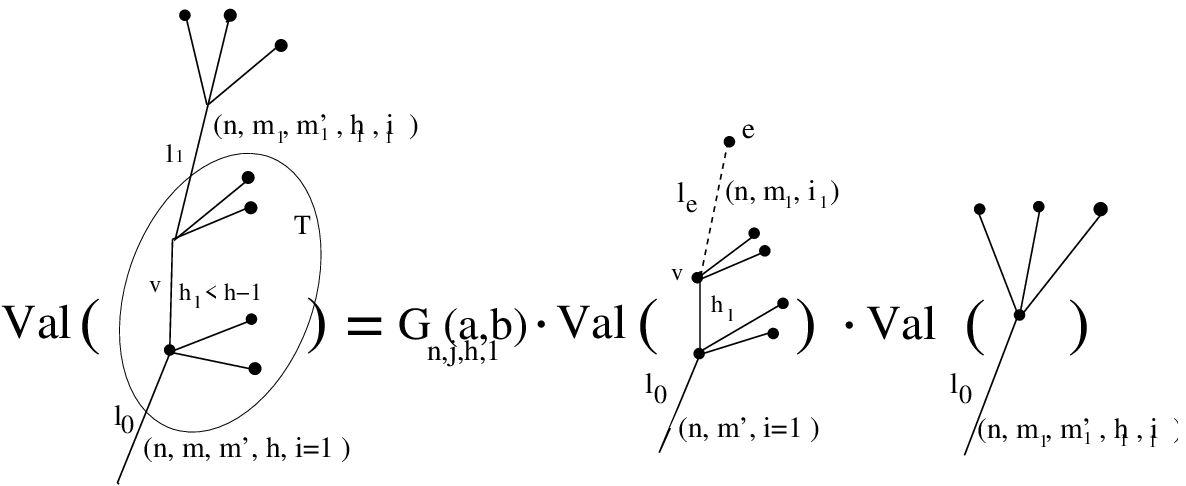}
\medskip
\line{\vtop{\line{\hskip1.5truecm\vbox{\advance\hsize by -3.1 truecm
\*
\noindent{\nota {\bf Figure 7}. We associate with the resonance $T$
(enclosed in an ellipse and such that $m=\L_{j}(a),m'=\L_{j}(b),
m_{1},m_{1}'\in \L_{j}$) the tree $\th_T\in \R_{h_{1},n,j}$,
and vice-versa.\vfil
}} \hfill} }}
\*
\endinsert

Let us now consider a tree $\th$ with a resonance $T$ whose exiting
line is the root line $\ell_{0}$ of $\th$, let $\th_{1}$ be the
tree atop the resonance.
Given a resonance $T$, there exists a 
unique $\th_T\in\R^\ka_{h,n,j}$, with $n=n_{\ell_{0}}$,
$j=j_{\ell_{0}}$ and $h=h_{T}$, such that (see Figure 7)
$$ \Val(\th)= g_{\ell_{0}} \, \Val(\th_T) \, \Val(\th_{1}) ,
\Eq(3.25) $$ 
so that we can call, with a slight abuse of language, $\Val(\th_T)$
the value of the resonance $T$.

\*\*
\0{\bf 3.5. Choice of the parameters $\LLLL_{\nnnnn,\jjjj}$}\*

\0With a suitable choice of the parameters $L_{n,j,h}$ the
functions $u_{n,m}^{(k)}$ can be rewritten as sum over 
``renormalised'' trees defined below.

\*
\0\sub(D11)
{\it We define the set of {\rm renormalised trees}
$\Th_{R,n,m}^{(k)}$ defined as the trees in 
$\Th_{n,m}^{(k)}$  with no resonances nor nodes with $s_v=1$.
In the same way we define $\R_{R,h,n,j}^{(k)}$.
We call $\R_{R,h,n,j}^{(k)}(a,b) $ the set of trees
$\th\in\R_{R,h,n,j}^{(k)}$ such that the entering line has
$m_{e}= \L_{j}(b)$ while the root line has $m_{\ell_{0}}'= \L_{j}(a)$.
Finally we define the sets $\Th_{R}^{(k)}$ and $\R_{R}^{(k)}$
as the sets of trees belonging to $\Th_{R,n,m}^{(k)}$ for some $n,m$
and, respectively, to $\R_{R,h,n,j}^{(k)}$ for some $h,n,j$.}
\*

We extend the notion of resonant line by including also the lines
coming out from a node $v$ with $s_v=1$. This leads to the
following definition.

\*
\0\sub(D12)
{\it A resonant line is either the root line of a resonance
(see Definition \secc(D9)) or the
line exiting a node $v$ with $s_{v}=1$.}
\*

The following result holds.

\*
\0\lm(L8)
{\it For all $k,n,m$ one has
$$ u_{n,m}^{(k)}= \sum_{\th\in \Th_{R,n,m}^{(k)}}\Val(\th) ,
\Eq(3.26) $$
provided we choose in \equ(3.13)
$$ \left\{ \eqalign{
L_{n,j,h}^{(k)}(a,b) & =- \sum_{h_{1}< h-1}
\sum_{\th\in \R_{R,h_{1},n,j}^{(k)}
(a,b)}\Val(\th) , \qquad (n,j) \in \O , \cr
L_{n,j,h}^{(k)}(a,b) & = 0 , \qquad (n,j) \notin \O , \cr} \right.
\Eq(3.27) $$ 
where $\R_{R,h_{1},n,j}^{(k)}(a,b) $ is
as in Definition \secc(D11).}
\*

\0{\bf Proof.} First note that by definition
$L_{n,j,h}= 0$ if $(n,j)\notin\O$.
We proceed by induction on $k$. For $k=1$ \equ(3.27) holds as
$\Th_{R,n,m}^{(1)}\equiv\Th_{n,m}^{(1)}$. Then we assume that
\equ(3.27) holds for all $r<k$. By \equ(3.12) one has
$U_{n,j,h,i}^{(k)}=  G_{n,j,h,i} F_{n,j}^\ka$for $i=-1,0$, and
$$ U_{n,j,h,1}^{(k)}=  G_{n,j,h,1} F_{n,j}^\ka +
G_{n,j,h,1}\left(\sum_{h_{2}=-1}^{\io} \sum_{i_{2}=1,0}\sum_{r=1}^{k-1}
L^{(r)}_{n,j,h} \, U^{(k-r)}_{n,j,h_{2},i_{2}} \right) ,
\Eq(3.28) $$ 
where $F^\ka_{n,j}$ is a function of the coefficients
$u_{n',m'}^{(r')}$ with $r'<k$.
By the inductive hypothesis each $u_{n',m'}^{(r')}$
can be expressed as a sum over trees in $\Th^{(r')}_{R,n',m'}$.
Therefore $(G_{n,j,h,i} F_{n,j}^\ka)(a)$ is given by the sum over 
the trees $\th\in \Th^{(k)}_{n,m}$, with $m=\L_{j}(a)$
and $s_{v_{0}}=3$ ($v_{0}$ is introduced in Definition \secc(D6)),
such that only the root line $\ell_{0}$ of $\th$ can be resonant.
Note that $\ell_{0}$ can be resonant only if $i=i_{\ell_{0}}=1$.
If $\ell_{0}$ is non-resonant then $\th\in \Th_{R,n,m}^\ka$,
so that the assertion holds trivially for $i\neq1$.

For $i=1$ we split the coefficients of $G_{n,j,h,1} F_{n,j}^\ka$
as sum of two terms: the first one, denoted $G_{n,j,h,1} J_{n,j}^{(k)}$,
is the sum over all trees belonging to $\Th_{R,n,m}$ for $m\in\L_{j}$
with $s_{v_{0}}=3$ and the second one is sum of trees with value
$$\Val(\th)=g_{\ell_{0}} \, \Val(\th_T) \, \Val(\th_{1}) ,
\Eq(3.29) $$
with $\th_T\in \R^{(r)}_{R,h_{1},n,j}$ and  $\th_{1}\in
\Th^{(k-r)}_{R,n,m'}$ with $m'=\L_{j}(b)$ for some $r$ and some $b$;
by definition of resonance we have  $h_{1}<h-1$.

We get terms of this type for all
$\th_T$ and $\th_{1}$ so that
$$ F_{n,j}^\ka(a) = J_{n,j}^\ka(a) +
\sum_{b=1}^{d_{j}} \sum_{h_{2}=-1}^{\io}\sum_{i_{2}=1,0}
\sum_{r=1}^{k-1}\sum_{h_{1}<h-1} \left(
\sum_{\th\in \R_{R,h_{1},n,j}^{(r)}(a,b)}
\Val(\th) \right) U_{n,j,h_{2},i_{2}}^{(k-r)}(b) ,
\Eq(3.30)$$ 
where the sum over $h_{1}<h-1$ of the terms between parentheses
gives $-L^{(r)}_{n,j,h}(a,b)$ by the first line in \equ(3.27)
Therefore all the terms but $J_{n,j}^{\ka}(a)$ in \equ(3.30)
cancel out the term between parentheses in \equ(3.28), and only
the term $G_{n,j,h,i} J_{n,j}^{(k)}(a)$ is left in \equ(3.28).
On the other hand $G_{n,j,h,i} J_{n,j}^{(k)}(a)$ is by definition
the sum over all trees in $\Th_{R,n,m}^\ka$, so that
the assertion follows also for $i=1$. \qed

\*

\0{\bf Remarks.} (1) The proof of Lemma \secc(L8) justifies why
we included into the definition of resonances (cf. Definition \secc(D9))
also the 2-resonances, even if the latter are not clusters.
Indeed in \equ(3.28) we have to sum also over $i_{2}=0$.\\ 
(2) Note that $\Val(\th)$ is a monomial of degree
$2k+1$ in $q$ for $\th\in\Th_{R,n.m}^{(k)}$, and it is a monomial
of degree $2k$ in $q$ for $\th\in\Th_{R,n,m}^{(k)}$.

\*

In the next section we shall prove that the matrices $L_{n,j,h}^\ka$
are symmetric (we still have to show that the matrices are
well-defined). For this we shall need the following result.

\*
\0\lm(L9) {\it For all trees $\th\in \RR_{R,h,n,j}(a,b)$ there exists a
tree $\th_{1}\in \RR_{R,h,n,j}(b,a)$ such that $\Val(\th)=\Val(\th_{1})$.}
\*

\0{\bf Proof.} Given a tree $\th\in \RR_{R,h,n,j}(a,b)$ consider the
path $\PPP=\PPP(\ell_{e},\ell_{0})$, and set $\PPP=\{\ell_{1},
\ldots,\ell_{N}\}$, with $\ell_{0} \succ \ell_{1} \succ \ldots \succ
\ell_{N} \succ \ell_{N+1}=\ell_{e}$. We construct a tree $\th_{1} \in
\R_{R,h,n,j}(b,a)$ in the following way.
\vskip.2truecm
\01. We shift the $\s_{\ell}$ labels down the path $\PPP$,
so that $\s_{\ell_{k}} \to \s_{\ell_{k+1}}$ for $k=1,\ldots,N$,
$\ell_{0}$ acquires the label $\s_{\ell_{1}}$, while
$\ell_{e}$ loses its label $\s_{\ell_{e}}$ (which becomes
associated with the line $\ell_{N}$).
\vskip.1truecm
\02. For all the lines $\ell\in\PPP$ we exchange the labels
$m_{\ell},m'_{\ell}$, so that $m_{\ell_{k}} \to m_{\ell_{k}}'$,
$m_{\ell_{k}}' \to m_{\ell_{k}}$ for $k=1,\ldots,N$,
while one has simply $m_{\ell_{0}}' \to m_{\ell_{e}}$ and
$m_{\ell_{e}} \to m_{\ell_{0}}'$ for the root and entering lines.
\vskip.1truecm
\03. For any pair $\ell_{1}(v),\ell_{2}(v)$ of lines not on the path
$\PPP$ and entering the node $v$ along the path,
we exchange the corresponding labels $\s_{\ell}$,
i.e. $\s_{\ell_{1}(v)} \to \s_{\ell_{2}(v)}$ and
$\s_{\ell_{2}(v)} \to \s_{\ell_{1}(v)}$.
\vskip.1truecm
\04. The line $\ell_{e}$ becomes the root line,
and the line $\ell_{0}$ becomes the entering line.
\vskip.2truecm

As a consequence of item 4. the ordering of nodes and lines along
the path $\PPP$ is reversed (in particular the arrows of all the
lines $\ell\in\PPP\cup\{\ell_{T},\ell_{T}^{1}\}$ are reversed).
On the contrary the ordering of all the lines and nodes outside
$\PPP$ is not changed by the operations above. This means that
all propagators and node factors of lines and nodes, respectively,
which do not belong to $\PPP$ remain the same. 

Then the symmetry of $M$, hence of the propagators,
implies the result. \qed

\*\*
\section(4,Bryuno lemma and bounds)

\0In the previous section we have shown that, with a suitable choice
of the parameters $L_{n,j}$, we can express the coefficients
$u^{(k)}_{n,m}$ as sums over trees belonging to $\Th^{(k)}_{R,n,m}$.
We show in this section that such expansion is indeed convergent if
$\h$ is small enough and $(\e,M)\in \DD(\g)$ (see Definition \secc(D3)).

\*\*
\0{\bf 4.1. Bounds on the trees in $\TTTT_\RRRR^{(\kkkk)}$}\*

\0Given a tree $\th\in\Th^{(k)}_{R,n,m} $, we call  $\SS(\th,\g)$
the set of  $(\e,M)\in \DD_{0}$ such that for all
$\ell\in L(\th)$ one has
$$ \left\{
\eqalign{
2^{-h_{\ell}-1}\g \le \; & |x_{n_{\ell},j_{\ell}}| \;
\le \; 2^{-h_{\ell}+1}\g ,
\hskip1.truecm h_{\ell}\neq -1 , \cr
& |x_{n_{\ell},j_{\ell}}| \; \ge \; \g ,
\hskip2.1truecm h_{\ell}= -1 , \cr} \right.
\Eq(4.1) $$
with $x_{n,j}$ defined in \equ(2.17), and
$$ \left\{
\eqalign{
& |y_{n_{\ell},j_{\ell}}| \; \le \; 2^{-2}\g ,
\hskip1.2truecm i_{\ell}=1 , \cr
2^{-3}\g \; \le \; & |y_{n_{\ell},j_{\ell}}| \; \le 
\; 2^{-1}\g ,
\hskip1.2truecm i_{\ell}=0 , \cr
2^{-2} \g \; \le \; & |y_{n_{\ell},j_{\ell}}| ,
\hskip2.6truecm i_{\ell}= -1 . \cr} \right.
\Eq(4.2) $$
with $y_{n,j}$ defined in \equ(2.8).
In other words we can have $\Val(\th)\neq 0$
only if $(\e,M)\in \SS(\th,\g)$.

We call $\DD(\th,\g)\subset \DD_{0}$ the set of $(\e,M)$ such that
$$ \left| x_{n_{\ell},j_{\ell}} \right|
\ge {\g \over |n_{\ell}|^{\t} }
\Eq(4.3) $$
for all lines $\ell\in L(\th)$ such that $i_{\ell}=1$, and
$$ \left| \d_{n_{\ell},j_{\ell}} - \d_{n_{\ell_{1}},j_{\ell_{1}}}
\right| \ge {\g \over |n_{\ell}- n_{\ell_{1}}|^{\t}}
\Eq(4.4) $$
for all pairs of  lines $\ell_{1}\prec \ell\in L(\th)$ such
that $n_{\ell}\neq n_{\ell_{1}}$, $i_{\ell},i_{\ell_{1}}=0,1$
and $\prod_{\ell'\in \PPP(\ell_{1},\ell)}
\s(\ell')\s(\ell_{1})=1$ (the last condition implies
that $|n_{\ell}-n_{\ell_{1}}|$ is bounded by the sum of
$|n_{v}|$ of the nodes $v$ preceding $\ell$ but not $\ell_{1}$).
This means that $\DD(\th,\g)$ is the set of $(\e,M)$ verifying
the Melnikov conditions \equ(2.2) and \equ(2.18) in $\th$.

In order to bound $\Val(\th)$ we will use the following result
(Bryuno lemma).

\*
\0\noindent\lm(L10)
{\it Given a tree $\th \in \Th^{(k)}_{R,n,m}$
such that  $\DD(\th,\g)\cap \SS(\th,\g)\neq \emptyset$, then the scales 
$h_{\ell}$ of $\th$ obey
$$ N_{h}(\th) \le  \max\{ 0 , c \, k(\th) 2^{(2-h)\b/\t} - 1 \} , 
\Eq(4.5) $$
where $ N_{h}(\th) $ is the number of lines $\ell$
with $i_{\ell}=1$ and scale $h_{\ell}$ greater or equal than $h$,
and $c$ is a suitable constant.}
\*

\0{\bf Proof.}
For $(\e,M)\in \DD(\th,\g)\cap \SS(\th,\g)$ both \equ(4.1)
and \equ(4.3) hold. Moreover by Lemma \secc(L5) one has
$|n| \le B k(\th)$. This implies  that one can have $N_{h}(\th) \ge 1$
only if $k(\th)$ is such that $k(\th) > k_{0} := B^{-1} 2^{(h-1)/\t}$.
Therefore for values $k(\th) \le k_{0}$ the bound \equ(4.5) is satisfied.

If $k(\th)>k_{0}$, we proceed by induction by assuming that
the bound holds for all trees $\th'$ with $k(\th')<k(\th)$.
Define $E_{h}:=c^{-1}2^{(-2+h)\b/\t}$: so we have
to prove that $N_{h}(\th)\le \max\{0,k(\th) E_{h}^{-1}-1\}$.
In the following we shall assume that $c$ is so large
that all the assertions we shall make hold true.

Call $\ell_{0}$ the root line of $\th$ and $\ell_{1},\ldots,\ell_{m}$
the $m\ge 0$ lines on scale $\ge h-1$ which are the closest
to $\ell_{0}$ and such that $i_{\ell_s}=1$ for $s=1,\dots,m$.

If the root line $\ell_{0}$ of $\th$ is on scale $< h$ then
$$ N_{h}(\th) = \sum_{i=1}^{m} N_{h}(\th_{i}) ,
\Eq(4.6) $$
where $\th_{i}$ is the sub-tree with $\ell_{i}$ as root line.

By construction $N_{h-1}(\th_{i})\ge 1$, so that $k(\th_{i})> B^{-1}
2^{(h-2)/\t}$ and therefore for $c$ large enough
(recall that $\b < \a \ll 1$) one has
$\max\{0,k(\th_{i}) E_{h}^{-1}-1\} =k(\th_{i}) E_{h}^{-1}-1$,
and the bound follows by the inductive hypothesis.

If the root line $\ell_{0}$ has $i_{\ell_{0}}=1$ and scale $\geq h$
then $\ell_{1},\ldots,\ell_{m}$ are the entering line of a cluster $T$.

By denoting again with $\th_{i}$ the sub-tree
having $\ell_{i}$ as root line, one has
$$ N_{h}(\th) = 1 + \sum_{i=1}^{m} N_{h}(\th_{i}) ,
\Eq(4.7) $$
so that, by the inductive assumption,
the bound becomes trivial if either $m=0$ or $m\ge 2$.

If $m=1$ then one has a cluster $T$ with two external lines
$\ell_{T}^{1}=\ell_{0}$ and $\ell_{T}=\ell_{1}$, such that
$h_{\ell_{1}}\ge h-1$ and $h_{\ell_{0}}\ge h $.
Then, for the assertion to hold in such a case, we have to prove
that $(k(\th)-k(\th_{1}))E_{h}^{-1}\ge 1$.
For $(\e,M)\in \SS(\th,\g) \cap \DD(\th,\g)$ one has
$$ \min \{ |n_{\ell_{0}}|,|n_{\ell_{1}}| \} |\ge 2^{(h-2)/\t} ,
\Eq(4.8) $$
and, by definition, one has $i_{\ell_{0}}=i_{\ell_{1}}=1$,
hence $|y_{n_{\ell_{0}},j_{\ell_{0}}}|,
|y_{n_{\ell_{1}},j_{\ell_{1}}}| \le \g/4$ (see \equ(4.2)),
so that we can apply Lemma \secc(L1).

We distinguish between two cases.
\vskip.2truecm
\01. If $n_{\ell_{0}} \not = n_{\ell_{1}}$, by Lemma \secc(L1)
with $s_{0}=s_{2}$ (and the subsequent Remark) one has
$$ \left| n_{\ell_{0}} \pm  n_{\ell_{1}} \right| \ge
{\rm const.} \min\{ |n_{\ell_{0}}|,|n_{\ell_{1}}|\}^{s_{2} / \t_{1}}
\ge {\rm const.} \min\{ |n_{\ell_{0}}|,|n_{\ell_{1}}|\}^{s_{2} / \t}
\ge {\rm const.} 2^{(h-2)s_{2}/\t^{2}} \ge E_{h} , $$  
where we have used that $s_{2}/\t^{2} \ge \b/\t$ for $\a$ small
enough. Therefore $B ( k(\th)-k(\th_{1}) ) \geq
\min_{a=\pm 1} |n_{\ell_{0}} +a n_{\ell_{1}}|\geq E_{h}$.
\vskip.2truecm
\02. If $n_{\ell_{0}} = n_{\ell_{1}}$, consider the
path $\PPP=\PPP(\ell_{1},\ell_{0})$. Now consider
the nodes along the path, and call $\ell_{i}$ the lines entering
these nodes and $\th_{i}$ the sub-trees which have such lines
as root lines. If $m_{i}$ denotes the momentum label $m_{\ell_{i}}$
one has, by Lemma \secc(L5), $|m_{i}|\le B k(\th_{i})$.

Call $\bar\ell$ the line on the path $\PPP\cup\{\ell_{1}\}$
closest to $\ell_{0}$ such that $i_{\bar\ell} \neq -1$ (that is all
lines $\ell$ along the path $\PPP(\bar\ell,\ell_{0})$ have $i_{\ell}=-1$).
\vskip.2truecm
\02.1. If $|n_{\bar\ell}|\le |n_{\ell_{0}}|/2$ then,
by the conservation law \equ(3.15) one has
$k(\th)-k(\th_{1})> B^{-1} |n_{\ell_{0}}|/2 \ge E_{h}$.
\vskip.2truecm
\02.2. If $|n_{\bar\ell}|\ge |n_{\ell_{0}}|/2$ we distinguish
between the two following cases.
\vskip.2truecm
\02.2.1. If $n_{\bar\ell} \neq n_{\ell_{0}}$ ($=n_{\ell_{1}}$)
then by Lemma \secc(L1) and \equ(4.2) one finds
$$ |n_{\bar\ell} \pm  n_{\ell_{0}}| \ge
{\rm const.} \min\{ |n_{\bar\ell}|,|n_{\ell_{0}}|\}^{s_{2} / \t}
\ge {\rm const.} 2^{-s_{2}/\t} 2^{(h-2)s_{2}/\t^{2}} > E_{h} .$$
\vskip.2truecm
\02.2.2. If $n_{\bar\ell}= n_{\ell_{0}}$
then we have two further sub-cases.
\vskip.2truecm
\02.2.2.1. If $j_{\ell_{0}} \neq j_{\bar\ell}$, then
$|m_{\bar\ell}-m_{\ell_{0}}'| \ge C_{2} p_{j_{\ell_{0}}}^{\b} \ge
C |n_{\ell_{0}}|^{\b}$, for some constant $C$. For all the lines
$\ell$ along the path $\PPP(\bar\ell,\ell_{0})$
one has $i_{\ell}=-1$, hence $m_{\ell}=m_{\ell}'$ (cf. the
Remark (3) after Definition \secc(D5)), so that
$|m_{\bar\ell}-m_{\ell_{0}}'| \le \sum_{i} |m_{i}| \le B
(k(\th)-k(\th_{1}))$, and the assertion follows
once more by using \equ(4.8).
\vskip.2truecm
\02.2.2.2. If $j_{\ell_{0}} = j_{\bar\ell}$ then $i_{\bar\ell}=0$
because $h_{\bar\ell} \le h-2$ and one would have $|h_{\bar\ell}-h|
=|h_{\bar\ell}-h_{\ell_{0}}|\le 1$ if $i_{\bar\ell}=1$.
As 2-resonances (as well as 1-resonances) are not possible
there exists a line $\ell'$ (again with $i_{\ell'}=0$ because
$h_{\ell'} \le h-2$), not on the path $\PPP(\bar\ell,\ell_{0})$, such
that $j_{\ell'}=j_{\ell_{0}}$ and $|n_{\ell'}|=|n_{\ell_{0}}| >
2^{(h-2)/\t}$; cf. condition (ii) in Definition \secc(D9).
In this case one has $k(\th)-k(\th_{1})> B^{-1} |n_{\ell'}| \ge E_{h}$.
\vskip.2truecm
This completes the proof of the lemma. \qed

\*

\0{\bf Remarks.}
(1) It is just the notion of 2-resonance and property (ii) in
Definition \secc(D9), which makes non-trivial the case 2.2.2.2.
in the proof of Lemma \secc(L10).\\
(2) Note that in the discussion of the case 2.2.2.1. we
have proved that $k(\th)-k(\th_{1}) \ge {\rm const.}
|n_{\ell_{0}}|^{\b}$ (using once more that $s_{2}/\t \ge
\b$ for $\a$, hence $\b$, small enough with respect to $s$).

\*

The Bryuno lemma implies the the following result.

\*
\0\lm(L11) {\it There is a positive constant $D_{0}$ such that
for all trees $\th\in \Th^{k}_{R,n,m}$ and for all 
$(\e,M)\in \DD(\th,\g) \cap \SS(\th,\g)$ one has
$$ \eqalign{
& (i) \quad \left| \Val(\th) \right| \le
D_{0}^{k} q^{2k+1} \Big( \prod_{h=1}^{\io} 2^{h N_{h}(\th)} \Big)
\prod_{\ell \in L(\th)} p_{j_{\ell}}^{-3s/4} , \cr
& (ii) \quad \left| \dpr_{\e} \Val(\th) \right| \le
D_{0}^{k} q^{2k+1} \Big( \prod_{h=1}^{\io} 2^{2 h N_{h}(\th)} \Big)
\prod_{\ell \in L(\th)} p_{j_{\ell}}^{-s_{2}-\a} , \cr
& (iii) \sum_{(n',j')\in \O}\sum_{a',b'=1}^{d_{j'}}
|\dpr_{M_{n',j'}(a',b')}\Val(\th)| \le D_{0}^{k} q^{2k+1} 
\Big( \prod_{h=1}^{\io} 2^{2 h N_{h}(\th)} \Big)
\prod_{\ell \in L(\th)} p_{j_{\ell}}^{-s_{2}-\a} . \cr}
\Eq(4.9) $$
if $|n|<Bk$ and $|m|< B k$, with $B$ given
as in Lemma \secc(L5), and $\Val(\th)=0$ otherwise.}
\*

\0{\bf Proof.}
By Lemma \secc(L5) we know that $\Th^{k}_{R,n,m}$ is empty if $|n|>Bk$
or $|m|> Bk$.  We first extract the factor $q^{2k+1}$ by
noticing that a renormalised tree of order $k$ has $2k+1$ end-points
(cf. the proof of Lemma \secc(L5)).

For $(\e,M)\in \DD(\th,\g) \cap \SS(\th,\g)$ the bounds \equ(4.5) hold.
First of all we bound all propagators $g_{\ell}$ such that $i_{\ell}
=-1,0$ with $16C_{1}^{1/2}\g^{-1}|p_{j_{\ell}}|^{-3s/4}$ according
to \equ(3.10). For the remaining $g_{\ell}$ we use the inequalities
\equ(4.1) due to the scale functions: by Lemma \secc(L3) (ii) one has
$|G_{n,j,h,1}(a,b)| \le \sqrt{d_{j}} p_{j}^{-s} |\d_{n,j} +
p_{j}^{-s_{1}}\n_{n,j}|$, so that we can bound the propagators
$g_{\ell}$ proportionally to $2^{h_{\ell}}|p_{j_{\ell}}|^{-3s/4}$.
This proves the bound (i) in \equ(4.9); notice that the product
over the scale labels is convergent.

When deriving $\Val(\th)$ with respect to $\e$ we get a sum of trees
with a distinguished line, say $\ell$, whose propagator $g_{\ell}$ is
substituted with $\dpr_\e g_{\ell}$ in the tree value.
For simplicity, in the following set
$j=j_{\ell}$, $h_{\ell}=h$ and $n=n_{\ell}$.

Let us first consider the case $i_{\ell}=-1,0$ (so that $g_{\ell}$ is
given by the first line of \equ(3.8)), and recall
Lemma \secc(L3) (ii) and (iii)).
Bounding the derivative $\dpr_\e g_{\ell}$
we obtain, instead of the bound on $g_{\ell}$, a factor
$$ C \g^{-1} {p_{j}^{s_{2}}|n | C_{1}^{1/2} p_{j}^{\a/2}
\over p_{j}^{s}|\d_{n,j} + p_{j}^{-s_{1}} \n_{n,j}|}
\le C C_{1}^{1/2} {16 \over \g^{2}} |n|  p_{j}^{-(s -2s_{2}-\a/2)} ,
\Eq(4.10) $$
arising when the derivative acts on $\bar \chi_{i}(y_{n,j})$ (here
and in the following factors $C\g^{-1}$ is a bound on the
derivative of $\chi$ with respect to its argument), and a factor
$$ { 2 |n| C_{1} p_{j}^{\a} \over
p_{j}^{s} (\d_{n,j} + p_{j}^{-s_{1}}\n_{n,j})^{2}}
\le 2 C_{1} {16^{2} \over \g^{2}} |n|  p_{j}^{-(s -2s_{2}-\a)} ,
\Eq(4.11) $$
arising when the derivative acts on the matrix
$(\d_{n,j} I+ p_{j}^{-s} \bar\chi(y_{n,j})M_{n,j})^{-1}$.

If $i_{\ell}=1$ then the propagator is given by the second line
in \equ(3.8), so that both summands arising from the derivation
of the function $\bar\chi_{1}(y_{n,j})$ and of the matrix
$(\d_{n,j} I+ p_{j}^{-s} \bar\chi(y_{n,j})M_{n,j})^{-1}$ are there,
and they are both bounded proportionally to $p_{j}^{-s_{2}- \a}|n |2^{2h}$
(recall that $s_{2}=(s-2\a)/4$). Moreover (see Lemma
\secc(L3) (iv)) there is also an extra summand containing a factor
$$ 2 C \g^{-1} C_{1}^{7/2} { |n| p_{j}^{7 \a/2} 2^{h+1}
\over p_{j}^{s} | \d_{n,j} + p_{j}^{-s_{1}} \n_{n,j}| }
\le {\rm const.} p_{j}^{-s + 4\a}|n|2^{2h},
\Eq(4.12) $$
arising when the derivative acts on $\chi_{h}(x_{n,j})$. Indeed, by
setting $A=(\d_{n,j}I+p_{j}^{-s}\bar\chi_{1}(y_{n,j})M_{n,j})^{-1}$,
so that $x=\|A\|^{-1}$, one has
$$ \dpr_{\e} x_{n,j} = {1\over d_{j}^{1/2} \|A\|^{3}}
\sum_{i,k,h,l=1}^{d_{j}} A(i,k) \, A(i,h) \, A(l,k)
\dpr_{\e} A^{-1}(h,l) ,
\Eq(4.13) $$
which implies \equ(4.12). For $\a\ll s$ we can bound $s - 4\a$
with $s_{2}+\a$.
 
Finally we can bound each $n=n_{\ell}$ with $Bk$ (see
Lemma \secc(L5)). All the undistinguished lines in the tree
(i.e. all lines $\ell'\neq\ell$ in $L(\th)$) can be bounded
as in item (i). This proves the bound (ii) in \equ(4.9).

The derivative with respect to $M_{n',j'}(a',b')$ gives a sum of
trees with a distinguished line $\ell$ (as in the previous
case (ii)), with the propagator $\dpr_{M_{n',j'}(a',b')} g_{\ell}$
replacing $g_{\ell}$. Notice that $\ell$ must carry the labels $n',j'$.
We have two contributions, one arising from the derivative of
the matrix and the other one (provided $i_{\ell}=1$) arising from
the derivative of the scale function $\chi_{h_{\ell}}$
(there is no contribution analogous to \equ(4.10) because
$y_{n,j}$ does not depend on $M$). By
reasoning as in the case (ii) we obtain a factor proportional
to $2^{2h}p_{j_{\ell}}^{-s_{2}-\a}$.

The sums over the labels $(n',j')\in\O$ and $a',b'=1,\ldots,d_{j'}$
can be bounded as follows. By Lemma \secc(L5) one has
$|n'| < Bk$. Then $j'$ must be such that $p_{j'}=O(n')$,
which implies that the number of values which $j'$ can assume
is at most proportional to $|n'|^{D-1}$,
and $a',b'$ vary in $\{1,\ldots,d_{j'}\}$, with $d_{j'}\le
C_{1} p_{j'}^{\a} \le {\rm const.} |n'|^{\a}$. Therefore
we obtain an overall factor proportional to $k^{1+(D-1)+2\a} \le
k^{1+D} \le C^{k}$ for some constant $C$.
Hence also the bound (iii) of \equ(4.9) is proved. \qed

\*\*
\0{\bf 4.2. Bounds on the trees in $\RR^{(\kkkk)}_{\RRRR}$}\*

\0Given a tree $\th\in \R_{R,h,n,j}$, we call 
$\widetilde{\SS}(\th,\g)$ set of  $(\e,M)\in \DD_{0}$ such that \equ(4.1)
holds for all $\ell\in L(\th)\setminus\{\ell_{e},\ell_{0}\}$,
and \equ(4.2) holds for all $\ell\in L(\th)$.
Let $\widetilde{\DD}(\th,\g)\subset \DD_{0}$ be the set of 
$(\e,M)$ such that \equ(4.3) holds for all $\ell\in
L(\th)\setminus\{\ell_{e},\ell_{0}\}$, and \equ(4.4) holds
for all pairs $\ell_{1} \prec \ell\in L(\th)$
such that
\vskip.1truecm
\0(i) $n_{\ell_{1}}\neq n_{\ell}$, $i_{\ell},i_{\ell_{1}}=0,1$
and $\prod_{\ell'\in \PPP(\ell_{1},\ell)}\s(\ell')\s(\ell_{1})=1 $;
\vskip.1truecm
\0(ii) either both $\ell,\ell_{1}$ are on the path $\PPP(\ell_{e},
\ell_{0})$ or none of them is on such a path.

\*

The following lemma will be useful.

\*
\0\lm(L12) {\it Given a tree $\th \in \RR^{(k)}_{R,\bar h,n,j}(a,b)$
such that $\widetilde{\DD}(\th,\g)\cap \widetilde{\SS}(\th,\g)
\neq \emptyset$ then there are two positive constants
$B_{2}$ and $B_{3}$ such that\\
(i)  a line $\ell$ on the path $\PPP(\ell_{e},\ell_{0})$
can have $i_{\ell}\neq -1$ only if $k\ge B_{2} |n|^{\b}$;\\
(ii) one has $k \ge B_{3} |m_a-m_b|^{\r}$ with
$1/\r= 1+\a/\b=1+D(1+D(D+2)!/2)$.}
\*

\0{\bf Proof.} (i) One can proceed very closely to case 2.
in the proof of Lemma \secc(L10), with $\ell_{e}$ and $\ell$ playing
the role of $\ell_{1}$ and $\bar\ell$, respectively -- see the
Remark (2) after the proof of Lemma \secc(L10). We omit the details.\\
\0(ii) By Lemma \secc(L2), for all $m_{a},m_{b}\in \L_{j}$ one
has $|m_{a}-m_{b}| \le C_{2} p_{j}^{\a+\b}$.
For $(n,j)\in \O$ this implies that $|m_a-m_b|\le {\rm const.}\,
|n|^{\a+\b} $. If $k\ge B_{2} |n|^{\b}$ then one has
$k \ge {\rm const.}\,|m_a-m_b|^{\b/(\a+\b)}$, the
statement holds true (recall that $\a/\b$ is given by \equ(2.4)).
If $k < B_{2} |n|^{\b}$ then by
item (i) all the lines $\ell$ on the path $\PPP(\ell_{e}, \ell_{0})$
have $i_{\ell}=-1$, hence $m_{\ell}=m'_{\ell}$. Then
by calling, as in the proof of Lemma \secc(L10),
$\th_{i}$ the sub-trees whose root lines enter the
nodes of $\PPP(\ell_{e},\ell_{0})$ and $m_{i}$ the
momentum label $m_{\ell_{i}}$, we obtain $| m_{\ell_{e}}-
m_{\ell_{0}}|\le \sum_{i} |m_{i}| \le B k$,
and the assertion follows once more. \qed

\*

The following generalisation of Lemma \secc(L10) holds.

\*
\0\lm(L13) {\it
Given  tree $\th \in \RR^{(k)}_{R,\bar h,n,j}$ such that
$\widetilde{\DD}(\th,\g)\cap \widetilde{\SS}(\th,\g)\neq \emptyset$
then the scales $h_{\ell}$ of $\th$ obey, for all  $h \le \bar h$,
$$ N_{h}(\th) \le  \max\{ 0 , c\,k(\th)\,2^{(2-h)\b/\t} \} ,
\Eq(4.14) $$
where $N_{h}(\th)$ and $c$ are defined as in Lemma \secc(L10).}
\*

\0{\bf Proof.} To prove the lemma we consider a slightly different
class of trees with respect to $\RR^{(k)}_{R,\bar h,n,j}$, which we
denote by $\RR^{(k)}_{R,\bar h}$. The differences are as follows:
\vskip.1truecm
\0(i) the root line has scale labels $h_{\ell_{0}} \le \bar h$
and $i_{\ell_{0}}\in\{-1,0,1\}$,
\vskip.1truecm
\0(ii) we remove the condition $n_{e}=n_{\ell_{0}}$,
$j_{e}=j_{\ell_{0}}$, and require only that  $|n_{e}|>2^{(\bar h-2)/\t}$.
\vskip.1truecm

Notice that, for all $\th\in \RR^{(k)}_{R,\bar h,n,j}$, among
the three sub-trees entering the root, two are in $\Th_R^{(k_{1})}$
and $\Th_R^{(k_{2})}$, respectively, and one is in
$\RR^{(k_{3})}_{R,\bar h_{1}}$, with $\bar h_{1}\leq \bar h$
(recall that by definition $h_{\ell}\leq \bar h$ for all
$\ell\in L(\th)$), and $k_{1}+k_{2}+k_{3}=k-1$. Hence we shall
prove \equ(4.14) for the trees $\th\in\RR^{(k)}_{R,\bar h}$,
for which we can proceed by induction on $k$.

For $(\e,M)\in \widetilde{\DD}(\th,\g)\cap \widetilde{\SS}(\th,\g)$
we have both \equ(4.1) and \equ(4.3) for all $\ell\in L(\th)\setminus
\{\ell_{0},\ell_{e}\}$. Moreover by Lemma \secc(L7)
we have $B k(\th)\geq |n_{\ell} + a n_{e}|$, where $a=0$ if $\ell$
is not on the path $\PPP=\PPP(\ell_{e},\ell_{0})$
and $a\in\{\pm1\}$ otherwise.

For $\ell$ not on the path $\PPP$ one can have $h_{\ell}\ge h$
only if $k(\th)$ is such that $k(\th) > k_{0}=B^{-1} 2^{(h-1)/\t}$
(cf. the proof of Lemma \secc(L10)). If all lines not along
the path $\PPP$ have scales $<h$, consider the line $\ell$
on the path $\PPP$ with scale $h_{\ell} \ge h$ which is
the closest to $\ell_{e}$ (the case in which such a line
does not exist is trivial, because it yields $N_{h}(\th)=0$)).
Then $\ell$ is the exiting line
of a cluster $T$ with $\ell_{e}$ as entering line.
Note that we have both $|n_{\ell}|\ge 2^{(h-1)/\t}$
and $|n_{e}|>2^{(\bar h-2)/\t}$, with $\bar h\ge h$.
As $T$ cannot be a resonance, if $n_{\ell}=n_{e}$ then either $j_{\ell}
\neq j_{e}$, so that
$$ k_{T} > \min\{ B_{2} |n_{e}|^{\b} , B^{-1} C_{2} |p_{j_{e}}|^{\b} \}
> {\rm const.} 2^{(h-1)\b/\t} $$
(cf. Lemma \secc(L12) (i) and the case 2.2.2.1. in the proof of
Lemma \secc(L10)), or $j_{\ell}=j_{e}$, so that
$$ k_{T} > B^{-1} 2^{(\bar h-2)/\t} \ge B^{-1} 2^{(h-2)/\t} $$
(cf. the case 2.2.2.2. in the proof of Lemma \secc(L10)).
If on the contrary $n_{\ell}\neq n_{e}$, by Lemma \secc(L1) one has
$B k(\th) \ge {\rm const.} \min\{|n_{\ell} \pm n_{e}|\}
\ge {\rm const.} 2^{(h-2)s_{2}/\t^{2}}$.
Therefore there exists a constant $\tilde B$ such that
for values $k(\th) \le \tilde k_{0}:= \tilde B^{-1}
2^{(h-1)s_{2}/\t^{2}}$ the bound \equ(4.14) is satisfied.

If $k(\th) > \tilde k_{0}$, we assume that
the bound holds for all trees $\th'$ with $k(\th')<k(\th)$.
Define $E_{h}=c^{-1} 2^{(-2+h)\b/\t}$: we want
to prove that $N_{h}(\th)\le \max\{0,k(\th) E_{h}^{-1}-1\}$.

We proceed exactly as in the proof of Lemma \secc(L10).
The only difference is that, when discussing the case 2.2.1,
one can deduce $|n_{\bar\ell} \pm n_{\ell_{0}}|\ge
{\rm const.} \min\{|n_{\ell_{0}}|,|n_{\bar\ell}|\}^{s_{2} / \t}
\ge {\rm const.} 2^{(h-2)s_{2}/\t^{2}}>E_{h} $ by using that
the quantity $n_{e}$ cancels out as the line $\bar\ell$
is along the path $\PPP$.  \qed

\*

The following result is an immediate consequence of the previous lemma.

\*
\0\lm(L14)
{\it For fixed $k$ the matrices $ L^\ka_{n,j}$ are symmetric;
moreover the following identity holds:
$$ L^\ka_{n,j}= - \bar \chi_{1}(y_{n,j})
\sum_{h=-1}^{\io} C_{h}(x_{n,j}) \sum_{\th\in
\R_{R,h,n,j}^\ka} \Val(\th) ,
\Eq(4.15) $$
where, by definition,
$$C_{h}(x)=\sum_{h_{1}=h+2}^{\io} \chi_{h}(x).
\Eq(4.16) $$  }
\*

\0{\bf Proof.} The previous analysis has shown that the matrices
$L_{n,j}^{(k)}$ are well-defined. Then the matrices are
symmetric by Lemma \secc(L9), where we have established a one to one
correspondence between the trees contributing to $L_{n,j}^{(k)}(a,b)$
and those contributing to $L_{n,j}^{(k)}(b,a)$ such that
the corresponding trees have the same value.
Identity \equ(4.15) follows from the definitions
\equ(3.13) and \equ(3.27). \qed
\*

\0{\bf Remark.} Notice that $C_{h}(x)=1$ when $|x|\le 2^{-h-2}\g$
and $C_{h}(x)=0$ when $|x|\ge 2^{-h-1}\g$.

\*
\0\lm(L15)
{\it Given a tree $\th \in \RR^{(k)}_{R,h,n,j}(a,b)$, for   
$(\e,M)\in \widetilde{\DD}(\th,\g)$ and $\s>0$ one has
$$ \eqalign{
& (i)\quad |\Val(\th)|\le  (Dq^{2})^{k}
2^{-h} \Big( \prod_{h'=-1}^{h} 2^{h' N_{h'}(\th)} \Big)
{\rm e}^{-\s|m_a-m_b|^{\r}}
\prod_{\ell \in L(\th)} p_{j_{\ell}}^{-3s/4}, \cr
&(ii)\quad  | \dpr_\e \Val(\th)|
\le  (Dq^{2})^{k}  2^{-h} |n|
\Big( \prod_{h'=-1}^{h} 2^{2 h' N_{h'}(\th)} \Big)
{\rm e}^{-\s |m_a-m_b|^{\r}}
\prod_{\ell \in L(\th)} p_{j_{\ell}}^{-s_{2}-\a}, \cr
& (iii) \sum_{(n',j')\in \O}\sum_{a',b'=1}^{d_{j'}}
|\dpr_{M_{n',j'}(a',b')}\Val(\th)| \cr
& \qquad \qquad \qquad \qquad
\le (Dq^{2})^{k} 2^{-h}
\Big( \prod_{h'=-1}^{h} 2^{2 h' N_{h'}(\th)} \Big)
{\rm e}^{-\s |m_a-m_b|^{\r}}
\prod_{\ell \in L(\th)} p_{j_{\ell}}^{-s_{2}-\a} . \cr}
\Eq(4.17) $$
for some constant $D$ depending on $\s$ and $\g$.}
\*

\0{\bf Proof.}  The proof follows the same lines as that
of Lemma \secc(L11). We first extract the factor $q^{2k}$
by noticing that a renormalised tree in $\RR^{(k)}_{R}$ has $2k$
end-points. To extract the factor $2^{-h}$ we recall that there is  
at least a line $\ell\neq \ell_{0}$ on scale $h_{\ell}=h$: then
$N_{h}(\th) \geq 1$ and by \equ(4.14) we obtain $k> {\rm const.}
2^{h\b/\t}$, so that $C^{k} 2^{-h} \ge 1$ for a suitable constant $C$.
To extract the factor ${\rm e}^{-\s |m_a-m_b|^{\r}}$ we use
Lemma \secc(L12) (ii) to deduce $\tilde C^{k} {\rm e}^{-\s|m_{a}-
m_{b}|^{\r}} \ge 1$. Hence the bound (i) in \equ(4.17) follows.

When applying the derivative with respect to $\e$ to $\Val(\th)$ we
reason exactly as in Lemma \secc(L11); the only difference is that we
bound $|n_{\ell}|< |n|+ Bk$, which provides in the bound \equ(4.17)
an extra factor $|n|$ with respect to the bound (ii) in \equ(4.9).

The derivative with respect to $M_{n',j'}(a',b')$ gives a sum
of trees with a distinguished line $\ell$ carrying the propagator
$\dpr_{M_{n',j'}(a',b')} g_\ell$ instead of $g_{\ell}$.
As in the case (iii) of Lemma \secc(L11) we have two contributions,
one when the derivative acts on the matrix and the other
(if $i_{\ell}=1$) when the derivative acts on $\chi_{h_{\ell}}$;
by the same arguments as in Lemma \secc(L11) (ii) we obtain a factor of
order $2^{2h}p_{j_\ell}^{-s_2-\a}$.

By Lemma \secc(L7) one has $\min\{|n'-n|,|n'+n|\} < Bk$,
so that the sum over $n'$ is finite and proportional to $k$.
The sum over $j',a',b'$ produces a factor proportional
to $|n'|^{(D-1)+2\a}$ -- reason as in the proof
of \equ(4.9) (iii) in Lemma \secc(L11).
This provides an overall factor of order $|n'|^{D}$.
If $k \ge B_{2} |n|^{\b}$ (with $B_{2}$
defined in Lemma \secc(L12)) this factor can be bounded
by $C^{k}$ for some constant $C$. If $k < B_{2} |n|^{\b}$
then, by Lemma \secc(L12) (i), one must have $i_{\ell}=-1$,
hence $m_{\ell}=m_{\ell}'$, for all lines $\ell$ on the
path $\PPP(\ell_{e},\ell_{0})$: then if $a'\neq b'$ necessarily
the line $\ell$, which the derivative is applied to,
is not on such a path, and the possible values
of $j',a',b'$ are bounded proportionally to $k^{D}$. If $a'=b'$
either $\ell\notin\PPP(\ell_{e},\ell_{0})$ -- and we can reason
as before -- or $\ell\in\PPP(\ell_{e},\ell_{0})$: in the last case
we use the conservation law \equ(3.16) of the momenta $(m_{\ell},
m_{\ell}')$, and we obtain again at most $k^{D}$ terms. \qed

\*

\0{\bf Remark.}
For any fixed $\s>0$ the constant $D$ in \equ(4.17) is proportional
to $\tilde C$, hence grows exponentially in $\s$. As we shall need
for $\tilde C$ to be at worst proportional to $1/\e_{0}$
(in order to have convergence of the series \equ(3.26)),
this means that $\s$ can be taken as large as $O(|\log \e_{0}|)$.

\*

We are now ready to prove the first part of Proposition 1.

\*
\0{\bf Proposition 1 (i)-(ii).}
{\it  There exist constants $c_{0}$, $K_{0}$, $Q_{0}$ and $\s$
such that the following bounds
hold for all $(\e,M)\in \DD(\g)$,  $q< Q_{0}$
and $\h\leq \h_{1}= c_{0} Q_{0}^{-2}$:
$$ \eqalign{
& \left| u_{n,m} \right| < K_{0} |\h| q^{3}
{{\rm e}^{-\s(|n|+|m|)}} ,
\qquad \left| L_{n,j} \right|_\s < K_{0} |\h| q^{2} , \cr
& \left| \dpr_\e L_{n,j} \right|_\s <
K_{0}|n|^{1+s_{2}} |\h| q^{2} ,
\qquad \left| \dpr_\h L_{n,j} \right|_\s <
K_{0} q^{2} , \cr}
\Eq(4.18) $$ 
for all $(n,j)\neq (1,1)$. Moreover the operator norm of the
derivative with respect to $M_{n,j}$ is bounded as
$$ \eqalign{
\Vert \dpr_{M }L\Vert_\s:
& = \sup_{A\in \BBB} {|\dpr_{M }L[A]|_\s\over |A|_\s} \cr
& \le \sup_{n,j\in\O}\sup_{a,b=1,\ldots,d_{j}}
\sum_{n',j'}\sum_{a',b'=1}^{d_{j'}}|
\dpr_{M_{n',j'}(a',b')}L_{n,j}(a,b)|
{\rm e}^{\s(|m_a-m_b|^{\r}-|m_{a'}-m_{b'}|^{\r})} <
K_{0} |\h| q^{2} , \cr}
\Eq(4.19) $$
where the space $\BBB$ is defined in the Remark after
the proof of Lemma \secc(L4).}
\*

\0{\bf Proof.} By definition   $\DD(\g)$ is contained in
all $\DD(\th,\g)$ and in all $\widetilde{\DD}(\th,\g)$, so that we can use
Lemma \secc(L11) and Lemma \secc(L15)  to bound the values of trees.
First we fix an unlabelled tree $\th$ and sum over the values of the
labels: we can modify independently all the end-point labels, the scales,
the type labels and the momenta $m_{\ell}$ if $i_{\ell}\neq -1$
(one has $m_{\ell}=m_{\ell}'$ for $i_{\ell}=-1$). 
Fixed  $(\e,M)$ and $(n_{\ell},j_{\ell})$ there are only
$d_{j_{\ell}}=O(p_{j_{\ell}}^\a)$ possible values for $m_{\ell}$.
This reduces the factors $p_{j_{\ell}}^{-s_{2}-\a}$
to $p_{j_{\ell}}^{-s_{2}}$ in the bounds \equ(4.9) and \equ(4.17).
By summing over the type and scale labels
$\{i_{\ell},h_{\ell}\}_{\ell\in L(\th)}$ (recall that after
fixing the mode labels and $\e$ there are only two possible values
for each $h_{\ell}$ such that $\Val(\th)\neq 0$),
we obtain a factor $4^{k}$, and summing over the possible
end-point labels provides another factor $2^{(D+1)(2k+1)}$.
Finally we bound the number of unlabelled trees of order $k$
by $\bar C^{k}$ for a suitable constant $\bar C$ \cita{HP}.
In \equ(4.9) we can bound
$$ \prod_{h=-1}^{\io} 2^{h N_{h}(\th) } =
\exp \Big( \log 2 \sum_{h=-1}^{\io} h N_{h}(\th) \Big) \le
\exp \Big( {\rm const.} k \sum_{h=-\io}^{\io}
h 2^{-h \b /2\t} \Big) \le C^{k}.
\Eq(4.20) $$
for a suitable constant $C$, and an analogous bound holds
for the products over the scales in \equ(4.17).

Since (see \equ(3.1) and \equ(3.19))
$$ u_{n,m}=\sum_{k=1}^{\io} \h^{k}
\sum_{\th\in \Th_{R,n,m}^{(k)}} \Val(\th),
\Eq(4.21) $$
and, by Lemma \secc(L5), $ \Th_{R,n,m}^{(k)}$ is empty
if $k < B^{-1}|n|$ or $k < B^{-1}|m|$,
we obtain the first bound in \equ(4.18).

Using \equ(4.17) (i), we bound the sum on $\th \in \R_{R,h,n,j}^{(k)}$
exactly in the same way. The main difference is
that $\R_{R,h,n,j}^{(k)}(a,b)$ is empty if $|m_{a}-m_{b}| >
B_{3}^{-1}k^{1/\r}$, by Lemma \secc(L12) (ii). Then
by Lemma \secc(L14), we obtain the second  bound in \equ(4.18).

As for the third bound in \equ(4.18), we have
$$ \eqalign{
& \dpr_\e L_{n,j} = - \bar\chi_{1}(y_{n,j})
\sum_{h=-1}^{\io} C_{h}(x_{n,j}) \sum_{\th\in \R_{R,n,j,h}^{(k)}}
\dpr_\e\Val(\th) \cr & \qquad - \bar\chi_{1}(y_{n,j})\sum_{h=-1}^{\io}
(\dpr_\e C_{h}(x_{n,j})) \sum_{\th\in \R_{R,n,j,h}^{(k)}}\Val(\th) -
(\dpr_\e\bar\chi_{1}(y_{n,j})) \sum_{h=-1}^{\io}
C_{h}(x_{n,j}) \sum_{\th\in \R_{R,n,j,h}^{(k)}}\Val(\th) , \cr}
\Eq(4.22) $$
where the first summand is treated, just like in
the previous cases, by using \equ(4.17) (ii) instead of
\equ(4.17) (i). In the other summands $\Val(\th)$ is bounded
exactly as in the previous cases,
but the derivative with respect to $\e$ gives
in the second summand an extra factor proportional to
$|n| 2^{h} p_{j}^{3\a}$ -- appearing only for those
values of $h$ such that $\chi_{h}(x_{n,j})$ is non-zero
(and for each value of $\e$ there are only two such values so that
the sum over $h$ is finite) -- and in the third summand
a factor proportional to $|n|\,p_{j}^{s_{2}}$. We omit the details,
which can be easily worked out by reasoning as for \equ(4.10)
and \equ(4.12) in the proof of Lemma \secc(L11). Finally we bound
$2^{h}$ by $C^{k}$ as in the proof of Lemma \secc(L15).

The fourth bound in \equ(4.18) follows trivially by noting
that to any order $k$ the derivative with respect to $\h$
of $\h^{k}$ produces $k\h^{k-1}$.

Finally, one can reason in the same way about the derivative
with respect to $M_{n,j}$, by using \equ(4.17) (iii),
so that \equ(4.19) follows. \qed

\*\*
\section(5,Whitney extension and implicit function theorems)

\0{\bf 5.1. Extension of $\UUUU$ and $\LLLL$}\*

\0In this section we extend the function $L_{n,j}$,
defined in  $\DD(\g)$, to the larger set $\DD_{0}$.

\*
\0\lm(L16)
{\it The following statements hold true.\\
(i) Given $\th\in \R^\ka_{R,h,n,j}$, we can extend
$\Val(\th)$ to a function, called $\Val^{E}(\th)$,
defined and $C^{1}$ in $\DD_{0}$, and $L_{n,j}(\h,\e,M;q)$
to a function $L_{n,j}^{E} \= L_{n,j}^{E}(\h,\e,M;q)$ such that
$$ L_{n,j}^{E} = - \bar\chi_{1}(y_{n,j})
\sum_{h=-1}^{\io} C_{h}(x_{n,j}) \sum_{k=1}^{\io} \h^{k}
\sum_{\th \in \RR_{R,n,j,h}^{(k)}} \Val^{E}(\th) ,
\Eq(5.1) $$
satisfies for any $(\e,M)\in \DD_{0}$ the same bounds
in \equ(4.18) and \equ(4.19) as $L_{n,j}^{E}(\h,\e,M;q)$ in $\DD(\g)$.
Furthermore $\Val(\th)=\Val^{E}(\th)$
for any $(\e,M)\in \widetilde{\DD}(\th,2\g) \subset \widetilde{\DD}(\th,\g)$
and $\Val^{E}(\th)=0 $ for 
$(\e,M)\in  \DD_{0}\setminus  \widetilde{\DD}(\th,\g)$.\\
(ii) In the same way, given $\th\in \Th^\ka_{R,n,m}$, 
we can extend $\Val(\th)$ to a function $\Val^{E}(\th)$
defined and $C^{1}$ in $\DD_{0}$, and $U_{n,j}(\h,\e,M;q)$
to a function $U_{n,j}^{E}(\h,\e,M;q)$ such that
$u_{n,m}^{E} \= u_{n,m}^{E}(\h,\e,M;q)$, given by
$$ u^{E}_{n,m}=\sum_{k=1}^{\io} \h^{k}
\sum_{\th\in \Th_{R,n,m}^{(k)}} \Val^{E}(\th) ,
\Eq(5.2) $$
satisfies for any $(\e,M)\in \DD_{0}$ the same bounds
in \equ(4.18) as $u_{n,m}$ in $\DD(\g)$.\\
Furthermore $\Val(\th)=\Val^{E}(\th)$
for any $(\e,M)\in \DD(\th,2\g) \subset \DD(\th,\g)$
and $\Val^{E}(\th)=0 $ for 
$(\e,M)\in  \DD_{0}\setminus  \DD(\th,\g)$.}
\*

\0{\bf Proof.} We prove first the statement for the
case $\th\in \RR^{(k)}_{R,h,n,j}$.
We use the $C^{\io}$ compact support function $\chi_{-1}(t):
\RRR\to \RRR^+$, introduced in Definition \secc(D4).
Recall that $\chi_{-1}(t)$ equals $0$ if $|t|<\g$ and $1$
if $|t|\ge 2\g$, and $|\dpr_t\chi_{-1}(t)|\le C \g^{-1}$,
for some constant $C$.

Given a tree $\th\in\R_{R,h,n,j}^{(k)}$, we define
$$ \Val^{E}(\th) = \!\!\!\!\! \!\!\!\!\!
\prod_{\ell\in L(\th)\setminus
\{\ell_{e},\ell_{0}\}\atop i_{\ell}=1} \!\!\!\!\!
\chi_{-1}(|x_{n_{\ell},j_{\ell}}| |n_{\ell}|^\t) \!\!\!\!\!
{\mathop{\prod}_{\ell_{1},\ell_{2}\in L(\th)}}^{\hskip-0.4truecm**}
\chi_{-1}(|\d_{n_{\ell_{1}},j_{\ell_{1}}}-
\d_{n_{\ell_{2},j_{\ell_{2}}}}| 
|n_{\ell_{1}}- n_{\ell_{2}}|^{\t_{1}})) \, \Val(\th) ,
\Eq(5.3) $$
where $\prod^{**}_{\ell_{1},\ell_{2}\in L(\th)}$ is the product 
on the pairs $\ell_{1}\prec \ell_{2}\in L(\th)$ such that
$\prod_{\ell\in \PPP(\ell_{1},\ell_{2})}\s(\ell)\s(\ell_{1})=1$,
$i_{\ell_{j}}=1,0$, $n_{\ell_{1}}\neq n_{\ell_{2}}$, and either both
$\ell_{1},\ell_{2}$ are on the path connecting $e$ to $v_{0}$ or both of
them are not on such a path. The sign $\prod_{\ell \in
\PPP(\ell_{1},\ell_{2})} \s(\ell)\s(\ell_{1})$ is such
that $ |n_{\ell_{1}}- n_{\ell_{2}}|\le n$.

By definition $\Val^{E}(\th)=\Val(\th)$ for $(\e,M)\in
\widetilde{\DD}(\th,2\g)$ as in this set the scale functions $\chi_{-1}$
in the above formula are identically equal to $1$.

By definition ${\rm supp}(\Val^{E}(\th))\subset \widetilde{\DD}(\th,\g)$,
as the scale functions $\chi_{-1}$ in the above formula are
identically equal to $0$ in the complement of $\widetilde{\DD}(\th,\g)$
with respect to $\DD_{0}$.

To bound the derivatives the only fact that prevents us from
simply applying \equ(4.17) (ii-iii) is the presence of the extra
terms due to the derivatives of the $\chi_{-1}$ functions.
Each factor of the first product in \equ(5.3), when derived,   
produces an extra factor proportional to $2^{h_{\ell}}
p_{j_{\ell}}^{3\a} |n_{\ell}|^{\t+1}$. Note that a summand of
this kind appears only if $i_{\ell}=1$ and $(\e,M)$ is such that
$$ 2^{-h_{\ell}-1} \g \le x_{n_{\ell},j_{\ell}}
\le {2\g \over |n_{\ell}|^{\t}} .
\Eq(5.4) $$
This implies $|n_{\ell}| < 2^{(h_{\ell}+1)/\t}$, so that the
presence of the extra factor simply produces, in \equ(4.17) (ii),
a larger constant $D$ and a larger exponent -- say $4$ -- instead
of $2$ in the factor $2^{2 h' N_{h'}(\th)}$.
Each factor of the second product produces an extra factor
$|n_{\ell_{1}} - n_{\ell_{2}}|^{\t_{1}+1}$,
which can be bounded by $C^{k}$.

Therefore the derivatives of
$L_{n,j}^{E}$ respect the same bounds \equ(4.18) as $L_{n,j}$
modulo a redefinition of the constants $c_{0}$, $K_{0}$.
As these bounds are uniform (independent of $(n,j)$), then
$L^{E}_{n,j}$ is a $C^1$ function of $(\e,M)$. 

We proceed in the same way for $\th\in \Th_{R,n,m}$:
$$ \Val^{E}(\th) = \!\!\!\!\! \! \prod_{\ell\in L(\th):  i_{\ell}=1}
\!\!\!\!\chi_{-1}(|x_{n_{\ell},m_{\ell}}| |n_{\ell}|^\t)
{\mathop{\prod}_{\ell_{1},\ell_{2}\in L(\th)}}^{\hskip-0.3truecm***}
\chi_{-1}(|\d_{n_{\ell_{1}},j_{\ell_{1}}}
\d_{n_{\ell_{2},j_{\ell_{2}}}}| |n_{\ell_{1}} -
n_{\ell_{2}}|^{\t_{1}}))\Val(\th) ,
\Eq(5.5) $$
where now the product $\P^{***}$ runs on the pairs of lines
$\ell_{1}\prec\ell_{2}$ such that
$\prod_{\ell\in \PPP(\ell_{1},\ell_{2})}\s(\ell)\s(\ell_{1})=1$,
$i_{\ell_{j}}=1,0$, and $n_{\ell_{1}}\neq n_{\ell_{2}}$. \qed

\*
\0{\bf Proposition 1 (iii).} {\it$L^{E}$ is differentiable in
$(\e,M)\in \DD_{0}$ and satisfy the bounds
$$ \eqalign{
& \left| \dpr_\e L^{E}_{n,j}(a,b) \right| <
C_{1} |n|^{1+s_{2}} {\rm e}^{-\s|m_a-m_b|^{\r} } |\h| q^{2} , \cr
\sum_{(n',j')\in \O}\sum_{a',b'=1}^{d_{j'}}
& \left| \dpr_{M_{n',j'}(a',b')}L^{E}_{n,j}(a,b) \right|
{\rm e}^{|m_{a}-m_{b}|^{\r}} < C_{1} |\h| , \cr}
\Eq(5.6) $$
where $C_{1}$ is a suitable constant.}
\*

\0{\bf Proof.} Simply combine the proof of Lemma \secc(L16)
with that of Proposition 1 (ii). \qed

\*\*
\0{\bf 5.2. The extended $\QQQQ$ equation}\*

\0Going back to \equ(2.12), we can extend it to all $\DD_{0}$ by using
$U^{E}_{n,j}$ instead of $U_{n,j}$ for all $(n,j)\neq (1,1)$;
we obtain the equation
$$ D^{s} q= f_{1,V}(u^{E})=\sum_{n_{1}+n_{2}-n_3=1
\atop m_{1}+m_{2}-m_3=V}
u^{E}_{n_{1},m_{1}}u^{E}_{n_{2},m_{2}}u^{E}_{n_3,m_3} .
\Eq(5.7) $$
The leading order is obtained for $n_{i}=1$ and $m_{i}\in \L_{1}$
for all $i=1,2,3$, namely at $\h=0$ we have a nonlinear
algebraic equation for $q$,
$$ D^{s} q = 3^{D} q^3\, ,
\Eq(5.8) $$ 
with solution  $q_{0}=\sqrt{D^{s} 3^{-D}}$.
We can now prove the following result.

\*
\0{\bf Proposition 1 (iv).}
{\it There exists $\h_{0}$ such that for all $|\h|\le \h_{0}$ and
$(\e,M)\in \DD_{0}$, equation \equ(5.7) has a solution $q^{E}(\e,M;\h)$,
which is analytic in $\h$ and $C^1$ in $(\e,M)$; moreover
$$|\dpr_\e q^{E}(\e,M;\h)|\le K |\h| , \qquad
\sum_{(n,j)\in\O}  \sum_{a,b=1}^{d_{j}}
\left| \dpr_{M_{n,j}(a,b)} q^{E}(\e,M;\h) \right|_{\s} \le K |\h| ,
\Eq(5.9) $$
for a suitable constant $K$.}
\*

\0{\bf Proof.}
Set $Q_{0}:=2q_{0}$. Then there exists $\h_{1}$ such that $u^{E}$ is
analytic in $\h,q$ for $|\h|\leq \h_{1}$ and $|q|\leq Q_{0}$
and $C^1$ in $(\e,M)$. By the implicit function theorem,
there exists $\h_{0}\le \h_{1}$ such that for all $|\h|\le \h_{0}$
there is a solution $q^{E}\=q^{E}(\h,\e,M)$ of the $Q$ equation
\equ(5.7) such that $|q^{E}| < 3q_{0}/2 $. By definition of the
extension $u^{E}$, the  equation \equ(5.7) coincides \equ(2.12)
on $\DD(2\g)$. The bounds on the derivatives
follow from Lemma \secc(L11) and Lemma \secc(L16). \qed

\*

We now define
$$ U^{E}_{n,j}(\h,\e,M)= U^{E}_{n,j}(\h,\e,M;
q^{E}(\h,\e,M))\,,\qquad  L^{E}_{n,j}(\h,\e,M)= L^{E}_{n,j}(\h,\e,M;
q^{E}(\h,\e,M)) .
\Eq(5.10) $$ 

\*
\0{\bf Proposition 1 (v).}
{\it There exists a positive constant $K_{1}$ such that
the matrices $L_{n,j}^{E}(\h,\e,M) $ satisfy the bounds
$$ \eqalign{
& |L^{E}(\h,\e,M)|_\s \leq |\h| K_{1} \,,
\quad |\dpr_\e L^{E}_{n,j}(\h,\e,M)|_{\s}
\leq |\h| K_{1} |n|^{1+s_{2}} , \cr
& \sum_{(n,j)\in \O}\sum_{a,b=1}^{d_{j}}
\left| \dpr_{M_{n,j}(a,b)} L^{E}(\h,\e,M) \right|_\s
{\rm e}^{-\s|m_{a}-m_{b}|^{\r}} \le |\h| K_{1} , \cr}
\Eq(5.11) $$
and the coefficients $U^{E}_{n,j}(\h,\e,M)$ satisfy the bounds
$$ \left| U^{E}_{n,j}(\h,\e,M) \right|
\le |\h|  K_{1} {\rm e}^{-\s (|n|+|p_{j}|^{1/2})} ,
\Eq(5.12) $$
uniformly for $(\e,M)\in \DD_{0}$.}
\*

\0{\bf Proof.} It follows trivially from the bounds \equ(5.9)
and from the bounds of Lemma \secc(L16). \qed

\*\*
\section(6,Proof of Proposition 2)

\0{\bf 6.1. Proof of Proposition 2 (i)}\*

\0Let us consider the compatibility equation \equ(2.11) where
$L_{n,j}=L^{E}_{n,j}(\h,\e,M)$. One can rewrite \equ(2.11) as
$$ \bar\chi_{1} (y_{n,j})M_{n,j}=  L^{E}_{n,j}(\h,\e,M)) \=
\h\bar\chi_{1} (y_{n,j}) \tilde L^{E}(\h,\e,M) ,
\Eq(6.1) $$
with $\tilde L^{E}(\h,\e,M)=O(1)$, so that \equ(6.1)
has for $\h=0$ the trivial solution $M_{n,j}=0$.

The bounds of Proposition 1 (v) imply that the Jacobian of the
application $\tilde L^{E}(\h,\e,M): \BBB\to \BBB$ is bounded in the
operator norm ($\BBB$ is defined in the Remark before Definition
\secc(D3)). Thus there exists $\h_{0},K_{2}$ such that, for
$|\h|\le \h_{0}$ and for all $(\e,M)\in \DD(2\g)$, we can apply
the implicit function theorem to \equ(6.1) and obtain a solution
$M_{n,j}(\h,\e)$, which satisfies the bounds
$$ |M_{n,j}|_\s \le K_{2} |\h| , \qquad |\dpr_{\e} M_{n,j}(\h,\e)|_{\s}
\le K_{2} |n|^{1+s_{2}} |\h| ,\qquad
|\dpr_{\h} M_{n,j}(\h,\e)|_{\s} \le K_{2} ,
\Eq(6.2) $$
for a suitable constant $K_{2}$.

Finally we fix $\e_{0}\le \h_{0}$,  $\h=\e$ and set
(with an abuse of notation) $M_{n,j}(\e)=M_{n,j}(\h=\e,\e)$,
so that, by noting that
$$ {\der \over \der \e} M_{n,j} (\e) = \dpr_{\h} M_{n,j} (\h,\e) +
\dpr_{\e} M_{n,j} (\h, \e) ,
\Eq(6.3) $$
we deduce from \equ(6.2) the bound \equ(2.26).

\*\*
\0{\bf 6.2. Proof of Proposition 2 (ii) -- measure estimates}\*

\0We now study the measure of the set \equ(2.27).
By definition this is given by the set of $\e\in \EE_{0}(\g)$
such that the further Diophantine conditions
$$ x_{n,j}(\e) :=\left\Vert (\d_{n,j}I+p_{j}^{-s}
\bar\chi_{1}(y_{n,j})\,M_{n,j}(\e))^{-1} \right\Vert^{-1}
\ge {2\g\over |n|^\t}
\Eq(6.4) $$
are satisfied for all $(n,j)\in\O$ such that $(n,j)\neq (1,1)$.
Recall that $(n,j)\in\O$ implies $n>0$.
By Lemma \secc(L3) (iii) one has
$$ x_{n,j}(\e) \ge \min_{i} |\l^{(i)}
(\d_{n,j}I+p_{j}^{-s} \bar\chi_{1}(y_{n,j})\,M_{n,j})|
= \min_{i}|\d_{n,j} + p_{j}^{-s+\a} \n_{n,j}^{(i)}(\e))| ,
\Eq(6.5) $$
since the matrices are symmetric. Recall that $p_{j}^{\a}
\n_{n,j}^{(i)}(\e)$ are the eigenvalues of $\bar\chi_{1}(y_{n,j})\,
M_{n,j}^{(i)}(\e)$ and that $d_{j} \le C_{1} p_{j}^{\a}$
(cf. Definition \secc(D2)).

Then we impose the conditions
$$ \left| \d_{n,j} + p_{j}^{-s+\a} \n_{n,j}^{(i)}(\e) \right|
\ge  {2\g\over n^\t} \qquad \forall (n,j) \in \O \setminus \{(1,1)\} ,
\quad i=1,\dots,d_{j} ,
\Eq(6.6) $$
and recall that $M_{n,j}= 0$ (i.e. $\n_{n,j}^{(i)}=0$) if $(n,j)
\notin \O$, so that for $(n,j) \notin \O$ the Diophantine
conditions \equ(6.6) are surely verified, by \equ(2.1).

Call $\CCCCCC$ the set of values of $\e\in\EE_{0}(\g)$ which
verify \equ(6.6). We estimate the measure of the subset of
$\EE_{0}(\g)$ complementary to $\CCCCCC$,
i.e. the set defined as union of the sets
$$ \III_{n,j,i}:= \left \{\e\in \EE_{0}(\g) : |\d_{n,j}+p_{j}^{-s+\a}
\n_{n,j}^{(i)}(\e) |\le {2\g\over n^\t} \right\}
\Eq(6.7) $$
for $(n,j)\in\O$ and $i=1,\ldots,d_{j}$.

Given $n$, the condition $(n,j)\in \O$ implies that $p_{j}$ can assume
at most $\e_{0} n+1$ different values -- cf. \equ(2.16).
On a $(D-1)$-dimensional sphere of radius $R$
there are at most $O(R^{D-1})$ integer points, hence
the number of values which $j$ can assume is bounded
proportionally to $n^{D-1} (1+\e_{0} n)$. Finally $i$ assumes
$d_{j}\le C_{1} p_{j}^{\a}$ values.

Since $\m\in \MM$ we have,
for $n \le (\g_{0}/4\e_{0})^{1/(\t_{0}+1)}$,
$$ \left| \d_{n,j} + p_{j}^{-s+\a} \n_{n,j}^{(i)}(\e) \right|
\ge |(D+\m) n- p_{j}-\m|- 2\e_{0} n
\ge {\g_{0}\over 2 n^{\t_{0}}} ,
\Eq(6.8) $$ 
so that we have to discard the sets $\III_{n,j,i}$ only for $n \ge
(\g_{0}/4\e_{0})^{1/(\t_{0}+1)}$.

Let us now recall that for a symmetric matrix $M(\e)$
depending smoothly on a parameter $\e$,
the eigenvalues are $C^1$ in $\e$ \cita{Ka}.
Then the measure of each $\III_{n,j,i}$ can be bounded from above by
$$ {4\g\over n^\t } \sup_{\e\in \EEsub_{0}(\g)}
\left| \left( {\der \over \der\e}
\left( \d_{n,j} + p_{j}^{-s+\a} \n_{n,j}^{(i)}(\e) \right)
\right)^{-1} \right| .
\Eq(6.9) $$
where one has
$$ \left| {\der \over \der\e} \left(
\d_{n,j}+p_{j}^{-s+\a} \n_{n,j}^{(i)}(\e) \right) \right| \ge {n \over 2} .
\Eq(6.10) $$
This can be obtained as follows. Proving \equ(6.10)
requires to find lower bounds for 
$$ \left| {\der \over \der \e} 
\left( - \o n + p_{j} + \m + p_{j}^{-s}
\bar\chi(y_{n,j})\l_{n,j}^{(i)}(\e) \right) \right| , $$
where $\l_{n,j}^{(i)}(\e)$ are the eigenvalues of $M_{n,j}(\e)$
(i.e. $\bar\chi_{1}(y_{n,j})\l_{n,j}^{(i)}(\e)=
p_{j}^{\a} \n_{n,j}^{(i)}(\e)$).
The eigenvalues $\l_{n,j}(\e)^{(i)}$ are $C^1$ in $\e$,
so that, by Lidskii's lemma \cita{Ka}, one has
$$ \left| {\der \over \der \e} \l_{n,j}^{(i)}(\e) \right|
\le d_{j} \left| {\der \over \der \e} M_{n,j} \right|_\io
\le C_{1} K_{2} \left( 1 + \e_{0} n^{1+s_{2}} \right) n^{\a} ,
\Eq(6.11) $$
where we have used \equ(6.3) and \equ(6.2). Since
$s_{2}+\a \le s-\a$, we obtain
$$ \left| {\der \over \der \e}
\left( - \o n + p_{j} + \m + p_{j}^{-s}
\bar\chi(y_{n,j}) \l_{n,j}^{(i)}(\e) \right) \right| \ge n/2  , $$
which implies \equ(6.10). 

Recall that $p_{j}$ is bounded proportionally to $n$.
Then for fixed $n$ we have to sum over ${\rm const.}
(1+\e_{0}n)n^{D-1}$ values of $j$ and over $d_{j} \le C_{1} p_{j}^{\a}
\le {\rm const.} n^{\a} \le {\rm const.} n$.

Therefore we have
$$ \eqalign{
\sum_{(n,j)\in\O} \sum_{i=1}^{d_{j}}
\hbox{meas} \left( \III_{n,j,i} \right) & \le {\rm const.}
\sum_{n \ge (\g/4\e_{0})^{-1/(\t_{0}+1)}} \g |n|^{D}
\left( {1 \over |n|^{\t+1} } + { \e_{0} \over |n|^{\t}}
\right) \cr & \le {\rm const.} \left( \e_{0}^{(\t-D)/(\t_{0}+1)} +
\e_{0}^{ 1+ (\t-D-1)/(\t_{0}+1)} \right) , \cr}
\Eq(6.12) $$
provided $\t>D+1$. Therefore the measure
is small compared to that of $\EE_{0}(\g)$ -- which is of
order $\e_{0}$-- if $\t>\max\{\t_{0}+D+1,D+1\}=\t_{0}+D+1$.

\*\*
\section(7,Generalisations and proof of Theorem 1)

\0{\bf 7.1. Equation \equ(1.4): proof of Theorem 1
in $\DDDDD$ $\maggiore$ 2}\*

\0In order to consider  equation \equ(1.4) we only need to make
few generalisations. By our assumptions
$$ f(x,u,\bar u)= g(x,\bar u)+\dpr_{\bar u} H(x,u,\bar u) ,$$
with $H$ real valued. For simplicity we discuss explicitly
only the case with odd $p$ in \equ(1.3) and $g$ odd in $u$.
Considering also even $p$ should require considering an expansion
in $\sqrt{\e}$: this would not introduce any technical difficulties,
but on the other hand would require a deeper change in notations.

We modify the tree expansion, analogously to what done in \cita{GP}.
The change of variables \equ(1.8) transforms each monomial
in \equ(1.3) into a monomial $\e^{(p_{1}+p_{2}-1)/2}
a_{p_{1},p_{2}}(x) \, u^{p_{1}}\bar u^{p_{2}}$; we can take into
account the contributions arising from $g(x,\bar u)$, by considering
the corresponding Taylor expansion and putting $p_{1}=0$ and $p_{2}
\ge 3$. All the other contributions are such $p_{1} a_{p_{1},p_{2}}=
(p_{2}+1) a_{p_{2}+1,p_{1}-1}$ (by the reality of $H$ and of
the functions $a_{p_{1},p_{2}}$).

Each new monomial produces internal nodes of order
$k_v = (p_{v,1}+p_{v,2}-1)/2 \in \NNN$, such that $k_{v} \ge 2$,
with $p_{v,1}+p_{v,2}$ entering lines among which $p_{v,1}$ have
sign $\s=1$ and $p_{v,2}$ have sign $\s=-1$; note that the case
previously discussed corresponds to $(p_{v,1},p_{v,2})=(2,1)$.
Hence, with the notations of Section 3.3, we can write 
$s_{v}=p_{v,1}+p_{v,2}$, with $s_{v}$ odd.

Each internal node $v$ has labels $k_{v},p_{v,1},p_{v,2},m_{v}$,
with the mode label $m_{v}\in \ZZZ^{D}$. The node factor associated
with $v$ is $a_{p_{v,1},p_{v,2},m_{v}}$, namely 
the Fourier coefficient with index $m_{v}$ in the Fourier expansion
of the function $a_{p_{v,1},p_{v,2}}$; by the analyticity assumption on
the non-linearity the Fourier coefficients decay exponentially
in $m$, that is
$$ \left| a_{p_{v,1},p_{v,2},m_{v}} \right| \le
A_{1} {\rm e}^{-A_{2}|m|} ,
\Eq(7.1) $$
for suitable constants $A_{1}$ and $A_{2}$.

The conservation laws \equ(3.15) and \equ(3.16) have to be
suitably changed. We can still write that $n_{\ell}$
is given by the right hand side, the only
difference being that $L(v)$ contain $s_{v}$ lines
(and each line $\ell\in L(v)$ has its own sign $\s(\ell)$).
On the contrary \equ(3.16) for $m_{\ell}'$ has to be changed
in a more relevant way: indeed one has
$$ m_{\ell}'= m_{v} + \sum_{\ell' \in L(v)} \s(\ell') m_{\ell'} ,
\Eq(7.2) $$
with $L(v)$ defined as before.

The order of any tree $\th$ is still defined as in \equ(3.18),
and, more generally, all the other labels are defined exactly
as in Section 3.3.

The first differences appear when one tries to bound the momenta
of the lines in terms of the order of the tree. In fact one has
$$ |E(\th)| \le 1 + \sum_{v\in V_{0}(\th)} \left( s_{v} - 1 \right) ,
\Eq(7.3) $$
which reduces to the formula given in the proof of Lemma \secc(L5)
only for $s_{v}\le 3$. One has $s_{v}=2k_{v}+1$, so that
$$ \sum_{v\in V_{0}(\th)} \left( s_{v} - 1 \right) =
2 \sum_{v\in V_{0}(\th)} k_{v} = 2 k ,
\Eq(7.4) $$
and one can still bound $|n_{\ell}| \le B k$ for
any tree $\th\in\Th^{(k)}$ and any line $\ell\in L(\th)$.

The conservation law \equ(7.1) gives, for any line $\ell\in L(\th)$,
$$ \max\{|m_{\ell}|,|m_{\ell}'|\} \le B k +
\sum_{v \in V_{0}(\th)} |m_{v}| ,
\Eq(7.5) $$
for some constant $B$. The bound in \equ(7.5) is obtained
by reasoning as in  the proof of Lemma \secc(L5);
the last sum is due to the mode labels of the internal nodes.
Thus the bound on $n_{\ell}$ in Lemma \secc(L5) still holds, while the
bounds on $m_{\ell},m_{\ell}'$ have to be replaced with \equ(7.5).
The same observation applies to Lemma \secc(L7).

Also Lemma \secc(L9) still holds. The proof proceeds as follows.
The tree $\th_{1}$ which one associates with
each $\th\in\RR_{R,n,j,h}(a,b)$ is the
tree in $\RR_{R,n,j,h}(b,a)$ defined as follows.
\vskip.2truecm
\01. As in the proof of Lemma \secc(L9).
\vskip.1truecm
\02. As in the proof of Lemma \secc(L9).
\vskip.1truecm
\03. Let $\bar v$ be a node along the path $\PPP=\PPP(\ell_{e},
\ell_{0})$ and let $\ell_{1},\dots,\ell_{s}$, with $s=s_{\bar v}$ be
the lines entering $\bar v$; suppose that $\ell_{1}$
is the line belonging to the path $\PPP \cup \{\ell_{e}\}$.
If $\s(\ell_{1})=1$ we change all the signs of the other lines,
i.e. $\s(\ell_{i})\to - \s(\ell_{i})$ for $i=2,\ldots,s$,
whereas if $\s(\ell)=-1$ we do not change the signs.
\vskip.1truecm
\024. As in the proof of Lemma \secc(L9).
\vskip.2truecm
Then one can easily check that the reality of $H$ implies that the
tree $\th_{1}$ is well defined (as an element of $\RR_{R,n,j,h}(b,a)$)
and has the same value as $\th$.

\*

\0{\bf Remarks.} (1) Note that item 3. above reduces to item 3.
of Lemma \secc(L9) if $s_{v}=3$ for each internal node $v$.\\
(2) If the node $\bar v$ has $p_{\bar v,1}=0$
(i.e. the monomial associated to it arises from the function
$g(x,\bar u)$) then the operation in item 3. is empty.

\*

Therefore we can conclude that the counterterms are still symmetric.

The analysis of Section \secc(4) can be performed almost
unchanged. Here we confine ourselves to show the few changes
that one has to take into account.

The first relevant difference appears in Lemma \secc(L10).
Because of the presence of the
mode labels of the internal nodes the bound \equ(4.5) on $N_{h}(\th)$
does not hold anymore, and it has to be replaced with
$$ N_{h}(\th) \le  \max \Big\{ 0 , c \Big(
k(\th) + \sum_{v\in V_{0}(\th)} |m_{v}| \Big)
2^{(2-h)\b/\t} - 1 \Big\} , 
\Eq(7.6) $$
for a suitable constant $c$. The proof of \equ(7.6)
proceeds as the proof of Lemma \secc(L10) in Section 4.
We use that in \equ(4.7) for $m=1$ one has
$$ k(\th) - k(\th_{1}) + \sum_{v\in V_{0}(\th)} |m_{v}| -
\sum_{v\in V_{0}(\th_{1})} |m_{v}| = k_{T} +
\sum_{v\in V_{0}(T)} |m_{v}| ,
\Eq(7.7) $$
and, except for item 2.2.2.1., we simply bound
the right hand side of \equ(7.7) with $k_{T}$.
The only item which requires a different argument is item 2.2.2.1.,
where instead of the bound $|m_{\bar\ell}-m_{\ell_{0}}| \le
\sum_{i}|m_{i}|$ we have, by \equ(7.5),
$$ \left| m_{\bar\ell} - m_{\ell_{0}} \right| \le 
\sum_{v\in \PPP(\bar\ell,\ell_{0})} |m_{v}| + \sum_{i} |m_{i}|
\le B k_{T} + \sum_{v\in V_{0}(T)} |m_{v}| \le \max\{B,1\}
\Big( k_{T} + \sum_{v\in V_{0}(T)} |m_{v}| \Big) , $$
where $v\in\PPP(\bar\ell,\ell_{0})$ means that the node $v$ is
along the path $\PPP(\bar\ell,\ell_{0})$ (i.e. $\ell_{v}\in
\PPP(\bar\ell,\ell_{0})\cup \ell_{0})$) and the sum over $i$
is over all sub-trees which have the root lines entering one
of such nodes.

\*

\0{\bf Remark.} Note that if the coefficients $a_{p_{1},p_{2}}(x)$
in \equ(1.3) are just constants (i.e. do not depend on $x$),
then $m_{v}\=0$ and \equ(7.6) reduces to \equ(4.5).

\*

Moreover in \equ(4.9) we have a further product
$$ \prod_{v\in V_{0}(\th)} A_{1} {\rm e}^{-A_{2}|m_{v}|} ,
\Eq(7.8) $$
while the product of the factors $2^{hN_{h}(\th)}$ can be written as
$$ \Big( \prod_{h=-1}^{h_{0}} 2^{h N_{h}(\th)} \Big)
\Big( \prod_{h=h_{0}+1}^{\io} 2^{h N_{h}(\th)} \Big) \le
2^{h_{0}k} \prod_{h=h_{0}+1}^{\io} 2^{h N_{h}(\th)} \Big) ,
\Eq(7.9) $$
with $h_{0}$ to be fixed, where the last product,
besides a contribution which can be bounded as in \equ(4.20),
gives a further contribution
$$ \prod_{h=h_{0}}^{\io} \prod_{v\in V_{0}(\th)}
2^{ c h |m_{v}| 2^{(2-h)\b/\t} } \le \prod_{v\in V_{0}(\th)}
\exp \Big( {\rm const.} |m_{v}| \sum_{h=h_{0}}^{\io}
h 2^{-h\b/\t} \Big) ,
\Eq(7.10) $$
so that we can use part of the exponential factors in \equ(7.8)
to compensate the exponential factors in \equ(7.10),
provided $h_{0}$ is large enough (depending on $\t$).

Another consequence of \equ(7.2) is in Lemma \secc(L12):
item (ii) has to be replaced with
$$ |m_{a}-m_{b}| \le {\rm const.} k^{1/\r} +
\sum_{v\in V_{0}(\th)} |m_{v}| ,
\Eq(7.11) $$
because each internal node $v$ contributes a momentum $m_{v}$
to the momenta of the lines following $v$. Up to this observation,
the proof of \equ(7.11) proceeds as in the proof of Lemma \secc(L12).

Therefore also the bound \equ(4.14) of Lemma \secc(L13) has to changed
into \equ(7.6), for all $h \le \bar h$. The proof proceeds
as that of Lemma \secc(L13) in Section 4, with the changes outlined
above when dealing with the case 2.2.2.1.

The property \equ(7.11) reveals itself also in the proof of Lemma
\secc(L15). More precisely, in order to extract a factor
${\rm e}^{-\s|m_{a}-m_{b}|^{\r}}$, we use that \equ(7.11) implies
(recall that $\r<1$)
$$ \left| m_{a} - m_{b} \right|^{\r} \le
C \Big( k + \sum_{v\in V_{0}(\th)} |m_{v}| \Big),
\Eq(7.12) $$
for a suitable $\r$-dependent constant $C$, so that we can write,
for some other constant $\tilde C$,
$$ {\rm e}^{\s |m_{a}-m_{b}|^{\r}} \le 
\tilde C^{k} \prod_{v\in V_{0}(\th)} {\rm e}^{\s |m_{v}|} ,
\Eq(7.13) $$
where $\s$ has to be chosen so small (e.g. $|\s|<A_{2}/4$, with
$A_{2}$ given in \equ(7.1)) that the last product in \equ(7.13) can
be controlled by part of the exponentially decaying factors
${\rm e}^{- A_{2}|m_{v}|}$ associated with the internal nodes. This
means, in particular, that $\s$ cannot be arbitrarily large when
$\e_{0}$ becomes small (cf. the Remark after the proof of Lemma \secc(L15)).

As in \equ(4.9) also in \equ(4.17) there are the further factors
\equ(7.8), which can be dealt with exactly as in the previous case.

Besides the issues discussed above, there is no other substantial
change with respect to the analysis of Sections \secc(4) to \secc(6).

\*\*
\0{\bf 7.2. Equation \equ(1.1) in dimension 2:
proof of Theorem 1 in $\DDDDD$ = 2}\*

\0We can consider more general nonlinearities in the case $D=2$,
that is of the form \equ(1.3) without the simplifying assumption \equ(1.4).
Indeed in such case the counterterms are $2\times 2$ matrices
(cf. Lemma \secc(L2)), so that we can bound $x_{n,j}$ by the absolute
value of the determinant of $\d_{n,j}I + \bar \chi_{1}(y_{n,j})
M_{n,j}p_{j}^{-s}$, which is a $C^{1}$ function of $\e$
(we have proved only $C^1$ but it should be obvious that we can
bound as many derivatives of $L^{E}_{n,j}$ as we need to, possibly
by decreasing the domain of convergence of the functions involved).

Set for notational simplicity $\overline M_{n,j}=\bar\chi_{1}(y_{n,j})
M_{n,j}$. Let us evaluate the measure of the Cantor set
$$ \EE_{1} = \left\{ \e \in \EE_{0}(\g) :
|\d_{n,j}^{2} +p_{j}^{-s} {\rm tr} \,
\overline M_{n,j} \d_{n,j} +p_{j}^{-2s} \det \overline M_{n,j}| \ge
2\g|n|^{-\t} \right\} ,
\Eq(7.14) $$
following the scheme of Section \secc(6). Here we are using
explicitly that for $D=2$ one has
$$ \det \left( \d_{n,j}I + p_{j}^{-s} \overline M_{n,j}\right) =
\d_{n,j}^{2} +p_{j}^{-s} {\rm tr} \,
\overline M_{n,j} \d_{n,j} + p_{j}^{-2s} \det \overline M_{n,j} ,
\Eq(7.15) $$
because $M_{n,j}$ is a $2 \times 2$ matrix.

We estimate the measure of the complement of $\EE_{1}$
with respect to $\EE_{0}(\g)$, which is the union of the sets
$$ \III_{n,j}:= \left\{ \e\in \EE_{0}(\g) :\; |\d_{n,j}^{2} +
p_{j}^{-s}a_{n,j} \d +p_{j}^{-2s}b_{n,j} |
\le {2\g\over |n|^\t} \right\} ,
\Eq(7.16) $$
where $a_{n,j} = {\rm tr} \, \overline M_{n,j}$ and
$b_{n,j} = \det \overline M_{n,j}$.

Given $n$ the condition $(n,j)\in \O$ implies that
$p_{j}$ can assume at most $\e_{0} n+1$ different values. On a
one-dimensional sphere of radius $R$ there are less than $R$ integer
points, so the number of values $j$ can assume is bounded
proportionally to $n(\e_{0} n+1)$.

Since $\m\in \MM$ we have for $|n| \le n_{0}
(\g_{0}^{2}/\e_{0})^{1/2\t_{0}}$, with some constant $n_{0}$,
$$ \left| \d_{n,j}(\d_{n,j} + p_{j}^{-s} a_{n,j})+p_{j}^{-2s}b_{n,j}
\right| \ge \left( |(D+\m) n- p_{j}-\m| - 2\e_{0} |n| \right)^{2} -
{\rm const.} \e_{0} \ge {\g_{0}^{2}\over 2|n|^{2\t_{0}} } ,
\Eq(7.17) $$ 
so that
$$ \left| \d_{n,j}(\d_{n,j} + p_{j}^{-s} a_{n,j})+p_{j}^{-2s}b_{n,j}
\right| \ge {\g \over |n|^{\t}} ,
\Eq(7.18) $$
provided $\g<\g_{0}^{2}/2$ and $\t>2\t_{0}$.

The measure of each $\III_{n,j}$ can be bounded from above by
$$ {2\g\over |n|^\t } \sup_{\e\in \EEsub_{0}(\g)} \left|
\left( {\der \over \der\e} \left( \d_{n,j}^{2} + p_{j}^{-s} a_{n,j}
\d_{n,j} + p_{j}^{-2s} b_{n,j} \right) \right)^{-1} \right| .
\Eq(7.19) $$

In order to control the derivative we restrict $\e$ to the Cantor set
$$ \EE_{2} = \left\{ \e \in \EE_{0}(\g) : 
\left| \d_{n,j} (2n +p_{j}^{-s}a_{n,j}'(\e))+ np_{j}^{-s}a_{n,j}+
p_{j}^{-2s}b_{n,j}'(\e) \right| \ge {\g \over |n|^{\t_{2}} }
\hbox{ for all } (n,j) \in
\O \right\} ,
\Eq(7.20) $$
with $a_{n,j}'(\e)=\der a_{n,j}(\e)/\der \e$ and
$b_{n,j}'(\e)=\der b_{n,j}(\e)/\der \e$. On this set we have
(recall that $p_{j}$ is bounded proportionally to $|n|$)
$$ \sum_{(n,j) \in \O} {\rm meas} \left( \III_{n,j} \right) \le
{\rm const.} \sum_{n \ge n_{0} (\g/\e_{0})^{1/{2\t_{0}}}}
{n + \e_{0}n^{2} \over |n|^{\t-\t_{2}}} \le {\rm const.}
\e_{0}^{(\t-\t_{2}-2)/2\t_{0}} ,
\Eq(7.21) $$
provided $\t> \t_{2}+3$. Hence ${\rm meas}( \III_{n,j})$ is small
with respect to $\e_{0}$ provided $\t>2\t_{0}+\t_{2}+2$.

Finally let us study the measure of $\EE_{2}$. The bounds \equ(6.2) --
and their proofs to deal with the second derivatives -- imply
$$ \eqalign{
& \left| a_{n,j} \right| , \left| b_{n,j} \right| \le C\e_{0} ,
\qquad \left| a'_{n,j} \right|, \left| b'_{n,j} \right| \le
C \left( 1+\e_{0}|n|^{1+s_{2}} \right) , \cr
& \left| a''_{n,j} \right| , \left| b''_{n,j} \right|
\le C \left( 1+\e_{0}|n|^{2+2s_{2}} \right) , \cr}
\Eq(7.22) $$
for some constant $C$.

Let us call $\III^1_{n,j}$ the complement of $\EE_{2}$ with
respect to $\EE_{0}(\g)$ at fixed $(n,j)\in\O$. As in estimating the
set $\III_{n,j}$ in \equ(7.16) we can restrict the analysis
to the values of $n$ such that $n>n_{1}(\g_{0}/\e_{0})^{1/2\t_{0}}$,
possibly with a constant $n_{1}$ different from $n_{0}$. Then we
need a lower bound on the derivative, which gives
$$ \left| |n|(2n +a_{n,j}' p_{j}^{-s})+
\d_{n,j} a_{n,j}''p_{j}^{-s}+b_{n,j}''p_{j}^{-2s} \right|
\ge {n^{2} \over 2} ,
\Eq(7.23) $$
(recall that $\d_{n,j}<1/2$). Hence we get
$$ \sum_{(n,j)\in \O} {\rm meas} \left( \III^1_{n,j} \right) \le
{\rm const.} \sum_{n \ge n_{1} (\g/\e_{0})^{1/2\t_{0}}}
{(n+\e_{0}n^{2})\over |n|^{\t_{2}-2}} \le
{\rm const.} \e_{0}^{(\t_{2}-4)/2\t_{0}} ,
\Eq(7.24) $$
provided $\t_{2}>5$. Again the measure is small
with respect to $\e_{0}$ provided $\t_{2}>2\t_{0}+4$.
For $\t_{0}>1$ this gives $\t_{2}>6$
and therefore $\t>2\t_{0}+\t_{2}+2>10$.

\*

\0{\bf Remark.} The argument given above applies only when $D=2$,
because only in such a case the matrices $M_{n,j}$ are of
finite $n$-independent size (cf. Lemma \secc(L2)).
A generalisation to the case $D>2$ should require some further work.

\*\*
\section(8,Proof of Theorems 2, 3 and 4)

\0Let us now consider \equ(1.9) with $\m=0$, under the conditions
\equ(1.3) if $D=2$ and both \equ(1.3) and \equ(1.4) if $D\geq 3$. 
Note that for $\m=0$ one has $\o=D-\e$.

The $Q$ subspace is infinite-dimensional, namely \equ(1.13))
is replaced by
$$ Q:=\{ (n,m)\in \NNN\times \ZZZ^D: \; Dn = |m|^2\},
\Eq(8.1) $$
so that $Q$ contains as many elements as the set of $m\in \ZZZ^D$
such that $|m|^2/D\in \NNN$. 

As in \cita{GP} our strategy will be as follows: first, we shall find
a finite-dimensional solution of the bifurcation equation, hence we
shall prove that it is non-degenerate in $Q$ and eventually we
shall solve both the $P$ and $Q$ equations iteratively.
 
A further difficulty comes from the separation of the resonant sites.
Indeed the conditions \equ(2.1) and \equ(2.2) are fulfilled now
only for those $(n,p)$ such that $Dn\neq p$.
This implies that Lemma \secc(L1) does not hold: given $p_i^{s_0}
|\o n_i- p_i|\leq \g/2$ for $i=1,2$ it is possible that $D(n_1-n_2)=
p_1-p_2$ and in such case we have at most $|p_1-p_2|\le \g/
\e p_2^{s_0}$, which in general provides no separation at all.
Hence we cannot use anymore the second Melnikov conditions.

We then replace Lemma \secc(L2) a by more general result
(cf. Lemma \secc(L20) below), due to Bourgain; consequently we deal
with a more complicated renormalised $P$ equations.

\*\*
\0{\bf 8.1. The $\QQQQ$ equations}\*

\0In \cita{GP} we considered the one-dimensional case and used the
integrable cubic term in order to prove the existence of
finite-dimensional subsets of $Q$ such that there exists a solution
of the bifurcation equation with support on those sets.

In order to extend this result to $D\ge 2$ we start by
considering  \equ(1.9) projected on the $Q$ subspace. We set
$u_{n,m}=q_{m}$ if $(n,m)\in Q$, so that the $Q$ equations become
$$ |m|^{2(1+s)} D^{-1} q_{m} =
\sum_{m_1+m_2-m_3=m \atop n_1+n_2-n_3= |m|^2/D}
u_{n_1,m_1}u_{n_2,m_2}u_{n_3,m_3} . $$
Setting $q_{m}=a_{m}+Q_{m}$, with $Q_{m}=O(\h)$,
the leading order provides a relatively simple equation,
as shown by Bourgain in \cita{Bo4}:
$$ |m|^{2(1+s)} D^{-1}
a_m= \sum_{m_1,m_2,m_3\atop {m_1+m_2-m_3=m\atop 
\la m_1-m_3,m_2-m_3\ra =0}} a_{m_1}a_{m_2}a_{m_3} ,
\Eq(8.2)$$
which will be called the {\it bifurcation equation}.
One can easily find finite-dimensional sets $\MMM$ such that\\
(i) if $m\in\MMM$ then $S_{i}(m)\in\MMM$ $\forall i =1,\ldots,D$
($S_{i}$ is defined in \equ(1.12)),\\
(ii) if $m_1,m_2,m_3\in \MMM$ and $\la m_1-m_3,m_2-m_3\ra =0$,
then $m_1+m_2-m_3 \in \MMM$.

\*

\0{\bf Remarks.}
(1) Condition (i) implies that $\MMM$ is completed described
by its intersection $\MMM_{+}$ with $\ZZZ_{+}^{D}$.\\
(2) Clearly \equ(8.2) admits a solution with support on sets
respecting (i) and (ii) above. An example is as follows.
For all $r$ the set $\MMM_+(r):=\{m\in \ZZZ_+^D: |m|= r\}$ is a
finite-dimensional set on which \equ(8.2) is closed.\\
(3) We look for a solution of \equ(8.2) which satisfies
the Dirichlet boundary conditions. Hence we study \equ(8.2)
as an equation for $a_{m}$ with $m\in\MMM_{+}$. 

\*

Finding non-trivial solutions of \equ(1.9) by starting
from solutions of the bifurcation equation like those of the
example may however be complicated,
so we shall prove the existence of solutions under
the following, more restrictive, conditions.

\*
\0\lm(L17) {\it There exist finite sets $\MMM_+\subset \ZZZ^D_+$ such
that $|m|^2$ is divided by $D$ for all $m\in \MMM_+$ and \equ(8.2)
is equivalent to
$$ \cases{
|m|^{2(1+s)}D^{-1} -2^{D+1}A^2 - ( 3^{D}-2^{D+1}) a_{m}^{2} =0 , &
$a_m \in \MMM_+$ , \cr
a_{m} = 0 , & $a_m \in \ZZZ_+^D \setminus \MMM_+$, \cr}
\Eq(8.3) $$
with $A^{2}:=\sum_{m\in\MMM_+}a_m^2$.}
\*

\0{\bf Proof.} The idea is to choose the $m\in \MMM_+$ so that
$|m|^{2} \in D \NNN$,  \equ(8.3) is equivalent to \equ(8.2) and has a
non-trivial solution. We choose $\MMM_+$ so that the following
conditions are fulfilled:\\
\0(a) setting $N := |\MMM_+|$ and  $\min_{m\in\MMM_{+}}|m|=|m_1|$
(this only implies a reordering of the elements of $\MMM_+$), we impose
$$ 2^{D+1}\sum_{m\in \MMM_+\setminus \{m_{1}\}}|m|^{2+2s}  \le
\left( 3^{D}+2^{D+1}(N-2) \right) |m_{1}|^{2+2s} ;
\Eq(8.4) $$
\0(b) the identity $\la m_1-m_3,m_2-m_3\ra =0$ can be verified only if
either $m_1-m_3=0$ or $m_2-m_3=0$ or $|(m_1)_i|=|(m_2)_i|=|(m_3)_i|$
for all $i=1,\ldots,D$ ($(m_j)_i$ is the $i$-th component
of the vector $m_j$).

An easy calculation shows that under conditions (b)
equation \equ(8.2) assumes the form
$$ a_m \left( |m|^{2(1+s)}D^{-1} -2^{D+1}A^{2} -
(3^{D}-2^{D+1}) a_{m}^{2} \right) =0 ,
\Eq(8.5) $$
and hence is equivalent to \equ(8.3). Now, in order to find a
non-trivial solution to \equ(8.5) we must impose
$$ |m_{1}|^{2+2s} = \min_{m\in \MMM_+} |m|^{2+2s} \ge 2^{D+1} A^{2} D,
\Eq(8.6) $$
with $A$ determined by
$$ D \left( 2^{D+1}(N-1)+ 3^{D} \right) A^2= M , \qquad
N = |\MMM_+| , \qquad M := \sum_{m\in\MMM_+} |m|^{2+2s} .
\Eq(8.7) $$
As in the one-dimensional case \cita{GP}, if we fix $N$ then
\equ(8.6) is equivalent to condition (a), i.e. \equ(8.4),
which is an upper bound on the moduli of
the remaining  $m_i\in \MMM_+ \setminus \{m_{1}\}$.
Then there exist sets of the type described above at least for $N=1$.

To complete the proof (for all $N\in\NNN$) we have still to show
that sets $\MMM_{+}$ verifying the conditions (a) and (b) exist. The
existence of sets with $N=1$ is trivial,
an iterative method of construction for any $N$ is then
provided in Appendix \secc(A3). \qed
 
\*

\0{\bf Remark.} The compatibility condition \equ(8.4) requires
for the harmonics of the periodic solution to be large enough,
and not too spaced from each other. Therefore, once we have proved
that the solutions of the bifurcation equation can be continued for
$\e\neq0$, we can interpret the corresponding periodic solutions as
perturbed wave packets. The same result was found in $D=1$ in \cita{GP}.

\*

We have proved that the bifurcation equation admits a
non-trivial solution
$$ q^{(0)}(x,t) = \sum_{m\in \MMM_+} q^{(0)}_{m} {\rm e}^{i{|m|^2
\over D}t} \left( 2 i \right)^{D} \prod_{i=1}^D \sin(m_ix_i),
\Eq(8.9) $$
with $q^{(0)}_{m}=a_{m}$ for $m\in\ZZZ_{+}^{D}$ and extended to
all $\ZZZ^{D}$ by imposing the Dirichlet boundary conditions.

We can set $q_m=q^{(0)}_m+Q_m$ for all $m\in \ZZZ^D$ and
split the $Q$ equation in a bifurcation equation \equ(8.3)
and a recursive linear equation for $Q_m$:
$$ |m|^{2+2s}D^{-1} Q_m- 2\!\!\!\!\!\!\!\!\!
\sum_{m_1,m_2,m_3\atop {m_1+m_2-m_3=m\atop 
\la m_1-m_3,m_2-m_3\ra =0}}\!\!\!\!\!\!\!\!\!\!\!
Q_{m_1}q^{(0)}_{m_2}q^{(0)}_{m_3} - \!\!\!\!\!\!\!
\sum_{m_1,m_2,m_3\atop {m_1+m_2-m_3=m\atop 
\la m_1-m_3,m_2-m_3\ra =0}}\!\!\!\!\!\!\!\!\!\!\!\!\!\!\!
q^{(0)}_{m_1}q^{(0)}_{m_2} Q_{m_3}= \!\!\!\!\!\!\!\!\!
{\mathop{\sum}_{m_1+m_2-m_3=m \atop
n_1+n_2-n_3= |m|^2/D}}^{\hskip-0.9truecm*}
u_{n_1,m_1}u_{n_2,m_2}u_{n_3,m_3}
\Eq(8.10) $$
where for all $(n,m)\in Q$ one has $u_{n,m}\equiv q_m$ and $*$
in the last sum means that the sum is restricted to the triples
$(n_i,m_i)$ such that if at least two of $u_{n_{i},m_{i}}$ are
$q^{(0)}_{m_{i}}$ then the label $(n,m)$ of the third one
must not belong to $Q$.

By using once more the Dirichlet boundary conditions, we can see
\equ(8.10) as an equation for the coefficients $Q_{m}$ with
$m\in\ZZZ_{+}^{D}$. In particular the left hand side yields
an infinite-dimensional matrix $J$ acting on $\ZZZ_{+}^{D}$. 
We need to invert this matrix.

\*
\lm(L18) {\it For all $D$ and for all choices of $\MMM_{+}$ as in
Lemma \secc(L17), one has that $J$ is a block-diagonal matrix,
with finite dimensional blocks, whose sizes are bounded from above
by some constant $M_{1}$ depending only on $D$ and $\MMM_{+}$.}
\*

The result above is trivial for $D=2$ and requires
some work for $D>2$, see Appendix \secc(A4).
In any case it is not enough to
ensure that the matrix $J$ is invertible. The following discussion
shows that at least for $N=1$ (and any $D\ge 2$) and for $N>4$ and $D=2$
 there exist sets $\MMM_{+}$ such
that the matrix $J$ is invertible outside a discrete set of values of $s$.

We can write  $J={\rm diag} \{|m|^{2+2s}/D -8 A \}  + Y$,
where $A$ is defined in \equ(8.7) and with $|Y|_{\io}$
bounded by a constant depending only on $D$ and $\MMM_+$.
Therefore for $M_{0}$ large enough we can write $J$ as
$$ J= \left( \matrix{ J_{1,1} &0 \cr
0 &J_{2,2} \cr} \right) , $$
where $J_{1,1}$ is a $M_{0}\times M_{0}$ matrix, and $J_{2,2}$ is --
by the definition of $M_0$ -- invertible. 

To ensure the invertibility of $J_{1,1}$ 
we notice that ${\rm det} J_{1,1}=0$  is an  analytic equation for the
parameter $s$, and therefore is either identically satisfied or  has
only a denumerable set of solutions with no accumulation points.
For all $s$ outside such denumerable set $J$ is invertible.    

Proving that for a given $\MMM_+$ the function ${\rm det} J_{1,1}$ is
not identically zero can be however quite lengthly. 

For $N=1$ and $\MMM_+= \{V\equiv (1,\dots,1)\}$ the Dirichlet boundary
conditions imply that we only need to consider those $m\in \ZZZ^D_+$ with 
strictly positive components. For all such $m$ either $m= V$
or $|m|^2>D$. This implies that $J_{1,1}$ has two diagonal
blocks: a $1\times 1$ block 
involving $\MMM_+$ and a block involving $m$ such that $|m|^2>D$.
The first block is trivially found to be non-zero. In the second block
the off-diagonal entries all depend linearly on $D^{2s}$, and
for all $m$ the diagonal entry with index $m$ is $|m|^{2(1+s)}/D$
plus a term depending linearly on $D^{2s}$:
therefore in the limit $s\to \io$ this block is invertible.
Hence det$J_{1,1}=0$ is not an identity in $s$.
 
If $N>1$ we restrict our attention to the case $D=2$ where we can describe
the matrix $J_{1,1}$ with sufficient precision. We have the following 

\*
\lm(L19) {\it For $D=2$ and $N>4$ consider $\MMM_+$ as
a point in $\CCC^{2N}$.\\
(i) The set of points $\MMM_+$ which either  do not respect Lemma
\secc(L17) or are such that ${\rm det} J_{1,1}=0$ identically in
$s$ is contained in a proper algebraic variety $\WW$.\\
(ii) Provided that $|m_1|$ is large enough one can always find integer
point which do not belong to $\WW$ and respect \equ(8.4) for all $s$
in some open interval.}
\*

The proof is in Appendix \secc(A5).

Therefore the forthcoming analysis applies without any further assumption
for $D=2$ and $N>4$, whereas one must assume that $J_{1,1}$ is invertible
to apply it to the other cases. Of course, given a set $\MMM_{+}$ verifying
the conditions of Lemma \secc(L17) one can check, through a finite number
of operations, whether $J_{1,1}$ is invertible, and, if it is, then the
analysis below ensure the existence of periodic solutions.

\*\*
\0{\bf 8.2. Renormalised $\PPPP$ equation}\*

\0The following Lemma (Bourgain lemma) will play a fundamental role
in the forthcoming discussion. A proof is provided in Appendix \secc(A6).

\*
\0\lm(L20) {\it For all sufficiently small $\a$
we can partition $\ZZZ^D=\cup_{j\in \NNN }\D_j$ so that, setting
$$ p_j= \min_{m\in \D_j} |m|^2, \qquad
\Phi(m)= (m,|m|^2) ,
\Eq(8.13) $$
there exist $j$-independent constants $C_{1}$ and $C_{2}$ such that
$$ |\D_j| \le C_{1} p_j^\a ,
\qquad {\rm dist}(\Phi(\D_i),\Phi(\D_j)) \ge
C_{2} \min \{ p_{i}^{\b} , p_{j}^{\b} \} ,
\qquad {\rm diam}(\D_{j}) < C_{1} C_{2} p_{j}^{\a+\b} ,
\Eq(8.14) $$
with $\b=\a/(1+2^{D-1}D!(D+1)!)D$.}
\*

\0{\bf Remarks.}
(1) For fixed $\e$, $\o n-|m|^2$ can be small only if
$n$ is the integer nearest to $|m|^2/\o$.\\
(2) For any $(m_{1},n_{1})$ and $(n_{2},m_{2})$ such that
$m_{1}\in \D_{j}$, $m_2\in \D_{j'}$ for $j'\neq j$, and $n_{i}$
is the integer nearest to $|m_{i}|^{2}/\o$, $i=1,2$, one has
$$ |m_{1} - m_{2}| + \left| n_{1} - n_{2} \right|
\ge C_{3} \min \{ p_{j}^{\b}, p_{j'}^{\b} \}
\Eq(8.15) $$
for some constant $C_{3}$ independent of $\o$.\\
(3) As in Lemma \secc(L2) also here one could prove that in fact
${\rm diam}(\Delta_{j}) < {\rm const.} p_{j}^{\a/D}$;
see Appendix \secc(A6) for details.

\*
\0\sub(D13)
{\it We call $\CC_{j}$ the sets of $(n,m)\in \ZZZ \times \ZZZ^{D}$
such that $m\in \D_j$, $D n\neq |m|^2$  and $-1/2+ (D-\e_{0}) n
\le |m|^2 \le D n + 1/2$.  
We set $\d_{n,m}= -\o n+|m|^2$ and $d_j= |\CC_{j}|$, and define
the $d_j$-dimensional vectors and the $d_{j} \times d_{j}$ matrices
$$ U_{j}=\{u_{n,m}\}_{(n,m)\in \CC_{j}} ,
\quad \DDDD_j={\rm diag}
\left\{ \left( {|m|^2\over p_{j}} \right)^{s}\d_{n,m}
\right\}_{(n,m)\in \CC_{j}} \!\!\!\!\!\!\!\!\! \!\!\!\!\!\!\!\!\! ,
\qquad \quad \widehat \chi_{1,j} = {\rm diag}
\left\{ \sqrt{\bar\chi_1(\d_{n,m})} \right\}_{(n,m)\in \CC_{j}}
\Eq(8.16) $$
parameterised by $j\in \NNN$.}
\* 

\0{\bf Remark.} Notice that for each pair $(n,m),(n',m')\in \CC_{j}$
we have $|(n,m)-(n',m') |\le  C( \e_{0}p_{j}/D + p_{j}^{2\a})$
for a suitable constant $C$.
\*
We define the renormalised $P$ equations
$$ \left\{
\eqalign{
& u_{n,m} = \h
{f_{n,m} \over |m|^{2s}\d_{n,m} } ,
\qquad \qquad \qquad \qquad \qquad \qquad
(n,m) \notin \bigcup_{j\in \NNN} \CC_{j},
\quad D n \neq |m|^{2} , \cr
& p_{j}^{s} \left( \DDDD_j + p_{j}^{-s}
\widehat M_{j} \right) U_j =
\h F_j + L_jU_j , \qquad \, j\in\NNN, \cr} \right.
\Eq(8.17) $$
where $\widehat M_{j}= \widehat \chi_{1,j} M_{j} \widehat \chi_{1,j}$,
and the parameter $\h$ and the counterterms $L_{j}$ will have
to satisfy eventually the identities
$$ \h = \e , \qquad \widehat M_{j} = L_{j}
\Eq(8.18) $$
for all $j\in\NNN$.

\*

\0{\bf Remark.} We note that $d_j$ can be as large as
$O(\e_{0}p_{j}^{1+\a})$, hence can be large with respect to $p_j$.
However for given $\e$ the matrix $A_{j}= \DDDD_j + p_{j}^{-s}
\widehat M_j $ is diagonal apart
from a $p_{j}^{\a} \times p_{j}^{\a}$ ($\e$-depending) block.
This implies that the matrix $A_{j}$ has at most $p_{j}^{\a}$
eigenvalues which are different from $|m|^{2s} \d_{n,m}$.
This can be proved as follows. Consider the entry $A_{j}(a,b)$,
with $a,b\in \CC_{j}$, with $a=(n_{1},m_{1})$ and $b=(n_{2},m_{2})$.
The non-diagonal part can be non-zero only if $\bar \chi_{1}
(\d_{n_{1},m_{1}}) \bar\chi_{1}(\d_{n_{2},m_{2}})\,M(a,b) \neq 0$,
which requires $|\d_{n_{i},m_{i}}|\le \g/4$ for $i=1,2$.
Therefore for fixed $\e$, $m_{1}$ and $m_{2}$ one has only
one possible value for each $n_{i}$, i.e. the integer
closest to $\o^{-1}|m_{i}|^2$. This proves the assertion
because $|\D_{j}|\le C_{1} p_{j}^{\a}$ and
for all $(n,m) \in \CC_{j}$ one has $m\in \D_{j}$.

\*

Definition \secc(D1) and Lemma \secc(L3) still hold, with $\ZZZ^{dD}$
replaced with $\ZZZ^{d(D+1)}$ in the definitions of $\AAA(\ul{m})$.
Definitions \secc(D2) (i)-(ii) can be maintained with $(n,j)$
replaced by $j$, while (iii) becomes
$$ x_{j} = \Vert \widehat \chi_{1,j}
(\DDDD_{j} + p_{j}^{-s} \widehat M_{j} )^{-1}
\widehat \chi_{1,j} \Vert^{-1} .
\Eq(8.19) $$
Finally there is no parameter $s_{2}$. Equivalently we can set
$s_{2}=0$, which leads to identify $y_{n,m}$ with $\d_{n,m}$
(cf. \equ(2.8)): this explains why there is no need to
introduce the further parameters $y_{n,m}$.

The main Propositions 1 and 2 in Section 2.4 still hold with
the following changes.
\vskip.2truecm
\01. $(n,j)\in \ZZZ\times \NNN$ (or $\O$) has to be always
replaced with $j\in \NNN$,
\vskip.1truecm
\02. In Proposition 1, $q$ (i.e. the solution of the $Q$ equation)
is not a parameter any more: it is substituted with the
solution, say  $q^{(0)}$, of the bifurcation equation \equ(8.2),
whose Fourier coefficients can be incorporated in the list of
positive constants given at the beginning of the statement.
\vskip.1truecm
\03. In Proposition 1 (i) the bound \equ(2.20) becomes
$$ \left| u_{n,m}(\h,M,\e\right) | \le
K_{0} |\h| {\rm e}^{-\s(|n|^{1/4}+|m|^{1/4})},
\Eq(8.20) $$ 
for some constant $K_{0}$,
namely we have only sub-analyticity  in space and time.
\vskip.1truecm
\04. In Proposition 1 (v) one must replace $s_{2}$ with $s$
in the first line of \equ(2.23) and in \equ(2.26),
and ${\rm e}^{-\s|m_a-m_b|^{\r}} $ with 
${\rm e}^{-\s(|(n_a,m_a)-(n_b,m_b)|^{\r}}$
in the second line of \equ(2.23), for a suitable constant $\r$.

\*\*
\0{\bf 8.3. Multiscale analysis}\*

\0The multiscale analysis follows in essence the same ideas
as in the previous sections,
but there are a few changes, that we discuss here.
It turns out to be more convenient to replace the functions
$\chi_{h}(x)$ with new functions $\tilde \chi_{h}(x)=\chi_{h}(32x)$,
in order to have $\tilde\chi_{-1}(x_{j})=1$ when
$\bar\chi_{1}(\d_{n,m})\neq 1$ for all $(n,m)\in \CC_{j}$.
This only provides an extra factor $32$ in the estimates.
For notational simplicity in the following we shall drop the tilde.

Let us call $A_{j}= \DDDD_{j}+p_{j}^{-s}\widehat M_{j}$. Note that
$$ 1= \bar\chi_{1}(\d_{n,m}) + \bar\chi_{0}(\d_{n,m}) +
\bar\chi_{-1}(\d_{n,m}) \qquad \forall (n,m)\in \CC_{j} .
\Eq(8.21) $$
Introduce a {\it block multi-index} $\vb$, defined as
a $d_{j}$-dimensional vector with components $\bb (a) \in \{1,0,-1\}$,
and set
$$ \bar\chi_{j,\vb} = \prod_{a=1}^{d_{j}}
\bar\chi_{\bb(a)} (\d_{n(a),m(a)}) .
\Eq(8.22) $$
For any $\vb$ we can consider the permutation $\pi_{\vb}$
which reorders $(\bb(1),\ldots,\bb(d_{j}))$ into
$(\bb_{\pi_{\vb}}(1),\ldots,\bb_{\pi_{\vb}}(d_{j}))$ in such a way
that the first $N_{1}$ elements are 1,
the following $N_{2}$ elements are 0, and the last
$N_{3}=d_{j}-N$, with $N=N_{1}+N_{2}$, elements are $-1$.
The permutation $\pi_{\vb}$ induces a permutation matrix
$P_{\vb}$ such that $P_{\vb} A_{j} P_{\vb}^{-1}$
can be written in the block form
$$ P_{\vb} A_{j} P_{\vb}^{-1} =
\left(
\matrix{A_{1,1} & A_{1,2} & A_{1,3} \cr
A_{1,2}^{T} & A_{2,2} & A_{2,3} \cr
A_{1,3}^{T} & A_{2,3}^{T} & A_{3,3} \cr} \right) ,
\Eq(8.23) $$
where the block $A_{1,1}$, $A_{2,2}$ and $A_{3,3}$ contain
all the entries $A_{j}(a,b)$ with $\bb(a)=\bb(b)=1$,
with $\bb(a)=\bb(b)=0$ and $\bb(a)=\bb(b)=-1$, respectively,
while the non-diagonal blocks are defined consequently.

Then for all $\vb$ such that $\bar\chi_{j,\vb}\neq0$ we can write
$$ A_{j}= P_{\vb} \left(
\matrix{ A_{1,1} & A_{1,2} & 0 \cr
A_{1,2}^{T} & A_{2,2} & 0 \cr
0 & 0 & A_{3,3} \cr} \right) P_{\vb}^{-1}
\Eq(8.24) $$
where we have used that if $\bar\chi_{j,\vb}\neq0$ then
the blocks $A_{1,3}$ and $A_{2,3}$ are zero. Furthermore,
for the same reason, the block $A_{3,3}$ is a diagonal matrix.
Note that $N \le C_{1} p_{j}^{\a}$ by the Remark after \equ(8.18).

The first $N\times N$ block of $A_{j}$ in general is not
block-diagonal, but it can be transformed into a block-diagonal
matrix. Indeed, we have
$$ A_{j} = S_{j,\vb} \widetilde A_{j,\vb} S_{j,\vb}^{T} ,
\qquad \widetilde A_{j,\vb} =
\left( \matrix{ \widetilde A_{1,1} & 0 & 0 \cr
0 & A_{2,2} & 0 \cr
0 & 0 & A_{3,3} \cr} \right) 
\qquad S_{j,\vb} = P_{\vb} \left(
\matrix{ I & B & 0 \cr 0 & I & 0 \cr 0 & 0 & I \cr} \right)
 ,
\Eq(8.25) $$
where
$$ \widetilde A_{1,1} = A_{1,1}-A_{1,2} A_{2,2}^{-1} A_{1,2}^{T} ,
\qquad B = A_{1,2} A_{2,2}^{-1} ,
\Eq(8.26) $$
while $I$ and $0$ are the identity and the null matrix
(in the correct spaces). Of course also the matrices $A_{i,j}$
depend on $\vb$ even if we are not making explicit such a dependence.

The invertibility of $A_{2,2}$ is ensured by the condition
$\bb(a)=0$ for the indices $a=N_{1}+1,\ldots,N$. The inverse
$A_{2,2}^{-1}$ can by bounded proportionally to $1/\g$
in the operator norm.
Then also $A_{j}$ can be inverted provided $\widetilde A_{1,1}$
is invertible, i.e. provided $\det \tilde A_{1,1} \neq 0$.
Hence in the following we shall assume that this is the case
(and we shall check that this holds true whenever it appears;
see in particular \equ(8.32) below).

Hence for all $\vb$ such that $\bar\chi_{j,\vb}\neq0$ we can write
$$ A_{j}^{-1} = S_{j,\vb}^{-T} \widetilde A_{j,\vb}^{-1} S_{j,\vb}^{-1} ,
\Eq(8.27) $$
and set
$$ \eqalign{
& \GGG_{j,\vb,1} = p_{j}^{-s} S_{j,\vb}^{-T}
\left( \matrix{ \widetilde A_{1,1}^{-1} & 0 & 0 \cr
0 & 0 & 0 \cr 0 & 0 & 0 \cr} \right) S_{j,\vb}^{-1} , \cr
& \GGG_{j,\vb,0} = p_{j}^{-s} S_{j,\vb}^{-T}
\left( \matrix{ 0 & 0 & 0 \cr
0 & A_{2,2}^{-1} & 0 \cr 0 & 0 & 0 \cr} \right) S_{j,\vb}^{-1} , \qquad
\GGG_{j,\vb,-1} = p_{j}^{-s} S_{j,\vb}^{-T}
\left( \matrix{ 0 & 0 & 0 \cr
0 & 0 & 0 \cr 0 & 0 & A_{3,3}^{-1} \cr} \right) S_{j,\vb}^{-1}  \cr}
\Eq(8.28) $$
so that \equ(8.27) gives
$$ p_{j}^{-s} A_{j}^{-1} =
\GGG_{j,\vb,-1} + \GGG_{j,\vb,0} + \GGG_{j,\vb,1}
\Eq(8.29) $$
for all $\vb$ such that $\bar\chi_{j,\vb}\neq 0$.
We can define $\GGG_{j,\vb,i}$ also for $\vb$ such that $\bar
\chi_{j,\vb}= 0$, simply by setting $\GGG_{j,\vb,i}=0$ for such $\vb$.
Then we define the propagators
$$ G_{j,\vb,i,h} = \cases{
\bar\chi_{j,\vb} \, \chi_{h}(x_{j}) \, \GGG_{j,\vb,1} , &
if $i=1$ and $\chi_{h}(x_{j}) \neq 0$, \cr
& \cr
\bar\chi_{j,\vb} \, \GGG_{j,\vb,i} , &
if $i=0,-1$ and $h=-1$, \cr
& \cr
0 , & otherwise , \cr}
\Eq(8.30) $$
so that we obtain
$$ \eqalign{
p_{j}^{-s} A_{j}^{-1} & =
p_{j}^{-s} \sum_{\vb} \bar\chi_{j,\vb} A_{j}^{-1} =
\sum_{\vb} \bar\chi_{j,\vb} \Big[ \left(
\GGG_{j,\vb,-1} + \GGG_{j,\vb,0} \right) +
\sum_{h=-1}^{\io} \chi_{h}(x_{j}) \, \GGG_{j,\vb,1} \Big] \cr
& = \sum_{\vb} \sum_{i=-1,0,1} \sum_{h=-1}^{\io} G_{j,\vb,i,h} , \cr}
\Eq(8.31) $$
which provides the multiscale decomposition.

\*

\0{\bf Remark.}
Only the propagator $G_{j,\vb,1,h}$ can produce small divisors,
because the diagonal propagator $G_{j,\vb,-1,-1}$ and
the non-diagonal propagator $G_{j,\vb,0,-1}$ have denominators
which are not really small.
We can bound $|G_{j,\vb,i,-1}|_{\s}$ for $i=-1,0$
by using a Neumann expansion, since by definition in the corresponding
blocks one has $|\d_{n,m}| \ge \g/8$ and $|M_{j}|_{\s}\le C \e_0,$.

\*

Hence we can bound the propagators as
$$ \left| G_{j,\vb,i,-1} \right|_{\s} \le C
\g^{-1} p_{j}^{-s} , \quad  i=0,-1 ,
\qquad \left| G_{j,\vb,1,h} \right|_{\io} \le
2^{h} C \g^{-1} p_{j}^{-s+\a} ,
\Eq(8.32) $$
for all $j\in\NNN$.

Recall that we are assuming $|J^{-1}|_{\s} \le C$ for
some $s$-dependent constant $C$.
 
We write the counterterms as
$$ L_{j}= \chi_{h}(x_{j}) \sum_{\vb}
\bar\chi_{j,\vb} L_{j,\vb,h} ,
\Eq(8.33) $$
where by definition $L_{j,\vb,h}(a,b)=0$ if either $\bb(a)=-1$
or $\bb(b)=-1$.

With this modifications to \equ(3.8) the multiscale expansion
follows as in Section 3.1, with $j=(n,m)$:
$$ U^\ka_{j} = \sum_{i=-1,0,1} \sum_{\vb}
\sum_{h=-1}^{\infty} U^\ka_{j,\vb,i,h} ,
\Eq(8.34) $$
with
$$ \left\{ \eqalign{
& u^\ka_{n,m} =  {f^\ka_{n,m}\over |m|^{2s}\d_{n,m} } ,
\qquad\qquad (n,m) \notin \bigcup_{j\in \NNN} \CC_{j} ,
\quad D n \neq |m|^{2} , \cr
& U^\ka_{j,\vb,i,h} =
G_{j,\vb,i,h} F^\ka_{j}+ \d_{i,1} G_{j,\vb,1,h}
\Big( \sum_{h_{1}=-1}^{\io}\sum_{\vb_{1}\neq \vec 0}
\sum_{i_{1}=0,1,-1} \sum_{r=1}^{k-1}L^{(r)}_{j,\vb,h}
U^{(k-r)}_{j,\vb_{1},i_{1},h_{1}} \Big) , \qquad
j \in \NNN, \cr
& u^{(k)}_{n,m} = q_m^\ka= J^{-1} \sum_{k_{1}+k_{2}+k_{3}=k}
{\mathop{\sum}_{m_1+m_2-m_3=m \atop
n_1+n_2-n_3= |m|^2/D}}^{\hskip-0.8truecm*} \quad
u_{n_1,m_1}^{(k_{1})} u_{n_2,m_2}^{(k_{2})} u_{n_3,m_3}^{(k_{3})} ,
\qquad Dn =|m|^{2} , \cr} \right.
\Eq(8.35) $$
where $*$ has the same meaning as in \equ(8.10).

\*\*
\0{\bf 8.4. Tree expansion}\*

\0We only give the differences with respect to Section 3.2.
\vskip.2truecm
\0(2) One has $(n_{v},m_{v})\in Q$ and the node factor is
$\h_{v}= q^{(0)}_{m_{v}}$.
\vskip.1truecm
\0(3)  We add a further label $r,p,q$ to the lines to evidence which
term of \equ(8.35) we are considering. We also associate
with each line $\ell$ a label $j_{\ell}\in \ZZZ_{+}$,
with the constraints $j_{\ell}\in\NNN$ if $\ell$ is a $p$-line
and $j_{\ell}=0$ otherwise.
\vskip.1truecm
\0(4) The momenta are:
$(n_{\ell},m_{\ell}),(n'_{\ell},m'_{\ell})\in \CC_{j_{\ell}}$
for a $p$-line,
$(n_{\ell},m_{\ell}),(n'_{\ell},m'_{\ell})\in Q$,
with $|m_{\ell}-m_{\ell}'| \le M_{1}$,
for a $q$-line, and finally
$(n_{\ell},m_{\ell})=(n'_{\ell},m'_{\ell}) \notin
\cup_{j\in\NNN} \CC_{j} \cup Q$ for an $r$-line.
For a $p$-line the momenta define the labels $a_{\ell},b_{\ell}
\in \{1,\ldots,d_{j}\}$, with $d_{j_{\ell}}=|\CC_{j_{\ell}}|$,
such that $(n_\ell,m_{\ell})= \CC_{j_{\ell}}(a_{\ell})$ and
$(n'_{\ell},m'_{\ell})=\CC_{j_{\ell}}(b_{\ell})$.
For a $q$-line the momenta define $a_{\ell},b_{\ell}$
such that $(n_{\ell},m_{\ell})= Q(a_{\ell})$ and
$(n'_{\ell},m'_{\ell})= Q(b_\ell)$.
\vskip.1truecm
\0(5) Each $p$-line carries also a {\it block label} $\vb_{\ell}$
with components $\bb_{\ell}(a)=-1,0,1$, where $a=1,\ldots,d_{j_{\ell}}$. 
\vskip.1truecm
\0(6) Both $r$-lines and $q$-lines $\ell$ have $i_{\ell}=-1$
and $h_{\ell}=-1$.
\vskip.1truecm
\0(7) One must replace $(n_{\ell},j_{\ell})$ with $j_{\ell}$.
Moreover if two lines $\ell$ and $\ell'$ have $j_{\ell}=j_{\ell'}$
then $|\bb_{\ell}(a)-\bb_{\ell'}(a)| \le 1$ and if $h_{\ell}\neq -1$
then  $\vb_{\ell}\neq \vec 0$ (by the definition of functions $\chi_h$). 
\vskip.1truecm
\0(8) One has $n_{\ell}=n_{v}$ instead of $n_{\ell}=1$ for lines
$\ell$ coming out from end-points.
\vskip.1truecm
\0(9) One must replace $(n_{\ell},j_{\ell})$ with $j_{\ell}$.
\vskip.1truecm
\0(10) Equation \equ(3.15) becomes
$$ n_{\ell}' = \s(\ell_{1}) n_{\ell_{1}} + \s(\ell_{2}) n_{\ell_{2}} +
\s(\ell_3) n_{\ell_{3}} =
\sum_{\ell'\in L(v)} \s(\ell') n_{\ell'} 
\Eq(8.36) $$
(that is $n_{\ell}$ is replaced with $n_{\ell}'$),
while \equ(3.16) does not change.
\vskip.1truecm
\0(11) The propagator $G_{\ell}$ of any line $\ell$ is given by
$g_{\ell}=G_{j_{\ell},\vb_{\ell},i_{\ell},h_{\ell}}(a_{\ell},b_{\ell})$,
as defined in\equ(8.30), if $\ell$ is a $p$-line, while it is given by
$g_{\ell}= J^{-1}(a_\ell,b_\ell)$ if $\ell$ is an $q$-line and
by $g_\ell = 1/\d_{n_{\ell},m_{\ell}}|m_\ell|^{2s}$ if $\ell$
is a $r$-line.
\vskip.1truecm
\0(12) The node factor for $s_{v}=1$ is $\h_{v}=
L^{(k_{v})}_{j_{\ell},\vb_{\ell},h_{\ell}}(a_v,b_v)$.
\vskip.2truecm

The set $\Th^{(k)}_{j}$ is defined as in Definition \secc(D7),
with $j$ instead of $n,m$, by taking into account also the new rules
listed above. This will lead to a tree representation \equ(3.19)
for \equ(8.35), which can be proved as for Lemma \secc(L6).

In Lemma \secc(L5) the estimate $|n_{\ell}|\le B k$ does not
hold any more because there is no longer conservation
of the momenta $n_{\ell}$ (i.e. \equ(3.17) has been replaced with
\equ(8.36)), and all the bounds on the momenta should be modified
into $|n_{\ell}|,|n_{\ell}'|,|m_{\ell}|,|m'_{\ell}|\le B k^{1+4\a}$
for some constant $B$. This can be proved by induction
on the order of the tree. The bound is trivially true to first order.
It is also trivially true if either the root line has $i=-1$ or it is
$q$-line or a $r$-line (one just needs to choose $B$ appropriately).
Suppose now that the root line is a $p$-line with $i \neq -1$:
call $v_{0}$ the node which the root line exits.
If $s_{v_{0}}=3$, call $\th_{1},\th_{2},\th_{3}$ the
sub-trees with root lines $\ell_{1},\ell_{2},\ell_{3}$, respectively,
entering the node $v_{0}$. We have $|(n_{\ell_{i}},m_{\ell_{i}})|
\le k_{i}^{1+4\a}$ by the inductive hypothesis, and by definition
$|(n'_{\ell},m'_{\ell})| \le \sum_{i=1}^{3} B k_{i}^{1+4\a}
\le B (k-1)^{1+4\a}$. Then $|(n_{\ell},m_{\ell})|\le B
(k-1)^{1+4\a} + C_{2} (k-1)^{2\a(1+4\a)}\le B k^{1+4\a}$.
If $s_{v}=1$ the proof is easier. 

\*\*
\0{\bf 8.5. Clusters and resonances}\*

\0Definition \secc(D8) of cluster is unchanged, while
Definition \secc(D9) of resonance becomes as follows.

\*
\0\sub(D14)
{\it We call {\rm 1-resonance} on scale $h\ge 0$ a cluster $T$ of
scale $h(T)=h$ with only one entering line $\ell_{T}$ and
one exiting line $\ell_{T}^{1}$ of scale $h_{T}^{(e)}>h+1$
with $|V(T)|>1$ and such that\\
(i) one has
$$ (a) \quad j_{\ell_{T}^{1}} =  j_{\ell_{T}} ,
\qquad (b) \quad p_{j_{\ell_T}} \ge 2^{(h-2)/\t} ,
\Eq(8.37) , $$
(ii) for all $\ell\in L(T)$ not on the path $\PPP(\ell_{T},
\ell_{T}^{1})$ one has $j_{\ell}\neq j_{\ell_{T}}$.\\
We call {\rm 2-resonance} a set of lines and nodes
which can be obtained from a 1-resonance by setting $i_{\ell_{T}}=0,-1$.
Resonances are defined as the sets which are either 1-resonances
or 2-resonances. Differently from \secc(D9)
we do not include among the resonant lines the lines exiting a
2-resonance.}
\*

Definition \secc(D10) is unchanged provided that we replace $(n,j)$
with $j$, we require $p_{j}\ge 2^{(h-2)/\t}$, we associate with
the node $e$ the labels $(n_{e},m_{e}) \in \CC_{j}$ and with
$\ell_{0}$ the labels $(n_{\ell_{0}},m_{\ell_{0}}) \in \CC_{j}$.

Since we do not have the conservation of the momentum $n$, Lemma
\secc(L7) does not hold in the same form: the bounds
have to be weakened into $|n_{\ell}|,|m_{\ell}|,|n_{\ell}'|,|m_{\ell}'|
\le B k^{1+4\a}$ for the lines $\ell$ not along the
$\PPP(\ell_{e},\ell_{0})$, and $|n_{\ell}|,|m_{\ell}|,|n_{\ell}'|,
|m_{\ell}'| \le B (|n|+k)^{1+4\a}$ for the lines along the path.

\*\*
\0{\bf 8.6. Choice of the counterterms}\*

\0The choice of the counterterm \equ(8.33) is not unique and
therefore is rather delicate.

Resonances produce contributions that make the power series to
diverge. We want to eliminate such divergences with a careful
choice of the counterterms.

The sets $\Th^{(k)}_{R,j}$ and $\RR^{(k)}_{R,h,j}$ are defined
slightly differently with respect to Definition \secc(D11).

\*
\0\sub(D15)
{\it We denote by $\Th_{R,j}^{(k)}$ the set of {\rm renormalised trees}
defined as the trees in  $\Th_{j}^{(k)}$  with the following differences:\\
(i) The trees do not contain any 1-resonance $T$ with
$\vb_{\ell_{T}^{1}}=\vb_{\ell_{T}}$.\\
(ii) If a node $v$ has $s_v=1$ then $\vb_{\ell} \neq \vb_{\ell'}$,
where $\ell$ and $\ell'$ are the lines exiting and entering,
respectively, the node $v$. The factor $\h_{v} = L^{(k_{v})}_{j_{\ell},
\vb_{\ell},h_{\ell}} $ associated with $v$
will be defined in \equ(8.43).\\
(iii) The propagators of any line $\ell$ entering any 1-resonance $T$
(recall that by (i) one has $\vb_{\ell_{T}}^{1} \neq \vb_{\ell_{T}}$),
where $\ell_{T}^{1}=\ell$, is
$$ g_{\ell} = \chi_{h_{\ell}}(x_{j_{\ell}})
\bar \chi_{j_{\ell},\vb_{\ell}}
\left( \GGG_{j_{\ell},\vb_{\ell_{T}},0}(a_{\ell},b_{\ell}) +
\GGG_{j_{\ell},\vb_{\ell_{T}},-1}(a_{\ell},b_{\ell}) -
\GGG_{j_{\ell},\vb_{\ell},0}(a_{\ell},b_{\ell}) -
\GGG_{j_{\ell},\vb_{\ell},-1}(a_{\ell},b_{\ell}) \right) ,
\Eq(8.38) $$
and the same holds for the propagator of any line $\ell$ with $i_\ell=1$
entering a node $v$ with $s_v=1$.
\*
In the same way we define $\R_{R,h,j}^{(k)}$.
We call $\R_{R,h,j}^{(k)}(a,b) $ the set of trees
$\th\in\R_{R,h,j}^{(k)}$ such that the entering line has
$m_{e}= \CC_{j}(a)$ while the root line has $m_{\ell_{0}}'= \CC_{j}(b)$.
Finally we define the sets $\Th_{R}^{(k)}$ and $\R_{R}^{(k)}$
as the sets of trees belonging to $\Th_{R,j}^{(k)}$ for some $j$
and, respectively, to $\R_{R,h,j}^{(k)}$ for some $h,j$.}
\*

By proceeding as in Section 3.5 we introduce the following matrices:
$$ \TT^{(k)}_{j,h}(a,b) = \sum_{h_{1}< h-1}
\sum_{\th\in \R^\ka_{R,j,h_{1}}(a,b)}\Val(\th) .
\Eq(8.40) $$
We use a different symbol for such matrices, as we shall see
that the counterterms will not be identified with the
matrices in \equ(8.38), even if they will be related to them.
We shall see that, by the analog of Lemma \secc(L9),
the matrices $\TT^{(k)}_{j,h}$ are symmetric.

To define the counterterms $L_{j}$ we note that, in order to
cancel at least the 1-resonances, we need  the following
condition:
$$ G_{j,\vb,1,h} \left( L_{j,\vb,h}^{(k)}
+\TT_{j,h}^{(k)} \right)
\GGG_{j,\vb,1} = 0 .
\Eq(8.41) $$
Moreover in order to solve the compatibility condition we need a
solution  $L_{j,\vb,h}(a,b)$ which is proportional to
$\bar\chi_{1}(\d_{n(a),m(a)})\bar\chi_{1}(\d_{n(b),m(b)})$, and
clearly the solution  $L_{j,\vb,h}^{(k)}
+\TT_{j,h}^{(k)}=0$ does not comply with this requirement.
However, since $G_{j,\vb,1,h}$ is not invertible, \equ(8.41) does not
imply $L_{j,\vb,h}^{(k)}=-\TT_{j,h}^{(k)}$; indeed there exists a
solution such that $L_{j,\vb,h}(a,b)\neq 0$ only  if $\bb(a)=\bb(b)=1$.
This solution does not cancel the resonances $T$
with $\vb_{\ell_{T}^{1}} \neq \vb_{\ell_{T}}$,
and does not even touch the 2-resonances.
Nevertheless, if \equ(8.41) holds, we shall see that we are left
only with 2-resonances and partially cancelled 1-resonances,
which admit better bounds (see \equ(8.20)).

By definition $L_{j,\vb,h}^{(k)}(a,b)=0$ if either $\bb(a)$ or
$\bb(b)$ is equal to $-1$. Then \equ(8.41) reduces to the following
equation for the matrix $X=L^{(k)}_{j,\vb,h}+ \TT_{j,h}^{(k)}$:
$$ \left( \matrix{ I \cr 0 \cr 0 \cr} \right)^{T} S_{j,\vb}^{-1} X
S_{j,\vb}^{-T} \left( \matrix{ I \cr 0 \cr 0 \cr}\right) =0
\qquad \Longrightarrow \qquad 
X_{1,1} - \left( BX_{1,2}^{T} + X_{1,2}B^{T} \right) + B X_{2,2}B^{T}=0 ,
\Eq(8.42) $$ where we define: $$P^{-1}_{\vb} X P_{\vb}= \left(\matrix{
X_{1,1} & X_{1,2} & X_{1,3} \cr
X_{1,2}^T & X_{2,2} & X_{2,3} \cr
X_{1,3}^T & X_{2,3}^T & X_{3,3} \cr}\right) $$
In \equ(8.42) there are two matrices which act as free parameters.
A (non-unique) solution is
$$P^{-1}_{\vb} L^{(k)}_{j,\vb,h}P_{\vb} = \left(
\matrix{ 
L^{(k)}_{1,1} & 0 & 0 \cr 0 & 0 & 0 \cr 0 & 0 & 0 \cr} \right) =
\left( \matrix{ I & - B & 0 \cr 0 & 0 & 0 \cr
0 & 0 & 0 \cr} \right) \TT_{j,h}
\left( \matrix{ I & 0 & 0 \cr -B & 0 & 0 \cr 0 & 0 & 0 \cr} \right) .
\Eq(8.43) $$
In this definition, $L_{j,\vb,h}(a,b) \neq 0$ only if
$\bb(a)=\bb(b)=1$, so that $L_{j,\vb,h}$ has the correct
factors $\bar\chi_{1}$. Moreover the 1-resonances with
$\vb_{\ell_{T}^{1}} = \vb_{\ell_{T}}$ are cancelled,
while the 2-resonances are untouched since
$L_{j,\vb,h}(G_{j,\vb,0,-1}+G_{j,\vb,-1,-1})=0$.

Let us now consider a 1-resonance $T$ with
$\vb_{\ell_{T}^{1}} \neq \vb_{\ell_{T}}$. We can write
$$ \eqalignno{
& G_{j,\vb,1,h} \left( L_{j,\vb,h}^{(k)} + \TT_{j,h}^{(k)} \right)
G_{j,\vb_{1},1,h_{1}}
& \eq(8.44) \cr
& \qquad \qquad =
G_{j,\vb,1,h} \left( L_{j,\vb,h}^{(k)} + \TT_{j,h}^{(k)} \right)
\chi_{h_{1}}(x_{j}) \bar\chi_{j,\vb_{1}}
\left( p_{j}^{-s} A_{j}^{-1} - \GGG_{j,\vb_{1},0} -
\GGG_{j,\vb_{1},-1} \right) \cr
& \qquad \qquad = G_{j,\vb,1,h}
\left( L_{j,\vb,h}^{(k)} + \TT_{j,h}^{(k)} \right)
\chi_{h_{1}}(x_{j}) \bar\chi_{j,\vb_{1}}
\left( \GGG_{j,\vb,0} + \GGG_{j,\vb,-1} -
\GGG_{j,\vb_{1},0} - \GGG_{j,\vb_{1},-1} \right) , \cr} $$
which does not vanish since $\vb_{\ell_{T}^{1}} \neq \vb_{\ell_{T}}$.
In that case we say that the 1-resonance is {\it regularised}.

Then Lemma \secc(L8) holds true, with $L^{(k)}_{n,j,h)}$
substituted with $T^{(k)}_{j,h}$, provided that in the
definition of renormalised trees (cf. Definition \secc(D11))
we add the condition that all 1-resonances $T$
with $\vb_{\ell_{T}^{1}} \neq \vb_{\ell_{T}}$ and all the nodes with
$s_v=1$ and $i_v=1$  are regularised.

Also Lemma \secc(L9) is still true, as the property for the
matrix to be symmetric depends only the non-linearity.

\*\*
\0{\bf 8.7. Bryuno Lemma in $\TTTT_\RRRR^{(\kkkk)}$}\*

\0The set $\SS(\th,\g)$ is defined by \equ(4.1),
provided we substitute $(n,j)$ with $j$ and $\g$ with
$\g/32$. \equ(4.2) is replaced by: 

$$ \left\{
\eqalign{
& |\d_{n(a),m(a)}| \; \le \; 2^{-2}\g ,
\hskip1.2truecm \bb_{\ell}(a)=1 , \cr
2^{-3}\g \; \le \; & |\d_{n(a),m(a)}| \; \le 
\; 2^{-1}\g ,
\hskip1.2truecm \bb_{\ell}(a)=0 , \cr
2^{-2} \g \; \le \; & |\d_{n(a),m(a)}| ,
\hskip2.6truecm \bb_{\ell}(a)= -1 , \cr} \right.
\Eq(8.45) $$
for all $j_{\ell}\ge 1$ and $a=1,\ldots,d_{j_{\ell}}$.

\0 For the definition of the set $\DD(\th,\g)$ we require
only the condition \equ(4.3), which becomes
$$ \left| x_{j_{\ell}} \right| \ge {\g \over p_{j_{\ell}}^{\t}} .
\Eq(8.46) $$

We define $N_{h}(\th)$ as the set of lines $\ell$ with $i_{\ell}=1$
and scale $h_{\ell}\ge h$, which do not enter any resonance.
Then, with this new definition of $N_{h}(\th)$,
Lemma \secc(L10) remains the same. The proof follows the same
lines as in Section 4.1, with the following minor changes.

In order to have a line on scale $h$ we need that $B k^{1+4\a} \ge
C p_{j_\ell}\ge C 2^{(h-1)/\t}$ for some constant $C$.
We proceed as in the proof of Lemma \secc(L10), up to \equ(4.8),
where again $n_{\ell_{i}}$ should be substituted with $p_{j_{\ell_{i}}}$
with $i=0,1$.
\vskip.2truecm
\01. If $j_{\ell_{1}}=j_{\ell_0}$ then, since $\ell_1$ by hypothesis
does not enter a (regularised) resonance, there exists a line $\ell'$
with $i_{\ell'} \in\{0,-1\}$, not along the path
$\PPP(\bar\ell,\ell_{0})$, such that $j_{\ell'}=j_{\ell_{0}}$.
By the Remark after Definition \secc(D13), we know
that $|n_{\ell'}|\ge|n_{\ell_{0}}/2| > 2^{(h-2)/\t}$.
In this case one has $(k(\th)-k(\th_{1})^{1+4\a}>
B^{-1} |n_{\ell'}| \ge E_{h}$.
\vskip.1truecm
\02. If $j_{\ell_{1}}\neq j_{\ell_0}$ then we call $\bar\ell
\in\PPP(\ell_{0},\ell_{1})$ the line with $i\neq -1$ which is the
closest to $\ell_0$.
\vskip.1truecm
\02.1. If $p_{j_{\bar\ell}}\le p_{j_{\ell_0}}/2$ then
$(k(\th)-k(\th_{1}))^{1+4\a}\ge C p_{j_{\ell_0}}$.
\vskip.1truecm
\02.2. If $p_{j_{\bar\ell}}> p_{j_{\ell_0}}/2$ then one
reasons as in case 2.2. of Lemma \secc(L10),
with the following differences.
\vskip.1truecm
\02.2.1. If $j_{\ell_{0}} \neq j_{\bar\ell}$, then
$|(n_{\bar\ell}, m_{\bar\ell})-(n_{\ell_{0}},m_{\ell_{0}})|
\ge {\rm const.} p_{j_{\ell_{0}}}^{\b}$. For all the lines $\ell$
along the path $\PPP(\bar\ell,\ell_{0})$
one has $i_{\ell}=-1$, hence either $n_{\ell}=n_{\ell}'$ and
$m_{\ell}=m_{\ell}'$ (if $\ell$ is a $p$-line) or
$|m_{\ell}-m_{\ell}'| \le M_{1}$ (if $\ell$ is a $q$-line), 
so that $|(n_{\bar\ell},m_{\bar\ell})-(n_{\ell_{0}},m_{\ell_{0}})|
\le 2 B (k(\th)-k(\th_{1}))^{1+4\a}$,
with the same meaning for the symbols as in Section 4.1,
and the assertion follows once more by using \equ(4.8).
\vskip.1truecm
\02.2.2. If $j_{\ell_{0}} = j_{\bar\ell}$ then
there are two further sub-cases.
\vskip.1truecm
\02.2.2.1. If $\bar\ell$ does not enter any resonance,
we proceed as in item 1.
\vskip.1truecm
\02.2.2.2. If $\bar\ell$ enters a resonance, then we continue up
to the next line $\tilde \ell$ on the same path with
$i\neq -1$. If $j_{\tilde \ell}\neq j_{\ell_{0}}$ the proof is
concluded as in 2.2.1. since $2B k^{1+4\a}\ge |(n_{\bar\ell},
m_{\bar\ell})-(n_{\tilde\ell},m_{\tilde\ell})| \ge C_{1} p_{j}^\b$.
Likewise -- using item 2.2.2.1-- the proof is concluded if
the line $\tilde\ell$  does not enter a resonance.
If $\tilde\ell$ enters a resonance with $j_{\tilde \ell}= j_{\ell_{0}}$,
we proceed until we reach a line with $i\neq -1$
which either has $j\neq j_{\ell_{0}}$ or does not enter a
resonance: this is surely possible, because by definition $\ell_{1}$
does not enter a resonance and $j_{\ell_{1}} \neq j_{\ell_{0}}$.
This completes the proof of the lemma.
\vskip.2truecm

Lemma \secc(L11) holds with $|n|,|m| \le B k^{1+4\a}$ and $q=1$,
and with $p_{j}^{-3s/4}$ in all the lines of \equ(4.9).
The proof is the same (recall that we can set $s_{2}=0$);
we only need to substitute $p_{j}^\a$ (which bounded the dimension of
the non-diagonal block) with $d_{j}$.
In (iii) the labels $(n',j')$ should be substituted by $j'$.

\*\*
\0{\bf 8.8. Bryuno Lemma in $\RR_\RRRR^{(\kkkk)}$}\*

\0The definitions of $\widetilde{\SS}(\th,\g)$
and $\widetilde{\DD}(\th,\g)$ are changed exactly
as $\SS(\th,\g)$ and $\DD(\th,\g)$, respectively,
in the previous Section 8.7.

\*
\0\sub(D16)
{\it We divide $\R_{R,h,j}$ into two sets $\R_{R,h,j}^{1}$
and $\R_{R,h,j}^{2}$: $\R_{R,h,j}^{1}$ contains all the trees such that
either $\PPP(\ell_0,\ell_e)=\emptyset$ or at least one line $\ell \in
\PPP(\ell_0,\ell_e)$ has $j_\ell \neq j$, and
$\R_{R,h,j}^{2}=\R_{R,h,j} \setminus \R_{R,h,j}^{1}$.
This naturally yields a decomposition
$\R_{R,h,j}^{(k)}=\R_{R,h,j}^{(k,1)}\cup\R_{R,h,j}^{(k,2)}$
for all $k\in\NNN$.}
\*

The two properties (i) and (ii) of Lemma \secc(L12) should
be restated as follows.
\vskip.2truecm
\0{\it (i) There exists a positive constant $B_{2}$ such that
if $k \le B_{2} p_{j}^{\b/(1+4\a)}$ then
$\R_{R,j,h}^{(k,1)}$ contains only trees with
$\PPP(\ell_{0},\ell_{e})=\emptyset$;
\vskip.1truecm
\0(ii) for all $\th\in \R_{R,h,j}^{(k,1)}(a,b)$ we have
$|(n(a),m(a))-(n(b),m(b))|^{\r}\le k$,
with $\r$ a constant depending on $D$.}
\vskip.2truecm

The proof of (i) can be obtained by reasoning as in the cases
2.1. and 2.2.1. of Section 8.7, while that of (ii)
proceeds as in the proof of Lemma \secc(L12) (ii).

For the trees in $\R_{R,h,j}^{(k,2)}$ all the lines $\ell$
along the path $\PPP(\ell_{0},\ell_{e})$ have
$j_{\ell}=j$, and we can bound the product of
the corresponding propagators as
$$ \Big( \prod_{\ell\in\PPP(\ell_{0},\ell_{e})}
4 C\g^{-1} p_{j_{\ell}}^{-s} \Big)
\exp \Big( - \s \sum_{\ell \in \PPP(\ell_{0},\ell_{e})}
| (n_{\ell},m_{\ell})-(n'_{\ell},m'_{\ell})|^{\r} \Big)
\le C^{k} {\rm e}^{ - \s|(n_{\ell_{0}},m_{\ell_{0}})-
(n_{\ell_{e}},m_{\ell_{e}})|^{\r}} ,
\Eq(8.47) $$
where the factor is due to regularisation of the propagators
with $i_{\ell}=1$ (see \equ(8.38)), and we have used \equ(8.32)
to bound $|G_{j,\vb,i,-1}|_{\s}$ for $i=0,-1$.
Hence also $|\Val(\th)|_{\s}$ is bounded by $C^{k}$.

Lemma \secc(L13) and properties (i) and (ii) of Lemma \secc(L15)
are modified exactly as the corresponding \secc(L10) and \secc(L11).
In \equ(4.17) (ii) $|n|$ should be substituted by $|p_{j}|^{1+4\a}$.
Finally \equ(4.17) (iii) should be replaced with
$$ \sum_{j'\in \NNN}\sum_{a',b'=1}^{d_{j'}}
|\partial_{M_{j'}(a',b')}\Val(\th)| \le
D^{k} 2^{-h} \Big( \prod_{h'=-1}^{h} 2^{2 h' N_{h'}(\th)} \Big)
\prod_{\ell} p_{j_{\ell}}^{3s/4} ,
\Eq(8.48) $$
which can be proved as follows.
\vskip.2truecm
\01. Let us first consider $\R_{R,j,h}^{1}$. We have no difficulty
in bounding the sums and derivatives applied on lines
$\ell\notin\PPP(\ell_{0},\ell_{e})$. By the analog of Lemma \secc(L12)
discussed above, if $B_{2} k\le p_{j}^{\b/(1+4\a)}$
then $\PPP(\ell_{0},\ell_{e})=\emptyset$ and we have no problem.
Otherwise we have at most $(2p_{j} +k)^{1+4\a}$
possible values of $(n,m)$ and $(n',m')$ which can be associated
with a line $\ell$ along the path $\PPP(\ell_{0},\ell_{e})$
and by our assumption one has $(2p_{j} +k)^{1+4\a}\le C^k$
for some constant $C$.
\vskip.1truecm
\02. If all the lines $\ell\in\PPP(\ell_0,\ell_e)$ have $j_\ell=j$
then the sums with $a'\neq b'$ contain at most $p_{j}^{2\a}$
terms, whereas the sums with $a'=b'$ contain at most $k$ terms,
since there are at most $k$ lines on $\PPP(\ell_0,\ell_e)$.
\vskip.2truecm

The rest of Section 4 is unchanged.
In Section 5.1 we remove the second Melnikov condition
(the $**$ and $***$ products) in \equ(5.3) and \equ(5.5). 

\*\*
\0{\bf 8.9. Measure estimates}\*

\0By definition we have to evaluate the measure of the set
$$ \left\{ \e : \left\Vert \hat\chi_{1,j}
(\DDDD_j+p_{j}^{-s} \widehat M_{j})^{-1}
\hat\chi_{1,j} \right\Vert^{-1} \ge {2\g\over p_{j}^\t}
\quad \forall j\in\NNN \right\} .
\Eq(8.49) $$
By Lemma \secc(L3) (iii) one has
$$ x_{j} \ge \min_{i=1,\ldots,d_{j}}
\left| \l^{(i)} (\DDDD_j+p_{j}^{-s} \widehat M_{j}) \right| ,
\Eq(8.50) $$
since the matrices are symmetric and the minimum is attained
for some $i$ such that $\bar\chi_1(\d_{n(i),m(i)}) \neq 0$.

The set \equ(8.49) contains the set
$$ \EE = \left \{ \e\in (0,\e_0) : \left| \l^{(i)}
(\DDDD_j+p_{j}^{-s} \widehat M_{j}) \right| \ge
{2\g\over p_{j}^\t} \quad \forall i=1,\ldots,d_{j}, \quad
\forall j\in\NNN ,\right\} .
\Eq(8.51) $$
We estimate the measure of the subset of $(0,\e_0)$ complementary
to $\EE$, i.e. the set defined as union of the sets
$$ \III_{j,i}:= \left\{ \e\in (0,\e_0) : \left| \l^{(i)} (\DDDD_j+
p_{j}^{-s} \widehat M_{j} \right| \le {2\g\over p_{j}^\t} \right\}
\Eq(8.52) $$
for $j\in \NNN$ and $i=1,\ldots,d_{j}$.

First we notice that if $|p_j| \le C/\e_{0}$, for an appropriately
small $C$, then
$$ \left| \l^{(i)}(\DDDD_{j} - p_j I +  p_{j}^{-s} \widehat M_j)
\right| \le \e_{0}^{2} p_j^{3} \le {p_j \over 2} ,
\Eq(8.53) $$
which implies that
$$ \l^{(i)}(\DDDD_j+p_{j}^{-s} \widehat M_{j}) \ge {p_j \over 2} , $$ 
so that we have to discard the sets $\III_{j,i}$ only for $p_j \ge 
C/\e_{0}$.

Let us now recall that for a symmetric matrix $M(x)$
depending analytically  on a parameter $x$,
the derivatives of the eigenvalues are:
$\partial_x \l^{(i)}(x)= \langle v_i,\partial_x M(x) v_i\rangle$,
where $v_i$ are the corresponding eigenvectors \cita{Ka}.

Since $\DDDD_j$ depends linearly -- and therefore analytically --
on $\e$ we consider $\l_{i}(x,\e) := \l^{(i)}(\DDDD_j(x)+p_{j}^{-s}
\widehat M_{j})$ with $x,\e$ independent parameters.

Clearly $|\partial_x \l_{i}(x,\e)| \ge p_j$, and by Lidskii's Lemma
$$ \left| \partial_\e \l_{i}(x,\e) \right| \le p_j^{-s}
\sum_{i=1}^{d_{j}} \left| \l^{(i)}(\partial_\e \widehat M_j) \right| .$$
Now $\widehat M_j$ is a $d_j\times d_j$ matrix which for each fixed
$\bar\e$ has only a nonzero block of size $p_{j}^{\a}$; the properties
of the functions $\bar\chi_{j,1}$ imply that also $\partial_\e
\widehat M_j$ has only a nonzero block of size $p_{j}^{\a}$.
So one has
$$ \left| \partial_\e \l_{i}(x,\e) \right| \le
C( 1 + \e_0 p_j^{1-s+5\a} ) , $$
for some constant $C$.

Then the measure of each $\III_{j,i}$ can be bounded from above by
$$ {4\g\over p_{j}^\t } \sup_{\e\in (0,\e_{0})}
\left| \left( {\der \over \der\e}
\l^{(i)}\left(\DDD_j(\e)+p_j^{-s} \hat M_j(\e) \right)
\right)^{-1} \right| \le {8\g\over p_{j}^{\t+1} }.
\Eq(8.54) $$

Therefore we have
$$ \sum_{j\in\nnn} \sum_{i=1}^{d_{j}}
\hbox{meas} \left( \III_{j,i} \right) \le {\rm const.}
\sum_{p \ge C/\e_{0}} \g p^{D+\a}
\left( {1 \over p^{\t+1} } \right)
\le {\rm const.} \left( \e_{0}^{(\t-D-\a)} \right) ,
\Eq(8.55) $$
provided $\t>D+1+\a$, so that the measure of the complementary of $\EE$
is small in $(0,\e_0)$ if $\t>D+1+\a$.

\*\*
\appendix(A1,Preliminary measure estimate)

\0We estimate the measure of the complement of $\EE_{0}(\g)$, defined
in \equ(2.2), with respect to the set $(0,\e_{0})$,
under the condition $\m\in \MM$.
For all $n,p\in \NNN$ we consider the set
$$ \III_{n,p} = \left\{ \e\in (0,\e_{0}) :
|\o n-p|\le {\g\over n^{\t_{1}}} \right\} .
\Eqa(A1.1) $$
The measure of such a set is bounded proportionally
to $|n|^{-(\t_{1}+1)}$. Moreover one has
$$ \sum_{n,p=1}^{\io} \hbox{meas}(\III_{n,p}) \le
\hbox{const.} \sum_{n=1}^{\io} |n|^{-(\t_{1}+1)} +
\hbox{const.} \, \e_{0} \sum_{n=1}^{\io} |n|^{-\t_{1}} ,
\Eqa(A1.2) $$
because the number of values that $p$ can assume is at most
$1+\e_{0}n$ (simply note that $|\o n-p| \ge 1/2$ if $p$ is not
the integer closest to $\o n$ and $|\o-D-\mu|\le \e_{0}$).

Finally we note that, by \equ(2.1), for
$n<(\g_{0}/2\e_{0})^{1/(\t_{0}+1)}$ one has
$$ \left| \o n-p\right| \ge \left| (D+\m)n-p\right| -
\e_{0}|n|\ge \g|n|^{-\t_{0}} ,
\Eqa(A1.3) $$
provided $\g\le\g_{0}/2$. Hence the sum in \equ(A1.2) can be
restricted to $n \ge (\g_{0}/2\e_{0})^{1/(\t_{0}+1)}$, so that
$$ \sum_{n,p} \hbox{meas}(\III_{n,p}) \le
\hbox{const.}\, \e_{0}^{\t_{1}/(\t_{0}+1)} +
\hbox{const.}\, \e_{0}^{1+ (\t_{1}-1)/(\t_{0}+1)} ,
\Eqa(A1.4) $$
which is infinitesimal in $\e_{0}$ provided $\t_{1}>\t_{0}+1$.

\*\*
\appendix(A2,Proof of the separation Lemma {\secc(L2)})

\0Let $D\in \NNN$ be fixed, $D \ge 2$. For all $D>d\geq 1$ and for
all $r>1$ let $S^{d}(r)$ be a $d$-sphere of radius $r$ and
$S^{d}_{0}(r)$ the sphere $S^{d}(r)$ centred at the origin.
Set $\d(\e,d) := 2\e/d(d+2)!$, and let us denote with $|A|$ the
number of elements of the finite set $A$.

\*
\0\lma(LA1) {\it 
For all $\e\ll 1$ one can define sets of integer points
$\L_{\a}=\L_\a(\e,r,D,d)$, with $\a=1,\ldots,N=N(\e,r,D,d)$, such that
$$ |\L_{\a}| \le C(D,d) \, \max\{r^\e,d+2\} , \qquad
S^{d}_{0}(r)\cap \ZZZ^{D} = \bigcup_{\a=1}^{N} \L_{\a} , \qquad
{\rm dist}(\L_{\a},\L_{\b}) \ge C'(D,d)\,r^{\d(\e,d)} ,
\Eqa(A2.1) $$
where $C(D,d)$ and $C'(D,d)$
are suitable $(\e,r)$-independent constants.}
\*

The proof of this lemma follows easily from the following result.

\*
\0\lma(LA2) {\it There exist constants $C$ and $C'$ such that
the following holds. Let $n_{1},\dots,n_{k} \in S^{d}(r)\cap \ZZZ^{D}$.
If for all $i=1,\ldots,k-1$ one has $|n_{i}-n_{i+1}|<
C r^{\d(\e,d)}$ then $k<C'\max\{r^\e,d+2\}$.}
\*

\0{\bf Proof.} Let us first recall some trivial facts:
\vskip.2truecm
\01. $| S^{d}(r)\cap \ZZZ^{D}|\leq \bar C(D,d)\,r^{d}$,
for some constant $\bar C(D,d)$;
\vskip.2truecm
\02. given $p$ linearly independent vectors $v_{1},\dots,v_{p} \in
\ZZZ^{D}$ the volume of the $p$-dimensional simplex they
identify is given by
$$ {1\over p!} |\det NN^T|^{1\over 2} , \qquad N = \left( \matrix {
v_{11} & \ldots & v_{1D} \cr
\ldots &\ldots &  \ldots \cr
v_{p1} & \ldots & v_{pD} \cr} \right) ,
\Eqa(A2.2) $$
and, since $N$ has integer coefficients, the volume of the simplex
is bounded from below as $1/p!$.
\vskip.2truecm

Let us fix also some notations. Given $p$ linearly independent
vectors connecting points in $S^{d}(r)$, consider the $p$-dimensional
simplex generated by these vectors. Suppose that the angles between the
vectors are small enough: the volume of the simplex is bounded
from above by the volume of the spherical cap in which the vectors
are contained. If $\G$ is the radius of the base of the cap,
then the volume of the spherical cap is of order $\G^{d+2}/r$;
see Figure 8.

\midinsert
\*
\insertplotttt{170pt}{120pt}{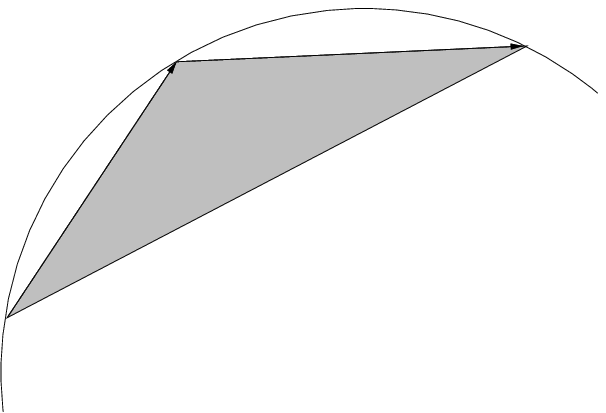}
\ins{160pt}{090pt}{$w_{1}$}
\ins{235pt}{125pt}{$w_{2}$}
\ins{215pt}{070pt}{$\G$}
\line{\vtop{\line{\hskip1.2truecm\vbox{\advance\hsize by -2.5 truecm
\noindent{\nota {\bf Figure 8.}
Simplex generated by two linearly independent vectors $w_{1}$
and $w_{2}$ which connects points on the sphere $S^{1}(r)$.
$\G$ is the basis of the spherical cap in which the two vectors
are contained. If the angle between the two vectors is small then
the volume (=area) of the cap is of order $\G^{3}/r$, with
$\G=O(|w_{1}|+|w_{2}|)$.\vfil
}} \hfill} }}
\*
\endinsert

Now we pass to the proof of the lemma. For $k\le d+1$ the assertion
is trivially satisfied, hence we can assume from now on $k\ge d+2$.
We proceed recursively. 
\vskip.2truecm
\0\underbar{Step 1.}
Consider three vectors in $\{n_{1},\ldots,n_{k}\}$ such that
the two difference vectors are linearly independent:
possibly reordering the vectors we can assume that they are
$n_{1}$, $n_{2}$ and $n_{3}$, and set $w_{1}=n_{2}-n_{1}$ and
$w_{2}=n_{3}-n_{2}$. The last two vectors connect points of some
$1$-sphere $S^{1}(r_{1})$, with $r_{1}\le r$.

Then there exists a constant $C_{1}$ such that
$\max\{|w_{1}|,|w_{2}|\} > C_{1} r_{1}^{1/3}$. The proof is by
{\it reductio ad absurdum}: consider the $2$-dimensional simplex
generated by $\{w_{1},w_{2}\}$; by the assumption on the distance
between the vectors and fact 2. we can say that there exist two
constants $D_{1}$ and $\bar D_{1}$ such that $r_{1} \le D_{1} \G^{3}
\le \bar D_{1} C_{1}^{3} r_{1}$, which is contradictory if $C_{1}$
is small enough.

One must have $r_{1} < A_{1} r^{\a_{1}\d}$, with $\d=\d(\e,d)$, for
$\a_{1}=3$ and a suitably large constant $A_{1}$: otherwise
$\max\{|w_{1}|,|w_{2}|\} > C_{1} r_{1}^{1/3} \ge A_{1} C_{1} r^{\d}$,
which is not possible if $A_{1} C_{1} > C$. By fact 1.
we have at most $\bar C(D,1)\,r_{1} \le \bar C(D,1)\,A_{1}r^{3\d}$
other integer vectors on $S^{1}(r_{1})\cap S^{d}(r)$.
\vskip.1truecm
\0\underbar{Step 2.}
Next, consider another vector (if any) in
$\{n_{1},\ldots,n_{k}\}$, say $n_{4}$,
such that $\{w_{1},w_{2},w_{3}\}$, with $w_{3}=n_{4}-n_{3}$,
are linearly independent and generate a simplex which is contained
in a $2$-sphere $S^{2}(r_{2})$ for some $r_{2}\le r$.
Of course the distance between $n_{4}$ and any vector
found in step 1. is bounded by $C r^{\d} + \bar B_{1} r_{1} \le
B_{1} r^{3\d}$, for suitable constants $\bar B_{1}$ and $B_{1}$.

Therefore we can prove, again by relying on fact 2., that
$\max\{|w_{1}|,|w_{2}|,|w_{3}|\} > C_{2} r_{2}^{1/4}$
for a suitable small constant $C_{2}$ (otherwise one would find
$r_{2} \le D_{2} \G^{4} \le \bar D_{2} C_{2}^{4} r_{2}$, hence a
contradiction for $C_{2}$ small enough).

Moreover, one must have $r_{2} < A_{2} r^{\a_{2}\d}$,
for $\a_{2}=4\a_{1}$ and a suitably large constant $A_{2}$;
otherwise $\max\{|w_{1}|,|w_{2}|,|w_{3}|\} > C_{2} r_{2}^{1/4}
\ge C_{2} A_{2} r^{\a_{2}\d/4} > B_{1} r^{\a_{1}\d}$,
which is not possible if $\a_{2}=4\a_{1}$ and $C_{2}A_{2}>B_{1}$.
By fact 1. we have at most $\bar C(D,2) \, r_{2}^{2} \le \bar C(D,2) \,
A_{2} r^{2\a_{2}\d}$ other integer vectors on $S^{2}(r_{2})\cap S^{d}(r)$.
\vskip.1truecm
\0\underbar{Step $j$.}
The proof is performed by induction. Assume that,
up to step $j-1$., we have found at most a $\bar C(D,1)\,
A_{1}r^{\a_{1}\d} + \bar C(D,2)\,A_{2}r^{2\a_{2}\d} + \ldots +
\bar C(D,j-1)\,A_{j-1}r^{j\a_{j-1}\d}$ vectors, with $\a_{i}=(i+2)!/2$
and suitably large constants $A_{i}$, such that the distance
between any two among these vectors is less than
$B_{j-1}r^{a_{j-1}\d}$ for a suitable constant $B_{j-1}$.

Moreover there are at least $j$ vectors,
which are linearly independent: we can assume are
$\{n_{1},\ldots,n_{j}\}$ and set $w_{i}=n_{i+1}-n_{i}$ for
$i=1,\ldots,j-1$. Suppose that there is at least another vector
$n_{j+1}$ on $S^{d}(r)$ which does not depend linearly on
$\{n_{1},\ldots,n_{j}\}$, and set $w_{j}=n_{j+1}-n_{j}$
(if there is no such vector the proof becomes easier).
Call $S^{j}(r_{j})$ the $j$-sphere which contain the $j$-simplex
generated by $\{w_{1},\ldots,w_{j}\}$. Once more fact 2. implies
that there is a constant $C_{j}$, small enough, such that
$\max\{|w_{1}|,\ldots,|w_{j}|\} > C_{j} r_{j}^{1/j}$.

One must have $r_{j} < A_{j} r^{\a_{j}\d}$ for $\a_{j}=(j+2)\a_{j-1}$
and $A_{j}$ suitably large: if this were not true then one would have
$C_{j}r_{j}^{1/j} \ge C_{j} A_{j} r^{\a_{j}\d/j} > B_{j-1}
r^{\a_{j-1}\d}$, hence a contradiction if $\a_{j}=j\a_{j-1}$
and $C_{j}A_{j} > B_{j-1}$. Hence the number of other vectors
that we have to add at this step is at most $\bar C(D,j)\,r_{j}^{j}
\le \bar C(D,j)\,A_{j} r^{j\a_{j}\d}$, and the distance between
all the points is bounded by $\bar B_{j}r_{j}^{j} \le B_{j}
r^{\a_{j}\d}$, for suitable constants $\bar B_{j}$ and $B_{j}$.
Hence the inductive hypothesis is satisfied.
\vskip.2truecm
The inductive estimate for $j=d$ yields the result,
provided one sets $\e=d(d+2)!\d/2$ and one chooses
$C$ and $1/C'$ small enough. This completes the proof. \qed

\*

\0{\bf Remarks.}
(1) A careful look at the proof of Lemma \secc(LA2) shows
that $C'=C'(D)$ is the maximum of $d\,C(D,d)\,A_{d}$ for $1\le d < D$,
hence $C'=(D-1)\,C(D,D-1)\,A_{D-1}$, whereas
$C=C(D)$ is obtained as the minimum between the constant $C$ and
the constants $B_{d}$ for $1\le d<D-2$, hence $C'=C$.
This shows that in Lemma \secc(LA1) one can
choose $C(D,d)$ and $C'(D,d)$ as functions of the only $D$.\\
(2) In the proof of Lemma \secc(LA2) the construction in
step 1. shows that if one takes three vectors $n_{1}$, $n_{2}$
and $n_{3}$ on a $1$-sphere $S^{1}(r_{1})$ then (with the notations
used in the proof of the lemma) one has
$\max\{w_{1},w_{2}\}>C_{1}r_{1}^{1/3}$. Therefore for $d=1$ these
sets $\L_{j}$ can be chosen in such a way that each set contains
at most two elements, and the distance between two distinct sets
on the same sphere $S^{1}(r)$ is larger than a universal
constant times $r^{1/3}$.

\*

Lemma \secc(LA1) implies that it is possible to decompose
the set $\ZZZ^{D} \cup S_{0}^{D}(r)$ as the union of sets
$\Delta$ such that ${\rm diam}(\Delta)< {\rm const.} r^{\d+\e}$
(cf. \cita{Bo2}, p. 399), and $|\Delta| < {\rm const.}r^{D(\d+\e)}$.
Hence, if we take $\a$ small enough and we set $\b=\d$ and
$\a=D(\d+\e)$, by using that $\e/\d=(d+2)!d/2$, Lemma \secc(L1) follows.

\*\*
\appendix(A3,Constructive scheme for Lemma {\secc(L17)})

\0Here we prove that the sets $\MMM_{+}$ verifying the
conditions (a) and (b) in the proof of Lemma \secc(L17) are non-empty.
The proof consists in providing explicitly a construction.\\
\01. Fix a list of parameters $\a_{2},\dots,\a_N \in \RRR$ such that
$\a_{i} < \a_{i-1}$ for $i=2,\ldots,N$, with $\a_{1}=1$, and
$$ 2^D\sum_{i=2}^N\a_i^{2+2s}  \le 3^D+2^D(N-2) .
\Eqa(A3.1) $$
\02. Given $r\in \RRR^+$ and for $i=2,\dots,N$ consider the regions
$\RR_i(r):=\{ x\in \ZZZ^D_+: \a_{i-1} r\le |x|\le \a_{i} r\}$ with
$r$ so big that it is not possible to cover any of the $\RR_i(r)$
with $3 N^22^{2D}$ planes and spheres.\\
\03. Choose an integer vector $m_{1}\in \ZZZ_+^D$  such that
$|m_{1}|^2 = r^{2}$ is divided by $D$, and construct the ``orbit''
$\OO(m_{1}):= \{ m\in \ZZZ^D: |m_{i}|=|(m_{1})_i|$.\\
\04. For each pair $m,m'\in \OO(m_{1})$ consider the two planes
orthogonal to $m-m'$ and passing respectively through $m$ and $m'$,
and the sphere which has $m-m'$ as diameter (there are at most
$3 \cdot 2^{D-1}(2^D-1)$ planes and spheres).\\
\05. Choose the second integer vector $m_{2}\in \RR_{2}(r)$ such that
$|m_{2}|^2$ divides $D$ and the orbit $\OO(m_{2})$ does not lie on
any of the planes and spheres defined at step 4.\\
\06. For each pair $m,m'\in \OO(m_{1}) \cup \OO(m_{2})$
proceed as in step 4. We have at most further $3 \cdot 2^{D}(2^{D+1}-1)$
planes and spheres.\\
\07. Then we proceed iteratively. When we arrive to $m_{N}$ we have
to remove at most $3 N 2^{D-1}(N2^{D}-1)$ planes and spheres.

\*\*
\appendix(A4,Blocks of the matrix $\titolo\JJJJ$)

\0Write $\MMM=\{m_{1},\ldots,m_{M}\}$, with $M=2^{D}N$,
and set $\VV=\{ v = (m,m') : m,m' \in\MMM , \,
m \neq m'\}$: clearly  $L:=|\VV|=M(M-1)$.
We call {\it alphabet} the set $\VV$ and {\it letters} the
elements (vectors) of $\VV$. We call {\it word} of length $\ell \ge 1$
any string $v_{1}v_{2}\ldots v_{\ell}$, with $v_{k}\in \VV$
for $k=1,\ldots,\ell$. Let us denote with $\AAA$ the set of
all words with letters in the alphabet $\VV$ plus the empty set
(which can be seen as a word of length $0$).

For $v\in \VV$ with  $v=(m_{i},m_{j})$
we write $v(1)=m_i$ and $v(2)=m_j$.
Given two words $a=v_{1}\ldots v_{n}$ and $b=v_{1}'\ldots v_{n'}'$
we can construct a new word $ab=v_{1}\ldots v_{n}v_{1}'\ldots v_{n'}'$
of length $n+n'$. Finally we can introduce a map $a \to w(a)$, which
associates with any letter $v\in\VV$ the vector $v(1)-v(2)$,
to any word $a=v_{1}\ldots v_{n}$ the vector $w(a)=
w(v_{1})+\ldots +w(v_{n})$ and finally $w(\emptyset)=0$.
We say that $a$ is a {\it loop} if $w(a)=0$.

\*

\0{\bf Remarks.} (1) Given a set $\MMM$ let $\VV$ be the corresponding
alphabet. If $|\MMM|=M$ then $|\VV|=L(M)=M(M-1)$. If we add an element
$m_{N+1}$ to $\MMM$ so to obtain a new set $\MMM'=\MMM\cup\{m_{N+1}\}$,
then the corresponding alphabet $\VV'$ contains all the
letters of $\VV$ plus other $2M$ letters. We can imagine that this
alphabet is obtained through $2M$ steps, by adding one by one the $2M$
new letters. In this way, we can imagine that the length $L$ of
the alphabet can be increased just by 1.\\
(2) By construction $w(v_{1}v_{2})=w(v_{2}v_{1})$.
In particular $w(a)$ depends only on the letters of $a$ (each with
its own multiplicity), but not on the order they appear within $a$.

\*

Define a matrix $J$, such that
\vskip.1truecm\noindent
(i) $J_{jk}=J(q_{j},q_{k})$, with $q_{j},q_{k}\in\ZZZ^{D}$,
\vskip.1truecm\noindent
(ii) $J(q,q')\neq 0$ if there exist $m_{1},m_{2} \in \MMM$ such that
$q-m_{1}=q'-m_{2}$ and $\la m'-m_{2},m_{1}-m_{2} \ra=0$, and
$J(q,q')=0$ otherwise.
\vskip.1truecm

A sequence $C=\{q_{0},q_{1},\ldots,q_{n}\}$ will be called a {\it chain}
if $J(q_{k-1},q_{k}) \neq 0$ for $k=1,\ldots,n$. We call $n=|C|$
the length of the chain $C$. A chain can be seen as a pair of a vector
and a word, that is $C=(q_{0};a)$, where $q_{0}\in\ZZZ^{D}$ and
$a=v_{1}\ldots v_{n}$, with $w(v_{k})=q_{k}-q_{k-1}$. Note that,
by definition of the matrix $J$, given a chain $C$ as above, one has
$$ q_{k} = q_{k-1} + w(v_{k}) , \qquad \la q_{k} - v_{k}(2),
w(v_{k}) \ra = 0 ,
\Eqa(A4.1) $$
for all $k=1,\ldots,n$.

\*
\0\lma(LA3) {\it Given a chain $C=(q_{0};a)$, if the word $a$ contains
a string $v_{0}a_{0}v_{0}$, with $v_{0}\in\VV$ and $a_{0}\in\AAA$,
then $\la w(v_{0}a_{0}),w(v_{0}) \ra =0$.}
\*

\0{\bf Proof.} As the word $a$ of $C$ contains the string
$v_{0}a_{0}v_{0}$, by \equ(A4.1) there exists $j \ge 1$ such that
$$ \la q_{j} - v_{0}(2) , w(v_{0}) \ra = 0 , \qquad
\la q_{j} + v_{0} + w(a_{0}) - v_{0}(2) ,w( v_{0}) \ra = 0 , $$
so that $\la w(v_{0})+w(a_{0}) , w(v_{0}) \ra = 0$. \qed

\*
\0\lma(LA4) {\it Given a chain $C=(q_{0};a)$, if the word $a$ contains
a string $a_{0}b_{0}a_{0}$, with $a_{0},b_{0}\in\AAA$ and $a_{0}$
containing all the letters of the alphabet $\VV$,
then $a_{0}b_{0}$ is a loop.}
\*

\0{\bf Proof.} For any $v\in\VV$ we can write $a_{0}=a_{1}va_{2}$,
with $a_{1},a_{2}\in\AAA$ depending on $v$. Then
$a_{0}b_{0}a_{0}=a_{1}va_{2}b_{0}a_{1}va_{2}$. Consider the
string $va_{2}b_{0}a_{1}v$: by Lemma \secc(LA3) one has
$\la w(va_{2}b_{0}a_{1}) , w(v) \ra = 0$. On the other hand
(cf. Remark (2) after the definition of loop) one has
$w(va_{2}b_{0}a_{1}) = w(a_{1}va_{2}b_{0})=w(a_{0}b_{0})$, so that
$\la w(a_{0}b_{0}), w(v) \ra = 0 $. As $v$ is arbitrary we conclude that
$$ \la w(a_{0}b_{0}) , w(v) \ra = 0 \quad \forall v \in \VV \qquad
\Longrightarrow \qquad w(a_{0}b_{0}) = 0 , $$
i.e. $a_{0}b_{0}$ is a loop. \qed

\*
\0\lma(LA5) {\it There exists $K$ such that if a word has
length $k \ge K$ then the word contains a loop. The value
of $K$ depends only on the number of letters of the alphabet.}
\*

\0{\bf Proof.} The proof is by induction on the length $L$
of the alphabet $\VV$ (cf. Remark (1) after the definition of loop).

For $L=1$ the assertion is trivially satisfied.
Assume that for given $L$ there exists an integer $K(L)$ such that
any word of length $K(L)$ containing at most $L$ distinct
letters has a loop: we want to show that then if the alphabet
has $L+1$ letters there exists $K(L+1)$ such that any word of the
alphabet with length $K(L+1)$ has also a loop.

Let $N(L)$ be the number of words of length $K(L)$ written
with the letters of an alphabet $\VV$ with $|\VV|=L+1$.
Consider a word $a=a_{1}\ldots a_{N(L)+1}$,
where each $a_{k}$ has length $K(L)$. We want to show by contradiction
that $a$ contains a loop. If this is not the case, by the
inductive assumption for each $k$ either $a_{k}$ contains a  loop
or it must contain all the $L+1$ letters. As all words $a_{k}$
have length $K(L)$ and there are $N(L)+1$ of them,
at least two words, say $a_{i}$ and $a_{j}$ with $i<j$,
must be equal to each other. Therefore we can write
$ a = a_{1}\ldots a_{i-1} a_{i} b a_{i} a_{j+1}
\ldots a_{N(L)+1}$, where $b=a_{i+1}\ldots a_{j-1}$
if $j>i+1$ and $b=\emptyset$ if $j=i+1$.
Hence $a$ contains the string $a_{i}ba_{i}$, with $a_{i}$
containing all the letters. Hence by Lemma \secc(LA4)
one has $w(a_{i}b)=0$, i.e. $a_{i}b$ is a loop. \qed

\*

\0{\bf Remark.} Note that the proof of Lemma \secc(LA5) implies
$$ K(L+1) \le K(L) \left( N(L)+1 \right) \le
\prod_{\ell=1}^{L} \left( N(\ell)+1 \right) ,
\Eqa(A4.2) $$
which provides a bound on the maximal length of the chains in terms
of the length of the alphabet $\VV$.

\*

Lemma \secc(L18) follows immediately from the results above,
by noting that all the spheres with diameter a vector $v(1)-v(2)$
with $v\in \VV$
are inside a compact ball of $\ZZZ^{D}$.

\*\*
\appendix(A5,Invertibility of $\JJJJ$ for $\DDDDD$=2)

\0In the following we assume $D=2$ and $N>4$.
We first prove that (i) implies (ii). As seen in
Appendix \secc(A3) condition \equ(8.4) is implied by
$$ \a_i\le \left( {|m_i|\over |m_1|} \right)^{2+2s}
\le \a_{i+1}\,,\quad \forall i=2,\ldots,N-1,
\Eqa(A5.1) $$
where the $\a_i>1$ are fixed in Appendix \secc(A3).

For $|m_1|$ large enough, \equ(A5.1) contains a $2N$-dimensional ball
of arbitrarily large radius. By definition an algebraic variety is
the set of solutions of some polynomial equations and therefore cannot
contain all the positive integer points of a ball provided the radius is
large enough (depending on the degree of the polynomial).

To prove (i) let us start with some notations. We consider
$\ZZZ^{2N}$ as a lattice in $\CCC^{2N}$, we denote $x=\{x_1,\dots,x_N\}
\equiv \MMM_+\in \CCC^{2N}$,
where each $x_i$ is a point in $\CCC^2$; we denote the points
in $\MMM$ still as $m_i\in \CCC^2$, and for each point $x_i\in\MMM_+$
we have the orbit $\OO(x_i)\in \MMM$ i.e. the four points in $\MMM$
obtained by changing the signs of the components of $x_i$. 

\*
\asub(D17){\it (i) Given two points $m_i,m_j$ in $\MMM$ we consider:
the circle with diameter $m_i-m_j$ (curve of type 1)
the two lines orthogonal to $m_i-m_j$ and passing
respectively through $m_i$ (curve of type 2) and
through $m_j$ (curve of type 3). Note that the curve is
identified by the couple $(m_{i},m_{j})$ and by the type label.
We call $\CC$ the finite set of distinct curves obtained in this way
for all couples $m_i\neq m_j$ in $\MMM$.\\
(ii) Let $C$ be a curve in $\CC$ identified by the couple
$(m_{i},m_{j})$. We say that a point $m'$ is g-linked by
$(m_{i},m_{j})$ to $m\in C$ if one has either (1) $m'= -m+m_i+m_j$,
if $C$ is a curve of type 1, or (2) $m'= m + (m_j-m_i)$,
if $C$ is a curve of type 2, or (3) $m'= m - (m_j-m_i)$,
if $C$ is a curve of type 3. Notice that in case (1) also $m'$ is
on the circle, while in cases (2) and (3) $m'$ is on a curve
of type 3 and 2, respectively.
We say that two points $m,m'\in\ZZZ_{+}^{2}$ are
linked by $(m_{i},m_{j})$ if there are two points $\bar m\in\OO(m)$
and $\bar m'\in\OO(m')$ such that $\bar m,\bar m'$
are g-linked by $(m_{i},m_{j})$.\\
(iii) Given $\MMM_+$ we consider the set $H$ of points $y_j\notin \MMM$
which lie on the intersection of two curves in $C$,
counted with their multiplicity. Set $r := |H|$: we
denote the list of intersection points as $y=\{y_1,\dots,y_{r(N)}\}$.
Note that $r$ depend only on $N$.}
\*

We first prove that the points $x\in \CCC^{2N}$ which do not
satisfy Lemma \secc(L17) lie on an algebraic variety.
As seen in Appendix \secc(A3), Lemma \secc(L17) is verified by
requiring  that if either a curve of type 1 contains three points in
$\MMM$ or a curve of type 2 or 3 contains two points in $\MMM$,
then such points are on the same orbit. It is clear
(see Appendix \secc(A3)) that this condition can be achieved
by requiring that $x$ does not belong to some proper algebraic variety,
say $\WW_{a}$, in $\CCC^{2N}$.

Let us now consider the set of points $x\in \ZZZ_+^{2N}$ where
${\rm det}J_{1,1}$ is identically equal to zero
(as a function of $s$); since $J_{1,1}$ is a block diagonal matrix
we factorise the single blocks and treat them separately.
The matrix $J_{1,1}$ has some simple blocks
which we can describe explicitly. Recall that
$$ 16A= c_1 \sum_{i=1}^{N}|x_i|^2\,,\;\qquad
2 a_{x_i}^2= (1-c_1) |x_i|^2 -c_1\sum_{j=1,\ldots,N
\atop j\neq i} |x_j|^2 ,
\Eqa(A5.2) $$
where $c_{1}=8/(8N+1)$.

\vskip.2truecm
\01. For all $m\in \ZZZ^2_+$ such that $m$ does not belong to any curve
$C\in\CC$ one has $Y_{m,m'}= 0$ for all $m'$; by considering the limit
$s\to \io$ one can easily check that $J_{m,m}=|m|^{2+2s}/2-8A=0$
is never an identity in $s$ (independently of the choice of $\MMM_+$).
\vskip.1truecm
\02. For all linked couples  $m,m'\in \ZZZ^2_+$
such that each point belongs to one and only one
curve  one has either a diagonal block $|m|^{2+2s}/2-8A-4a_{x_i}^2$
for some $x_i\in \MMM_+$ if $m=m'$, or a $2\times 2$ matrix
$$ \left(\matrix{ |m|^{2+2s}/2 -8A & -2{a_{m_i}a_{m_j}}\cr
-2{a_{m_i}a_{m_j}}&|m'|^{2+2s}/2 -8A }\right) $$
if $m \neq m'$ and $(m_i,m_j)$ is the couple linking $m'$ to $m$.
In both cases a trivial check of the limit $s\to \pm\io$ will
ensure that the determinant is not identically null. 
\vskip.1truecm
\03. There is a block matrix containing all and only the elements of
$\MMM_{+}$. Such a matrix is easily obtained by differentiating
the left hand side of \equ(8.5):
$$ - 2 \,\left(\matrix{a_{m_1}\!\!\!\!\!\!\! & 0 &\dots  &0
\cr  0 &\ddots &\ddots &\vdots\cr \vdots & \ddots &\ddots & 0
\cr 0 &\dots &0 &a_{m_N}\!\!}\right)\left(\matrix{9 & 8 & \dots & 8
\cr  8 &\ddots &\ddots &\vdots\cr \vdots & \ddots &\ddots & 8
\cr 8 &\dots &8 &9}\right)\left(\matrix{a_{m_1}\!\!\!\!\!\!\! & 0
&\dots  &0\cr  0 &\ddots &\ddots &\vdots\cr \vdots & \ddots &
\ddots & 0 \cr 0 &\dots &0 &a_{m_N}\!\!}\right).$$
Since all the $a_{m_i}$ are non-zero we only need to prove
that the matrix in the middle is invertible, which is trivially
true since the determinant is an odd integer. 
\vskip.2truecm

We now have considered all those blocks in $J_{1,1}$ whose
invertibility can be easily checked directly. We are left with
the intersection points in $H':=H \cap \ZZZ_{+}^{2} \setminus \MMM_{+}$
and all those points $m'$ which are linked to some $y_{j} \in H'$.
We call $\tilde J$ the restriction of $J_{1,1}$ to such points;
the crucial property of $\tilde J$ is that it is a $K\times K$
matrix with $K$ bounded by above by some constant depending only on $N$.
 
We will impose that $\tilde J$ is invertible at $s=0$ by requiring
that $\MMM_+$ does not lie  on an appropriate algebraic variety
in $\CCC^{2N}$.

By definition the points in $H$ (and the points linked to them) 
are algebraic functions of $x\in \CCC^{2N}$.
 By construction $\tilde J_{m,m}-Y_{m,m}=|m|^{2+2s}/2-8A$
and moreover  $Y_{m,m'}$
contains a contribution  $-2 a_{m_i}a_{m_j}$  for each couple
$(m_i,m_j)$ linking $m'$ to $m$.
We want to prove that
for $s=0$ the equation ${\rm det} \, \tilde J =0$ (which is an equation for 
$x\in \CCC^{2N}$) defines a proper algebraic variety, say $\WW_f$, in $\CCC^{2N}$.

We consider the space $ \CCC^T := \CCC^{2N}\times
\CCC^N\times \CCC^{2r}$ and, with an abuse of notation, we denote
 the generic point in $\CCC^T$ by $(x,a,y)=$ $
( x_1,\ldots,x_N,a_{x_1},\ldots,a_{x_N},$ $y_1,\dots,y_{r})$ 
(therefore  we consider $(x,a,y)$ as independent variables). Note that
 ${\rm det} \, \tilde J =0$ is a polynomial equation in $\CCC^T$.
We call $\WW_{b}$ the algebraic variety  defined by requiring both
that the $a_{x_i}$ satisfy \equ(A5.2) and that each $y_j$ lies
on at least two curves of $C$ ($\WW_{b}$ is equivalent to a
finite number of copies of $\CCC^{2N}$). 

We now recall a standard theorem in algebraic geometry which states:
{\it Let $W$  be  an algebraic variety  in $\CCC^{n+m}$ and let $\Pi$
be the projection $\CCC^{n+m}\to \CCC^{n}$ then $\overline{\Pi(W)}$
is an algebraic variety} (clearly it may be the whole $\CCC^{n}$!)
We set $n=2N$ (the first $2N$ variables), $m=2r+N$ and  apply the stated 
theorem to $\overline{\Pi(\WW_{b}\cap \WW_{f})}$; we now only need to prove
 that the algebraic variety we have
obtained is proper; to do so it is convenient to treat separately
the invertibility conditions of each single block of $\tilde J$. 

The first step is to simplify as far as possible the structure of the
intersections and therefore of the matrix $\tilde J$. The simplest
possible block involving an intersection point $y_j$ is such that
\vskip.2truecm
\0(i) only two curves in $\CC$ pass through $y_j$;
\vskip.1truecm
\0(ii) the two points linked to $y_j$ (by the couples of points
in $\MMM$ identifying the curves) are not intersection points.
\vskip.2truecm

Such a configuration gives either a $3 \times 3$ matrix
or a $2\times 2$ matrix -- if one of the curves is either
an horizontal or vertical line or a circle centred at the origin.

\*
\asub(D18){\it We say that a curve $C\in\CC$ depends on
the two -- possibly equal -- variables $x_i,x_{j} \in \CCC^2$ if
$C$ is identified by the couple $(m_{i},m_{j})$, 
such that $m_i\in \OO(x_i) $ and $m_{j}\in \OO(x_{j})$.}
\*

The negation of (i) is that $y_j$ is on (at least) three curves
of $\CC$: such condition defines a proper algebraic
variety in $\CCC^{T}$, say $\WW_j$. We now consider
the projection of $\WW_b\cap \WW_j$ on $\CCC^{2N}$: 
its closure is an algebraic variety and either it is
proper or the triple intersection occurs for any choice of $x$
(which unfortunately can indeed happen due to the symmetries introduced
by the Dirichlet boundary conditions).

\0 Three curves in $\CC$ depend on at most six variables in $\CCC^2$.
If four or more of such variables are different then at least
one variable, say $x_{k}$, appears only once. By moving $x_{k}$
in $\CCC^2$ we can move arbitrarily one of the curves, while the
other two (which do not depend on  $x_{k}$) remain fixed.
This implies that the triple intersection cannot hold true
for all values of $x_{k}$ and thus $\overline{\Pi(\WW_{b}\cap \WW_j)}$
is a proper variety in $\CCC^{2N}$.
 
In the same way the negation of (ii) is that one point
linked to $y_j$ lies on (at least) two curves of $\CC$
(one curve is fixed by the fact that the point is linked to $y_j$);
again the intersection is determined by six points in $\MMM$
and the same reasoning holds. 

We call $\WW_{c}$ the variety in $\CCC^{T}$
defined by the union of all those $\WW_{j}$ such that
$\overline{\Pi(\WW_{b}\cap \WW_j)}$ is proper.

In $\WW_{b} \setminus \WW_{c}$ we can now classify the possible
blocks appearing in $\tilde J$ (notice that only 
intersection points which are integer valued have to be
taken into account when constructing the blocks in $\tilde J$).
\vskip.2truecm
\01. We have a list of at most $3\times 3$  blocks
corresponding to the intersection points of type (i)-(ii).
Such intersection points are identified by two curves which can
depend on at most four different variables 
$x_{i_k}$ with $k=1,\dots,4$.
\vskip.2truecm
\02. There are more complicated blocks corresponding to multiple
intersections (or intersection points linked to each other),
which occur for all $x\in \CCC^{2N}$ due to symmetry.
As we have proved above the curves defining such intersections
depend on at most three different variables $x_{i_k}$.

In any given block, call it $B_{h}$, the contribution from $Y$ involves
only terms of the form $-2 a_{m_i}a_{m_j}$ such that $m_i,m_j\in \cup_{k=1}^4\OO(x_{i_k})$.
Each $a_{m_j}$ depends on all the components of $x$;
in particular, $a_{m_j}^2$ can be written as a term depending
only on the $x_{i_k}$ plus the term $-{1\over2}
c_1 \sum_{j\neq i_1,\dots,i_4} x_{j}^2$.
Since by hypothesis $N>4$ and $k\leq 4$
the second sum is surely non-empty.

Finally one has the  diagonal contributions (from $J-Y$): $|y_{j}|^2- {1\over2}
c_1 \sum_{j\neq i_1,\dots,i_4}x_j^2 +z$, where $z$ is a polynomial
function in the $x_{i_k}$'s.

In the limit $\sum_{j\neq i_1,\dots,i_4}x_j^2 \to \io$
the terms depending on the $x_{i_k}$'s become irrelevant
and we are left with a matrix (of unknown size) whose entries,
apart from the common factor $ - {1\over2} 
c_1\sum_{j\neq i_1,\dots, i_4}x_j^2$, are integer numbers.
It is easily seen that these numbers are odd on the diagonal,
 while all the off-diagonal terms are even; indeed $Y$ contributes only 
even entries while $J-Y$ is diagonal and odd due to the term $8A$.
Thus the determinant (apart form the common factors) is odd and
hence the equation ${\rm det} B_{h}=0$ is not an identity on $\WW_{b}$.
If we call $\WW_{h}$ the variety in $\CCC^T$ defined by ${\rm det}
B_{h} =0$ then $\overline{\Pi(\WW_{b}\cap \WW_{h} )}$ is surely
proper. Finally we call $\WW_{d}$ the union of all the
$\WW_{h}$ and set $\WW_{f}= \WW_{d} \cup \WW_{c}\cup \WW_{c}$.

\*\*
\appendix(A6,Proof of the separation Lemma {\secc(L20)})

\0The following proof is adapted from \cita{Bo4}.
Given $\e>0$ define $\d=\d(\e,D)=\e/2^{D-1}D!(D+1)!$.
Then Lemma \secc(L20) follows from the results below.

\*
\0\lma(LA9) {\it Let $x\in\RRR^{d}$.
Assume that there exist $d$ vectors $\D_{1},\ldots,\D_{d}$,
which are linearly independent in $\ZZZ^{d}$, and such
that $|\D_{k}| \le A_{1}$ and $|x\cdot\D_{k}| \le A_{2}$ for all
$k=1,\ldots,d$. Then $|x| \le C(d) A_{1}^{d-1} A_{2}$
for some constant $C(d)$ depending only on $d$.}
\*

\0{\bf Proof.} Call $\b_{k}\in[0,\p/2]$ the angle between $\D_{k}$
and the direction of the vector $x$. Without any loss of generality
we can assume $\b_{k} \ge \b_{d}$ for all $k=1,\ldots,d-1$.
Set $\b_{d}'=\p/2-\b_{d}$. One has $\b_{d}'>0$ because
$\D_{1},\ldots,\D_{d}$ are linearly independent.

Consider the simplex generated by the vectors
$\D_{1},\ldots,\D_{d}$. By the fact 2. in the proof of Lemma
\secc(L20) one has, for some $d$-dependent constant $C(d)$,
$$ 1 \le C(d) \left| \D_{1} \right| \left| \D_{2} \right| \ldots
\left| \D_{d} \right| \left| \sin\a_{1,2} \right|
\left| \sin \a_{12,3} \right| \ldots
\left| \sin \a_{1\ldots (d-1),d} \right| ,
\Eqa(A6.1) $$
where $\a_{1\ldots (j-1),j}$, $j\ge 2$, is the angle between the
vector $\D_{j}$ and the plane generated by the vectors
$\D_{1},\ldots,\D_{j-1}$. Hence
$$ 1 \le C(d) A_{1}^{d-1} |\D_{d}| |\sin\a_{1\ldots (d-1),d}| .
\Eqa(A6.2) $$ 

Moreover one has
$$ \left| x \cdot \D_{d} \right| =
\left| x \right| \left| \D_{d} \right| \left| \cos \b_{d} \right| =
\left| x \right| \left| \D_{d} \right| \left| \sin \b_{d}' \right|
\ge \left| x \right| \left| \D_{d} \right|
\left| \sin \a_{1\ldots(d-1),d} \right| ,
\Eqa(A6.3) $$
so that, from \equ(A6.2) and \equ(A6.3), we obtain $|x|A_{1}^{-(d-1)}
\le C(d) A_{2}$, so that the assertion follows. \qed

\*
\0\lma(LA10) {\it There exist constants $C$ and $C'$ such that
the following holds. Let $n_{1},\ldots,n_{k} \in \ZZZ^{D}$ be a
sequence of distinct elements such that $|\Phi(n_{j})-\Phi(n_{j+1})|
\le C r^{\d}$. Then $k \le C' \max \{r^{\e},D+2\}$.}
\*

\0{\bf Proof.}
Since the vectors $n_{j}$ are on the lattice $\ZZZ^{D}$ there exist
a constant $C_{1}(D)$ and $j_{0} \le k/2$ such that
$|n_{j_{0}}| > C_{1}(D) k^{1/D}$. Set $\D_{j}=n_{j}-n_{j_{0}}$.
By assumption one has $|\Phi(n_{j})-\Phi(n_{j+1})| \le C r^{\d}$,
hence $|\Phi(n_{j})- \Phi(n_{j_{0}})| \le C (j-j_{0}) r^{\d}$ for all
$j_{0}+1 \le j \le k$. Then $|\Phi(n_{j})-\Phi(n_{j_{0}})|
\le A_{1}:=C J_{1} r^{\d}$ for all $j_{0}+1 \le j \le j_{0}+J_{1}$.
Fix $J_{1}=k^{1/\a(D)}$, with $\a(n)=2n(n+1)$.
By using that
$$ \Phi(n_{j})-\Phi(n_{j_{0}})= \left( \D_{j},2\D_{j}
\cdot n_{j_{0}}+ |\D_{j}|^{2} \right) ,
\Eqa(A6.4) $$
we find $|\D_{j}| \le A_{1}$ and $|n_{j_{0}}\cdot\D_{j}| \le
A_{2} :=A_{1}^{2}$ for all $j_{0}+1 \le j \le j_{0}+J_{1}$.

If ${\rm Span} \{\D_{j_{0}+1},\ldots,\D_{j_{0}+J_{1}}\}=D$
then by Lemma \secc(LA9) one has $|n_{j_{0}}| \le C(D) A_{1}^{D+1}$.
Then, for this relation to be not in contradiction with
$|n_{j_{0}}| > C_{1}(D) k^{1/D}$, we must have $C_{1}(D) k^{1/D} <
C(D) A_{1}^{D+1}$, hence $k \le C_{2}(D) \, r^{\a(D)\,\d}$
for some constant $C_{2}(D)$.

If ${\rm Span} \{\D_{j_{0}+1},\ldots,\D_{j_{0}+J_{1}}\} \le D-1$ then
there exists a subspace $H_{1}$ with ${\rm dim}(H_{1})=D-1$ such that
$n_{j} \in n_{j_{0}}+ H_{1}$ for $j_{0}+1\le j\le j_{0}+J_{1}$.
Choose $j_{1}\le J_{1}/2$ such that $P_{H_{1}}n_{j_{1}} :=
n_{j_{1}}- n_{j_{0}}\in H_{1}$ satisfies $|P_{H_{1}}n_{j_{1}}| >
C(D-1) J_{1}^{1/(D-1)}$, and fix $J_{2}=J_{1}^{1/\a(D-1)}$.
Redefine $\D_{j}=n_{j}-n_{j_{1}}$ for $j\ge j_{1}+1$,
$A_{1}=CJ_{2}r^{\d}$ and $A_{2}=A_{1}^{2}$: by reasoning as in the
previous case we find again $|\D_{j}| \le A_{1}$
and $|n_{j_{1}}\cdot\D_{j}| \le A_{2}$
for all $j_{0}+1 \le j \le j_{0}+J_{1}$.

If ${\rm Span} \{\D_{j_{1}+1},\ldots,\D_{j_{1}+J_{2}}\}=D-1$
then by Lemma \secc(LA9) one has $|n_{j_{1}}| \le C(D-1) A_{1}^{D}$,
which implies $C_{1}(D-1) J_{1}^{1/D-1} < C(D) A_{1}^{D}$.
By using the new definition of $A_{1}$, we obtain $J_{1} \le
C_{2}(D-1) \, r^{\a(D-1)\,\d}$, hence $k \le C_{3}(D) \,
r^{\a(D-1)\,\a(D)\,\d}$ for some other constant $C_{3}(D)$.

If ${\rm Span} \{\D_{j_{1}+1},\ldots,\D_{j_{0}+J_{2}}\} \le D-2$ then
there exists a subspace $H_{2}$ with ${\rm dim}(H_{1})=D-1$ such that
$n_{j} \in n_{j_{1}}+ H_{2}$ for $j_{1}+1\le j\le j_{1}+J_{2}$.
Then we iterate the construction until either we find
$k \le C_{n+2}(D) \,r^{\a(D-1)\ldots \a(D-n) \,\d}$ for
some $n \le D-1$ and some constant $C_{n+2}(D)$ or we arrive
at a subspace $H_{D-1}$ with ${\rm dim}(H_{D-1})=1$.

In the last case the vectors $\D_{j_{D-2}+1},\ldots, \D_{j_{D-2}+
J_{D-1}}$, with $J_{D-1}=J_{D-2}^{1/\a(2)}$, are linearly dependent
by construction, so that they lie all on the same line.
Therefore, we can find at least $J_{D-1}/2$ of them, say the
first $J_{D-1}/2$, with decreasing distance from the origin.
If we set $n_{j_{D-2}+1}=a$, $n_{j_{D-2}+J_{D-1}/2}=b$,
and $n_{j_{D-2}+1} - n_{j_{D-2}+J_{D-1}}/2=c$, and sum over
$j_{D-2}+1 \le j \le j_{D-2}+J_{D-1}/2$ the inequalities
$$ \left| n_{j} - n_{j-1} \right| +
|n_{j}|^{2} - |n_{j-1}|^{2}
\le {\rm const.} \left| \Phi(n_{j})-\Phi(n_{j-1}) \right| \le
{\rm const.} \, C r^{\d} ,
\Eqa(A6.5) $$
we obtain
$$ |c| + |c|^{2} \le |c| + |a|^{2} - |b|^{2} \le
{\rm const.} \, C r^{\d} { J_{D-1} \over 2} ,
\Eqa(A6.6) $$
where $|c| \ge J_{D-1}/2$. Hence $J_{D-1} \le (Cr^{\d})^{2}$.

By collecting together all the bound above we find $k \le C_{D}(D)
\, r^{2 \a(D)\ldots \a(2) \, \d}$, so that, by defining $C'=C_{D}(D)$
and using that $\e/\d=\a(D)\ldots \a(2)=2^{D-1}D!(D+1)!$,
the assertion follows. \qed

\*
\0\lma(LA11) {\it There exist constants $\e'$, $\d'$, $C$ and $C'$ such
that the following holds. Given $n_{0} \in \ZZZ^{D}$ there exists a set
$\D \subset \ZZZ^{D}$, with $n_{0}\in\D$, such that
${\rm diam}(\D) < C' r^{\e'}$ and $|\Phi(x)-\Phi(y)| > C'r^{\d'}$
for all $x\in\D$ and $y\notin\D$.}
\*

\0{\bf Proof.} Cf. \cita{Bo4}, p. 399, which proves the assertion
with $\e'=\d+\e$ and $\d'=\d=\d(\e,D)$. \qed

\*
\0\lma(LA12) {\it Let $\D$ be as in Lemma \secc(LA11). There exists
a constant $C''$ such that one has $|\D| \le C'' r^{D(\e+\d)}$.}
\*

\0{\bf Proof.} The bound follows from Lemma \secc(LA11)
and from the fact that ${\rm diam}(\D) < C' r^{\e}$,
by using that the points in $\D$ are distinct lattice
points in $\RRR^{D}$. \qed


\rife{Bo1}{1}{
J. Bourgain,
{\it Construction of quasi-periodic solutions for Hamiltonian
perturbations of linear equations and applications to nonlinear PDE},
Internat. Math. Res. Notices
{\bf 1994}, no. 11, 475ff., approx. 21 pp. (electronic). }
\rife{Bo2}{2}{
J. Bourgain,
{\it Construction of periodic solutions of nonlinear
wave equations in higher dimension},
Geom. Funct. Anal.
{\bf 5} (1995), 629--639. }
\rife{Bo3}{3}{
J. Bourgain,
{\it Periodic solutions of nonlinear wave equations},
Harmonic analysis and partial differential equations
(Chicago, IL, 1996), 69--97, Chicago Lectures in Math.,
Univ. Chicago Press, Chicago, IL, 1999. }
\rife{Bo4}{4}{
J. Bourgain,
{\it Quasi-periodic solutions of Hamiltonian perturbations of 2D
linear Schr\"odinger equations},
Ann. of Math. (2)
{\bf 148} (1998), no. 2, 363--439. }
\rife{Bo5}{5}{
J. Bourgain,
{\it Green's function estimates for lattice Schr\"odinger
operators and applications},
Annals of Mathematics Studies 158,
Princeton University Press, Princeton, NJ, 2005. }
\rife{C}{6}{
B. Connes,
{\it Sur les coefficients des s\'eries trigonom\'etriques
convergentes sph\'eriquement},
C. R. Acad. Sci. Paris S\'er. A-B
{\bf 283} (1976), no. 4, Aii, A159--A161. } 
\rife{CW}{7}{
W. Craig, C.E. Wayne,
{\it Newton's method and periodic solutions of nonlinear
wave equations},
Comm. Pure Appl. Math.
{\bf 46} (1993), 1409--1498. }
\rife{E}{8}{
L.H. Eliasson,
{\it Absolutely convergent series expansions for quasi periodic motions},
Math. Phys. Electron. J.
{\bf 2} (1996), Paper 4, 33 pp. (electronic). }
\rife{EK}{9}{
L.H. Eliasson, S. Kuksin, 
{\it KAM for non-linear Schr\"odinger equation},
preprint, 2006. }
\rife{Ga}{10}{
G. Gallavotti,
{\it Twistless KAM tori},
Comm. Math. Phys.
{\bf 164} (1994), no. 1, 145--156. } 
\rife{GY}{11}{
J. Geng, J. You,
{\it A KAM theorem for Hamiltonian partial differential equations
in higher dimensional spaces},
Comm. Math. Phys.
{\bf 262} (2006), no. 2, 343 - 372. }
\rife{GM2}{12}{
G. Gentile, V. Mastropietro,
{\it Construction of periodic solutions of the nonlinear wave equation
with Dirichlet boundary conditions by the Linsdtedt series method},
J. Math. Pures Appl. (9)
{\bf 83} (2004),  no. 8, 1019--1065. }
\rife{GMP}{13}{
G. Gentile, V. Mastropietro, M. Procesi,
{\it  Periodic solutions for completely resonant nonlinear
wave equations
with Dirichlet boundary conditions},
Comm. Math. Phys. {\bf 256} (2005), no. 2, 437-490. }
\rife{GP}{14}{
G. Gentile, M. Procesi,
{\it Conservation of resonant periodic solutions for the
one-dimensional non linear Schr\"odinger equation},
Comm. Math. Phys. {\bf 262} (2006), no. 3, 533-553. }
\rife{HP}{15}{
F. Harary, E.M. Palmer,
{\it Graphical enumeration},
Academic Press, New York-London, 1973. }
\rife{Ka}{16}{
T. Kato,
{\it Perturbation theory for linear operators}, 
Springer-Verlag, Berlin-New York, 1976. }
\rife{K1}{17}{
S.B. Kuksin,
{\it Nearly integrable infinite-dimensional Hamiltonian systems},
Lecture Notes in Mathematics 1556, Springer, Berlin, 1994. }
\rife{KP}{18}{
S.B. Kuksin, J.P\"oschel,
{\it Invariant Cantor manifolds of quasi-periodic oscillations
for a nonlinear Schr\"odinger equation},
Ann. of Math. (2)
{\bf 143} (1996), no. 1, 149--179. }
\rife{LS}{19}{
B.V. Lidski\u\i, E.I. Shul$^{\prime}$man,
{\it Periodic solutions of the equation $u_{tt}-u_{xx}+u^{3}=0$},
Funct. Anal. Appl.
{\bf 22} (1988), no. 4, 332--333. }
\rife{W}{20}{
C.E. Wayne,
{\it Periodic and quasi-periodic solutions of nonlinear wave
equations via KAM theory},
Comm. Math. Phys.
{\bf 127} (1990), no. 3, 479--528. }
\rife{Wh}{21}{
H. Whitney,
{\it Analytic extensions of differential
functions defined in closed sets},
Trans. Amer. Math. Soc.
{\bf 36} (1934), no. 1, 63--89. }
\rife{Y}{22}{
X. Yuan,
{\it A KAM theorem with applications to
partial differential equations of higher dimension},
Preprint, 2006. }

\biblio

\bye